\documentclass[10pt,a4]{article}

\usepackage{amscd}
\usepackage{amssymb}
\usepackage{amsmath}
\usepackage{amsthm}

\newtheorem{thm}{Theorem}[section]
\newtheorem{lemma}[thm]{Lemma}
\newtheorem{corollary}[thm]{Corollary}
\newtheorem{prop}[thm]{Proposition}

\theoremstyle{definition}

\newtheorem{rem}[thm]{Remark}

\newtheorem{defn}[thm]{Definition}

\makeatletter

\newcommand{\isom}{\overset{\sim}{\rightarrow}}

\def\le{{\leqq}}
\def\ge{{\geqq}}

\title{Classification of two dimensional split trianguline 
representations of $p$-adic fields.}

\author{Kentaro Nakamura}

\begin{document}
\maketitle
\pagestyle{plain}
\footnote{2008 Mathematical Subject Classification 11F80 (primary), 11F85, 11S25 (secondary).
Keywords: $p$-adic Hodge theory, trianguline representations, B-pairs.}
\begin{abstract}
The aim of this article is to classify two dimensional 
split trianguline representations of $p$-adic fields.
This is a generalization of a result of Colmez who classified
two dimensional split trianguline representations of $\mathrm{Gal}(\bar{\mathbb{Q}}_p/\mathbb{Q}_p)$ for $p\not= 2$ by using $(\varphi,\Gamma)$-modules over Robba ring.
In this article, for any prime $p$ and for any $p$-adic field $K$, we classify two dimensional split trianguline representations of $\mathrm{Gal}(\bar{K}/K)$ by using $B$-pairs defined by Berger.

\end{abstract}

\tableofcontents
\section{Introduction.}

\subsection{Background.}
Let $p$ be a prime number. In this article, we study two dimensional 
trianguline representations of any $p$-adic field $K$, i.e. a finite extension of $\mathbb{Q}_p$. Trianguline representation is a class of $p$-adic 
representations of $G_K:=\mathrm{Gal}(\bar{K}/K)$, which is turning out to be 
an important notion. 
Trianguline representation was defined by Colmez and is defined by using 
$(\varphi,\Gamma_K)$-modules over Robba ring $B^{\dagger}_{\mathrm{rig},K}$  which is noncanonically isomorphic 
to a ring of Laurent power series which converge in some annulus and is equipped with Frobenius $\varphi$ action and with $\Gamma_K:=\mathrm{Gal}
(K_{\infty}/K)$ actions. Here $K_{\infty}:=K(\zeta_{p^{\infty}})$ is the extension of $K$ obtained by adjoing $p^n$-th roots of unity $\zeta_{p^{n}}$ for every 
$n\in\mathbb{N}$. $(\varphi,\Gamma_K)$-modules over $B^{\dagger}_{\mathrm{rig},K}$ are defined as finite free $B^{\dagger}_{\mathrm{rig},K}$-modules with semi-linear $\varphi$ and $\Gamma_K$ actions. 
By the work of Kedlaya, the notion of slopes of $\varphi$-modules over Robba ring is very important and we say that a $(\varphi,\Gamma_K)$-module over $B^{\dagger}_{\mathrm{rig},K}$ is $\acute{\mathrm{e}}$tale if it is pure 
of slope zero as a $\varphi$-module. 
By the works of Fontaine, Cherbonnier-Colmez and Kedlaya, the category of $p$-adic representations of $G_K$ is equivalent to the category of $\acute{\mathrm{e}}$tale $(\varphi,\Gamma_K)$-modules over $B^{\dagger}_{\mathrm{rig},K}$. This equivalence enables us to see the category of $p$-adic representaion of $G_K$ as a full subcategory of the category of $(\varphi, \Gamma_K)$-modules over $B^{\dagger}_{\mathrm{rig},K}$ (without slope conditions). We say a $p$-adic 
representation $V$ of $G_K$ is split trianguline if $D_{\mathrm{rig}}(V)$, the $(\varphi, \Gamma_K)$-module corresponding to $V$, is a succesive extension of rank one 
objects in the catetgory of $(\varphi,\Gamma_K)$-modules over $B^{\dagger}_{\mathrm{rig},K}$. And we say $V$ is 
trianguline if $V$ is split trianguline after making a finite 
extension of coefficients. 
Trianguline representations 
can be seen as generalizations of ordinary 
representations. But many interesting irreducible representations 
can be trianguline. For example, all semi-stable representations  are trianguline. Moreover, there are many interesting trianguline representations 
which are not de Rham. For example, when $K=\mathbb{Q}_p$, Kisin showed that 
two dimensional $p$-adic representations of $G_{\mathbb{Q}_p}$attached to finite slope overconvergent modular forms are trianguline, this fact was a key for his $p$-adic Hodge theoretic approach to the study of Coleman-Mazur eigencurves ($\cite{Ki}$). Recently Bella$\ddot{\i}$che-Chenevier generalized Kisin's method and studied 
higher dimensional trianguline representations of $G_{\mathbb{Q}_p}$ to study 
higher dimensional eigenvarieties ($\cite{Bel-Ch}$). 
Another example is Colmez's $p$-adic Langlands correspondence for $\mathrm{GL}_2(\mathbb{Q}_p)$. In his theory of $p$-adic Langlands correspondence, trianguline representations are at the heart. When $K=\mathbb{Q}_p$ and $p\not= 2$, he classified two dimensional split trianguline representations of $G_{\mathbb{Q}_p}$ and, based on this classification, he proved that the sets of points corresponding to trianguline representations are Zariski dense in deformation spaces of two dimensional $p$-adic representaions of $G_{\mathbb{Q}_p}$. By combining some constructions of Colmez and of 
Berger-Breuil, we get the $p$-adic Langlands correspondence for trianguline 
representations ($\cite{Be-Br}$, $\cite{Co04}$, $\cite{Co07b}$). The correspondence for trianguline representations and Zariski density of trianguline points played essential roles in his construction of $p$-adic Langlands correspondence for $\mathrm{GL}_2(\mathbb{Q}_p)$ ($\cite{Co08}$). 

The main purpose of this article is to classify completely two dimensional 
split trianguline representations of $G_K$ for any finite extension $K$ of $\mathbb{Q}_p$ for any prime $p$ (we don't need to assume that $p\not= 2$). We determine the parameter space of all split trianguline representations. 
And we also determine the parameter space of potentially cristalline 
or potentially semi-stable split trianguline representations, then we explicitly describe the filtered $(\varphi,N,G_K)$-modules associated to tnem. 
Currently, the only interesting examples of tringuline representations are in the case $K=\mathbb{Q}_p$. The author of this article hopes that this article will be useful to give many interesting examples of trianguline representations in general $K\not =\mathbb{Q}_p$ case. For example, he wants to know whether $p$-adic representaions attached to finite slope overconvergent Hilbert 
modular forms are trianguline or not.
In this article, he could not attack the problem about Zariski density of trianguline representations. He doesn't know whether the set of 
trianguline points is Zariski dense or not in general $K\not =\mathbb{Q}_p$ case. He wants to study these problems in future works.

\subsection{Contents of this article.}
In section 1, we define trianguline representations 
and split trianguline representations 
by using $B$-pairs instead of $(\varphi,\Gamma_K)$-modules 
over $B^{\dagger}_{\mathrm{rig},K}$ for some technical reasons. 
$B$-pair was defined by Berger ($\cite{Be07}$). We write
$B_e:=B^{\varphi=1}_{\mathrm{cris}}$. 
An $E$-$B$-pair (the $E$-coefficient version of $B$-pair, here $E$ is a finite 
extension of $K$ which contains Galois closure of $K$) is a pair $W:=(W_e,W^+_{\mathrm{dR}})$ where $W_e$ is a finite free 
$B_e\otimes_{\mathbb{Q}_p}E$-module with a continuous semi-linear $G_K$-action such that $W^+_{\mathrm{dR}}\subseteq W_{\mathrm{dR}}:=B_{\mathrm{dR}}\otimes_{B_e}W_e$ is a 
$G_K$-stable $B^+_{\mathrm{dR}}\otimes_{\mathbb{Q}_p}E$-lattice of $W_{\mathrm{dR}}$. In $\cite[\mathrm{Theorem}\,2.2.7]{Be07}$, Berger established an equivalence between the category of $E$-$B$-pairs and 
the category of $E$-($\varphi,\Gamma_K$)-modules over $B^{\dagger}_{\mathrm{rig},K}$ which are also the $E$-coefficient version of ($\varphi,\Gamma_K$)-modules.
The category of $E$ representations of $G_K$ ($E$-coefficient version of 
$p$-adic representations) is embedded in the 
category of $E$-$B$-pairs by $V\mapsto W(V):=(B_e\otimes_{\mathbb{Q}_p}V, B^+_{\mathrm{dR}}\otimes_{\mathbb{Q}_p}V)$.
So, for defining trianguline representations, we can use both $B$-pairs and 
$(\varphi,\Gamma_K)$-modules. In this article, we choose to use $B$-pairs for 
some technical reasons. Then, we say that an $E$-$B$-pair $W$ is split trianguline if $W$ is a succesive extension of rank one $E$-$B$-pairs. We say that $E$-representation $V$ is split trianguline if $W(V)$ is split trianguline. We say that an $E$-$B$-pair $W$ (resp. an $E$-representation $V$) is trianguline if $W\otimes_{E}E'$ (resp. $V\otimes_E E'$) is split trianguline for a finite extension $E'$ of $E$. The main purpose of this article 
is to classify completely two dimensional split trianguline $E$-representations  of $G_K$. The classification is done by the following steps,
\begin{itemize}
\item[Step 1]: Classification of rank one $E$-$B$-pairs,
\item[Step 2]: Calculation of the dimensions of extension groups $\mathrm{dim}_E\mathrm{Ext}^1(W_2, W_1)$ for any rank one $E$-$B$-pairs $W_1$, $W_2$,
\item[Step 3]: Determination of the conditions for $W$ to be $\acute{\mathrm{e}}$tale, here $W$ is an extension of $W_1$ by $W_2$ as in Step2,
\item[Step 4]: Classification of de Rham split trianguline 
$E$-representations.
\end{itemize}
The step 1 is done in section 1, the step 2 is done in section 2, 
the step 3 is done in section 3, the step 4 is done in section 4.

From now on, we explain these steps more precisely. 

For step 1, let $\delta:K^{\times}\rightarrow E^{\times}$ be a continous character 
with respect to $p$-adic topology on both sides. Then we can define a rank one $E$-$B$-pair $W(\delta)$ as follows. Let $\pi_K$ be a uniformizer of $K$. 
We decompose $\delta$ into $\delta:=\delta_0\delta_1$ such that $\delta_0|_{\mathcal{O}_K^{\times}}=\delta|_{\mathcal{O}_K^{\times}}$, $\delta_0(\pi_K):=1$, $\delta_1|_{\mathcal{O}_K^{\times}}$ is trivial character, $\delta_1(\pi_K):=\delta(\pi_K)$. Then $\delta_0:K^{\times}\rightarrow \mathcal{O}_E^{\times}$ is a unitary character, so by local class field 
theory we get a character $\tilde{\delta}_0:G_K^{\mathrm{ab}}\rightarrow \mathcal{O}_E^{\times}$ such that $\delta_0=\tilde{\delta}_0\circ \mathrm{rec}_K$, here $\mathrm{rec}_K:K^{\times}\rightarrow G_K^{\mathrm{ab}}$ is the reciprocity 
map such that $\pi_K$ is mapped to a lifting of the inverse of 
$q$-th power Frobenius (here, $q:=p^f$, $f:=[K_0:\mathbb{Q}_p]$, $K_0$ is the 
maximal unramified extension of $\mathbb{Q}_p$ in $K$). 
 For $\delta_1$, we define a rank one $E$-$B$-pair $W(\delta_1)$ such that $D(\delta_1)$, the $(\varphi,\Gamma_K)$-module corresponding to $W(\delta_1)$, is defined by $D(\delta_1):=B^{\dagger}_{\mathrm{rig},K}\otimes_{\mathbb{Q}_p}Ee_{\delta_1}$, $\varphi^f(e_{\delta_1})
:=\delta_1(\pi_K)e_{\delta_1}$, $\gamma(e_{\delta_1}):=e_{\delta_1}$ for any 
$\gamma\in \Gamma_K$. We define $W(\delta):=
W(\delta_0)\otimes W(\delta_1)$. We can show that $W(\delta)$ does not depend 
on the choice of $\pi_K$ (Remark $\ref{-4}$). Then the main theorem of section 1 is the 
following.
\begin{thm}$(\mathrm{Theorem}\, \ref{15})$
Let $W$ be a rank one $E$-$B$-pair of $G_K$. Then there exists unique continuous character $\delta:K^{\times}\rightarrow E^{\times}$ such that $W\isom W(\delta)$.
\end{thm}
\begin{rem}
We can see this theorem as a natural generalization of one to one correspondence $\{\delta:K^{\times}\rightarrow \mathcal{O}^{\times}_E: $continuous character$\}\isom \{\tilde{\delta}:G^{\mathrm{ab}}_K\rightarrow \mathcal{O}^{\times}_E:$ continuous character$\}$, which is induced by local class field theory.
This theorem is also a generalization of $\cite[\mathrm{Proposition}\,3.1]{Co07a}$.
\end{rem}

For step 2, from this theorem, it suffices to calculate $\mathrm{dim}_E\mathrm{Ext}^1(W(\delta_2), W(\delta_1))$ for any continuous characters $\delta_1,\delta_2:K^{\times}\rightarrow E^{\times}$. In section 2, for this purpose we define Galois cohomology of $E$-$B$-pairs. Let 
$W:=(W_e,W^+_{\mathrm{dR}})$ be an $E$-$B$-pair. Put $W_{\mathrm{dR}}:=
B_{\mathrm{dR}}\otimes_{B_e}W_e$. Then we defne the Galois cohomology $\mathrm{H}^{*}(G_K,W)$ of $W$ as the Galois cohomology of the complex $W_e\oplus W^+_{\mathrm{dR}}\rightarrow W_{\mathrm{dR}}:(x,y)\mapsto x-y$, here $W_e\oplus W^+_{\mathrm{dR}}$ sits in zero-th part of this complex. This cohomology groups are finite dimensional $E$-vector spaces. We can show in the usual way that 
there is a natural isomorphism $\mathrm{Ext}^1(B_E, W)\isom \mathrm{H}^1(G_K, W)$, here $B_E:=(W_e\otimes_{\mathbb{Q}_p}E,W^+_{\mathrm{dR}}\otimes_{\mathbb{Q}_p}E)$ is the trivial $E$-$B$-pair. By Bloch-Kato's fundamental 
short exact sequence $0\rightarrow \mathbb{Q}_p\rightarrow 
B_e\oplus B^+_{\mathrm{dR}}\rightarrow B_{\mathrm{dR}}\rightarrow 0$, for any 
$E$-representation $V$, we have a natural isomorphism $\mathrm{H}^*(G_K,V)\isom \mathrm{H}^*(G_K, W(V))$. So we can see this cohomology as a natural generalization of Galois cohomology of $p$-adic representations. As in the classical case, this cohomology satisfies Euler-Poincar$\acute{\mathrm{e}}$ characteristic formula and Tate local duality theorem. For this, we review the results of Liu concerning these theorems for $(\varphi,\Gamma_K)$-modules over $B^{\dagger}_{\mathrm{rig},K}$ ($\cite{Li}$). By using his theorems, we can calculate all the extension groups that we want. For any embedding $\sigma:K\hookrightarrow E$ and $k\in\mathbb{Z}$, we define $\sigma(x)^k:K^{\times}\rightarrow E^{\times}:y\mapsto \sigma(y)^k$. We define $N_{K/\mathbb{Q}_p}:K^{\times}\rightarrow \mathbb{Q}_p^{\times}:y\mapsto\prod_{\sigma:K\hookrightarrow E}\sigma(y)$, $|-|:\mathbb{Q}^{\times}_p\rightarrow E:p\mapsto\frac{1}{p},\,a\mapsto 1$ for any $a\in
\mathbb{Z}_p^{\times}$, $|N_{K/\mathbb{Q}_p}(x)|:=|-|\circ N_{K/\mathbb{Q}_p}$. Then the main result of section 2 is as follows.
\begin{thm}$(\mathrm{Teorem}\,\ref{29})$
Let $\delta_1,\delta_2:K^{\times}\rightarrow E^{\times}$ be continuous 
characters. 

Then $\mathrm{dim}_E\mathrm{Ext}^1(W(\delta_2),W(\delta_1))$ is 
equal to 
\begin{itemize}
\item[$\mathrm{(1)}$] $[K:\mathbb{Q}_p]+1$ if $\delta_1/\delta_2= \prod_{\sigma:K\hookrightarrow E}\sigma(x)^{k_{\sigma}}$ such that $k_{\sigma}\in\mathbb{Z}_{\le 0}$ for any $\sigma$,
\item[$\mathrm{(2)}$]$[K:\mathbb{Q}_p]+1$ if $\delta_1/\delta_2=|N_{K/\mathbb{Q}_p}(x)|\prod_{\sigma:K\hookrightarrow E}\sigma(x)^{k_{\sigma}}$ such that $k_{\sigma}\in\mathbb{Z}_{\ge 1}$ for any $\sigma$,
\item[$\mathrm{(3)}$]$[K:\mathbb{Q}_p]$ otherwise.
\end{itemize}
\end{thm}
\begin{rem}
This theorem is a generalization of dimension fromula 
of Galois cohomology of one dimensional $p$-adic representations. 
Also this is a generalization of $\cite[\mathrm{Theorem}\,2.9]{Co07a}$.
\end{rem}

From this theorem we can determine all the split trianguline $E$-$B$-pairs 
of rank two. Let $W(s)$ be a split trianguline $E$-$B$-pair which is 
an extension corresponding 
to $s\in\mathbb{P}_E(\mathrm{Ext}^1(W(\delta_2),\allowbreak W(\delta_1)))$, here for any 
finite dimensional $E$-vector space $M$ we put $\mathbb{P}_E(M):=\{[v]|v\in M-\{0\}, [v]=[v'] \iff v'=av$ for some $a\in E^{\times}\}$. Then the isomorphism class of $W(s)$ as $E$-$B$-pair depends only on $s$. 

For step 3, we must determine the conditions on ($\delta_1$, $\delta_2$) and $s$ for $W(s)$ to be $\acute{\mathrm{e}}$tale, i.e. to be $W(s)\isom W(V(s))$ for an  $E$-representation $V(s)$. In section 3, we determine all the conditions 
by 
using Kedlaya's slope filtaration theorem of $\varphi$-modules over Robba ring. The idea is essentially same as Colmez's one when $K=\mathbb{Q}_p$ and $p\not=2$ (see the 
proof of $\cite[\mathrm{Proposition}\,4.7]{Co07a}$), but 
we have to deal with all the additional complications which come from 
working with $K\not =\mathbb{Q}_p$. In fact, in general $K\not= \mathbb{Q}_p$ case, the parameter space of split 
trianguline representations is more complicated than that of $K=\mathbb{Q}_p$ case. For any two continuous characters $\delta_1,\delta_2:K^{\times}\rightarrow E^{\times}$, we put $S(\delta_1,\delta_2):=\mathbb{P}_E(\mathrm{Ext}^1(W(\delta_2), W(\delta_1)))$. We put $S^+:=\{(\delta_1,\delta_2)| \delta_1,\delta_2:K^{\times}\rightarrow E^{\times}$ continuous characters such that $\mathrm{val}_p(\delta_1(\pi_K))+\mathrm{val}_p(\delta_2(\pi_K))=0, \mathrm{val}_p(\delta_1(\pi_K))\ge 0\}$ (here $\mathrm{val}_p$ is a valuation of $E$ such that $\mathrm{val}_p(p):=1$). For any $(\delta_1,\delta_2)\in S^+$, in section 3 we will explicitly define a certain subspace ${S'}^{non-\acute{\mathrm{e}}t}(\delta_1,\delta_2)\subseteq 
S(\delta_1,\delta_2)$ which corresponds to non $\acute{\mathrm{e}}$tale split 
trianguline $E$-$B$-pairs. All these spaces $S^+$ and ${S'}^{non-\acute{\mathrm{e}}t}(\delta_1,\delta_2)$ naturally appear when we consider the slope zero conditions by using Kedlaya's slope filtration theorem. Then our main result of section 3 is as follows.
\begin{thm}$(\mathrm{Lemma}\,\ref{33},\,\mathrm{Theorem}\, \ref{34})$
Let $\delta_1,\delta_2:K^{\times}\rightarrow E^{\times}$ be continuous 
characters. Let $W(s)$ be the split trianguline $E$-$B$-pair corresponding 
to $s\in S(\delta_1,\delta_2)$.
\begin{itemize}
\item[$\mathrm{(1)}$]If $W(s)$ is $\acute{\mathrm{e}}$tale, then 
$(\delta_1,\delta_2)\in S^+$.
\item[$\mathrm{(2)}$]The following conditions are equivalent.
\begin{itemize}
\item[$\mathrm{(i)}$]$W(s)$ is $\acute{\mathrm{e}}$tale, i.e. $W(s)\isom W(V(s))$ for an $E$-representation $V(s)$.
\item[$\mathrm{(ii)}$]$s\notin {S'}^{non-\acute{\mathrm{e}}t}(\delta_1,\delta_2)$.
\end{itemize}
\end{itemize}
\end{thm}
\begin{rem}
This theorem is a generalization of $\cite[\mathrm{Proposition}\,4.7]{Co07a}$.
 The space $\mathcal{S}_+^{\mathrm{ncl}}$ in his paper $\cite[0.2]{Co07a}$ 
corresponds to 
$\sqcup_{(\delta_1,\delta_2)\in S^+}S^{non-\acute{\mathrm{e}}t}(\delta_1,\delta_2)$ in this article. Moreover we can determine the conditions when we have $V(s)\isom V(s')$ for distinct parameters 
$s\in S(\delta_1,\delta_2)\setminus {S'}^{non-\acute{\mathrm{e}}t}(\delta_1,\delta_2)$, 
$s'\in S(\delta'_1,\delta'_2)\setminus {S'}^{non-\acute{\mathrm{e}}t}(\delta'_1,\delta'_2)$
 under certain conditions (Theorem $\ref{36}$).
\end{rem}

For step 4, we have to determine the conditions on ($\delta_1$,$\delta_2)\in S^+$ and $s\in S(\delta_1,\delta_2)\setminus {S'}^{non-\acute{\mathrm{e}}t}(\delta_1,\delta_2)$ for $V(s)$ to be potentially semi-stable or potentially cristalline. 
In $p$-adic representations case, Bloch-Kato finite cohomology is useful 
for this kind of problems. We define Bloch-Kato cohomology 
for $B$-pair $W$ as follows (Definition $\ref{d}$). By the definition of $\mathrm{H}^*(G_K, W)$, we have natural maps $\mathrm{H}^*(G_K, W)\rightarrow \mathrm{H}^*(G_K, W_e)
\rightarrow \mathrm{H}^*(G_K, B_{\mathrm{cris}}\otimes_{B_e}W_e)\rightarrow 
\mathrm{H}^*(G_K, B_{\mathrm{dR}}\otimes_{Bt_e}W_e)$. As in the classical case,  we define $\mathrm{H}^1_{?}(G_K,W):=\mathrm{Ker}(\mathrm{H}^1(G_K,W)\rightarrow \mathrm{H}^1(G_K,B_{*}\otimes_{B_e}W_e))$, here when $?=e$ (resp. $?=f$, resp.  $?=g$) then $*=e$ (resp. $*=\mathrm{cris}$,  resp. $*=\mathrm{dR}$). 
By calculating these, in section 4 we will explicitly define the parameter spaces $S^{\acute{\mathrm{e}}t}_{\mathrm{cris}}(\delta_1,\delta_2)$, $S_{\mathrm{st}}(\delta_1,\delta_2) \subseteq S(\delta_1,\delta_2)$ for any $(\delta_1,\delta_2)\in S^+$ such that $\delta_i=\prod_{\sigma:K\hookrightarrow E}\sigma(x)^{k_{i,\sigma}}\tilde{\delta}_i$ such that $k_{i,\sigma}\in\mathbb{Z}$ and $\tilde{\delta}_i$ are locally constant characters for i=1,2.. Then the conditions for $V(s)$ to be potentially semi-stable are as follows.
\begin{thm}$(\mathrm{Lemma}\,\ref{38}, \mathrm{Proposition}\,\ref{41})$
Let $s\in S(\delta_1,\delta_2)\setminus {S'}^{non-\acute{\mathrm{e}}t}(\delta_1,\delta_2)$. Let $V(s)$ be the split trianguline $E$-representation corresponding to $s$. 
\begin{itemize}
\item[$\mathrm{(1)}$]If $V(s)$ is potentially semi-stable, then 
$\delta_i=\prod_{\sigma:K\hookrightarrow E}\sigma(x)^{k_{i,\sigma}}\tilde{\delta}_i$ such that $k_{i,\sigma}\in\mathbb{Z}$ for any $\sigma$ and $\tilde{\delta}_i$ are locally constant characters for i=1,2.
\item[$\mathrm{(2)}$]If $(\delta_1,\delta_2)\in S^+$ satisfies the condition in 
$\mathrm{(1)}$. Then the following conditions are equivalent.
\begin{itemize}
\item[$\mathrm{(i)}$]$V(s)$ is potentially criatalline
\item[$\mathrm{(ii)}$]$s\in S^{\acute{\mathrm{e}}t}_{\mathrm{cris}}(\delta_1,\delta_2)$.
\end{itemize}
\item[$\mathrm{(3)}$]If $(\delta_1,\delta_2)\in S^+$ satisfies the condition in 
$\mathrm{(1)}$. Then the following conditions are equivalent.
\begin{itemize}
\item[$\mathrm{(i)'}$]$V(s)$ is potentially semi-stable and not potentially cristalline.
\item[$\mathrm{(ii)'}$]$s\in S_{\mathrm{st}}(\delta_1,\delta_2)$.
\end{itemize}
\end{itemize}
\end{thm}

\begin{rem}
These parameter spaces $S^{\acute{\mathrm{e}}t}_{\mathrm{cris}}(\delta_1,\delta_2)$, $S_{\mathrm{st}}(\delta_1,\delta_2)$ are generalizations of $\mathcal{S}_+^{\mathrm{cris}}$ or 
$\mathcal{S}_+^{\mathrm{st}}$ defined in $\cite[0.2]{Co07a}$. 
We can also see these spaces as the parameter spaces of weakly admissible filtrations 
of a $(\varphi,N,G_K)$-module corresponding to $V(s)$. 
Moreover we can explicitly calculate the filtered $(\varphi,N,G_K)$-modules associated to $V(s)$ as above ($\mathrm{Theorem}\ref{42}$, $\mathrm{Theorem}\ref{45}$).
\end{rem}

By this theorem, we complete the classification of two dimensional 
split trianguline $E$-representations. These are main contents of this article.

In the appendix, we study a relation between two dimensional potentially semi-stable trianguline representations and local Langlands correspondence for $\mathrm{GL}_2(K)$.
Let $V$ be a two dimensional potentially semi-stable $E$-representation of $G_K$. Fontaine defined a two dimensional Weil-Deligne representation $\bar{D}_{\mathrm{pst}}(V):=D_{\mathrm{pst}}(V)\otimes_{K_0^{\mathrm{un}}\otimes_{\mathbb{Q}_p}E}\bar{K}$ of $K$ from the filtered $(\varphi, N,G_K)$-module 
$D_{\mathrm{pst}}(V):=\cup_{K\subseteq L,\mathrm{finite}}(B_{\mathrm{st}}\otimes_{\mathbb{Q}_p}V)^{G_L}$. By local Langlands correspondence for $\mathrm{GL}_2(K)$, we can attach an irreducible smooth admissible representation $\pi(\bar{D}_{\mathrm{pst}}(V)^{ss})$ of $\mathrm{GL}_2(K)$ (here $\bar{D}_{\mathrm{pst}}(V)^{ss}$ is the Frobenius semi-simplification of $\bar{D}_{\mathrm{pst}}(V)$).
Irreducible smooth admissible representations of $\mathrm{GL}_2(K)$ are classified 
into supercuspidal ones and non supercuspidal ones (i.e. one dimensional representations, principal series or special series).
Then the main result of the appendix is the following.
\begin{thm}$(\mathrm{Theorem}\, \ref{19})$
Let $V$ be a potentially semi-stable $E$-representation.
Then the following conditions are equivalent.
\begin{itemize}
\item[$\mathrm{(1)}$]$V$ is trianguline, i.e. $V\otimes_{E}E'$ is split trianguline for some finite extension $E'$ of $E$.
\item[$\mathrm{(2)}$]$\pi(\bar{D}_{\mathrm{pst}}(V)^{ss})$ is non supercuspidal.\end{itemize}
\end{thm}

\subsection*{Notation.}
Let $p$ be a prime number. $K$ is a finite extension of $\mathbb{Q}_p$. 
$\bar{K}$ is a fixed algebraic closure of $K$.
$K_0$ is the maximal unramified extension of $\mathbb{Q}_p$ in $K$.
$K^{\mathrm{nor}}$ is the Galois closure of $K$ in $\bar{K}$.
$G_K:=\mathrm{Gal}(\bar{K}/K)$ is the absolute Galois group of $K$ equipped with profinite topology. 
$\mathcal{O}_K$ is the integer ring of $K$. $\pi_K\in\mathcal{O}_K$ is a uniformizer of $K$. $k:=\mathcal{O}_K/\pi_K\mathcal{O}_K$ is the residue field of $K$. $q=p^f:=\sharp k$ is the order of $k$.
$K_{\infty}:=K(\zeta_{p^{\infty}})$ is the extension of $K$ obtained by adjoining $p^n$-th 
roots of unity $\zeta_{p^n}$ for every $n\in\mathbb{N}$. $K'_0$ is the maximal unramified extension of $\mathbb{Q}_p$ in $K_{\infty}$, 
so $K_0\subset K'_0$.
$H_K:=\mathrm{Gal}(\bar{K}/K_{\infty})$, $\Gamma_K:=G_K/H_K=\mathrm{Gal}(K_{\infty}/K)$.
$\chi:G_K\rightarrow \mathbb{Z}_p^{\times}$ is the $p$-adic cyclotomic character which factors through the inclusion $\Gamma_K\hookrightarrow \mathbb{Z}_p^{\times}$ (i.e. $g(\zeta_{p^n})=\zeta_{p^n}^{\chi(g)}$ for any $p^n$-th roots 
of unity $\zeta_{p^n}$ and for any $g\in G_K$).
$\mathbb{C}_p:=\widehat{\bar{K}}$ is the $p$-adic completion of $\bar{K}$, which is 
an algebraically closed $p$-adically complete field. $\mathcal{O}_{\mathbb{C}_p}$ is  its integer ring.   
$E$ is a finite extension of $\mathbb{Q}_p$ such that $K^{\mathrm{nor}}\subset E$. In this paper, we will write $E$ as coefficient of representations.
$\chi_{\mathrm{LT}}:G_K\rightarrow \mathcal{O}_K^{\times}\hookrightarrow E^{\times}$ is the Lubin-Tate character associated to the uniformizer $\pi_K$. 
$\mathrm{rec}_K:K^{\times}\rightarrow G_K^{\mathrm{ab}}$ is the reciprocity map of local class field theory such that $\mathrm{rec}_K(\pi_K)$ is a lifting of the inverse of $q$-th power Frobenius of $k$, then $\chi_{\mathrm{LT}}\circ\mathrm{rec}_K:K^{\times}\rightarrow \mathcal{O}^{\times}_{K}$ satisfies $\chi_{\mathrm{LT}}\circ \mathrm{rec}_K(\pi_K)=1$ and 
 $\chi_{\mathrm{LT}}\circ\mathrm{rec}_K|_{\mathcal{O}_K^{\times}}=id_{\mathcal{O}_K^{\times}}$. 

\subsection*{Acknowledgements.}
This paper is originally the author's doctor thesis 
at University of Tokyo. The author would like to 
express his sincere gratitude to his thesis advisor 
Prof. Atsushi Shiho for valuable discussion and reading the 
manuscript carefully.
This paper could never have existed without his constant 
encouragement. And the author also thank Laurent 
Berger for many useful suggestions and for reading the revised 
version carefully, especially for the 
suggestion about B$\acute{\mathrm{e}}$zout property of $B_e\otimes_{\mathbb{Q}_p}E$, which makes many arguments in 
the first version of this article easier.

\section{$B$-pairs and $(\varphi,\Gamma_K)$-modules.}
In this section, we recall the definitions and some properties of $B$-pairs and $(\varphi, \Gamma_K)$-modules. In particular, we recall the results of Berger concerning the equivalence between 
the category of $B$-pairs and the category of $(\varphi,\Gamma_K)$-modules over Robba ring.
\subsection{Review of $p$-adic period rings and the definition of $B$-pair.}
We begin this section by recalling some $p$-adic period rings ($\cite{Fo94a},\cite{Fo94b},\cite{Be02}$) and by recalling 
the definition of $B$-pair ($\cite{Be07}$).
First let $\widetilde{\mathbb{E}}^+:=\varprojlim_n\mathcal{O}_{\mathbb{C}_p}=
\varprojlim_n\mathcal{O}_{\mathbb{C}_p}/p\mathcal{O}_{\mathbb{C}_p}$, where the limits are taken with respect to $p$-th power maps. It is known that $\widetilde{\mathbb{E}}^+$ is a complete valuation ring of characteristic $p$ whose valuation
is defined by $\mathrm{val}(x):=\mathrm{val}_{p}(x^{(0)})$  (here $x=(x^{(n)})\in \varprojlim_n\mathcal{O}_{\mathbb{C}_p}$ and $\mathrm{val}_p$ is the valuation on $\mathbb{C}_p$ such that $\mathrm{val}_p(p)=1$).
In this paper, we fix a system of $p^n$-th roots of unity $\{\varepsilon^{(n)}\}_{n\ge 0}$ 
such that $\varepsilon^{(0)}=1$, $(\varepsilon^{(n+1) })^p=\varepsilon^{(n)}$, 
$\varepsilon^{(1)}\not= 1$. Then $\varepsilon:=(\varepsilon^{(n)})$ is an element of $\widetilde{\mathbb{E}}^+$ such that $\mathrm{val}(\varepsilon -1)=p/(p-1)$. $\widetilde{\mathbb{E}}:=\widetilde{\mathbb{E}}^+[1/(\varepsilon -1)]$ is the fraction field of 
$\widetilde{\mathbb{E}}^+$, which is known to be an algebraically closed complete valuation field of characteristic $p$ containing the subfield $\mathbb{F}_p((\varepsilon -1))$. $G_K$ acts on these rings in natural way. We put  $\widetilde{\mathbb{A}}^+:=W(\widetilde{\mathbb{E}}^+)$, $\widetilde{\mathbb{A}}:=W(\widetilde{\mathbb{E}})$, where, for a ring $R$,
$W(R)$ is the Witt ring of $R$. We put $\widetilde{\mathbb{B}}^+:=\widetilde{\mathbb{A}}^+[\frac{1}{p}]$, 
$\widetilde{\mathbb{B}}:=\widetilde{\mathbb{A}}[\frac{1}{p}]$. These rings also have natural continuous $G_K$-actions and Frobenius 
actions $\varphi$, here the topology of these rings are defined by the 
$p$-adic topology on $\widetilde{\mathbb{A}}^+$ and $\widetilde{\mathbb{A}}^+$. Then we have a continuous $G_K$-equivariant surjection $\theta:\widetilde{\mathbb{A}}^+
\rightarrow \mathcal{O}_{\mathbb{C}_p}:\sum_{k=0}^{\infty}p^k[x_k]\mapsto \sum_{k=0}^{\infty} p^k x_k^{(0)}$, where $[\,\, ]:\widetilde{\mathbb{E}}^+\rightarrow \widetilde{\mathbb{A}}^+$ is the
Teichm$\ddot{\mathrm{u}}$ller character.
By inverting $p$, we get a surjection $\widetilde{\mathbb{B}}^+\rightarrow \mathbb{C}_p$.
We put $B_{\mathrm{dR}}^+:=\varprojlim_n\widetilde{\mathbb{B}}^+/(\mathrm{Ker}(\theta))^n$, which is
a complete discrete valuation ring with residue field $\mathbb{C}_p$ and 
is equipped with the projective limit topology of the $p$-adic topology 
on $\widetilde{\mathbb{B}}^+/(\mathrm{Ker}(\theta))^n$. 
Let $A_{\mathrm{max}}$ be the $p$-adic completion of $\widetilde{\mathbb{A}}^+[\frac{[\tilde{p}]}{p}]$, where $\tilde{p}:=(p^{(n)})$ is an element in $\widetilde{\mathbb{E}}^+$ such that $p^{(0)}=p, (p^{(n+1) })^p=p^{(n)}$. We put $B_{\mathrm{max}}^+:=A_{\mathrm{max}}[\frac{1}{p}]$. $A_{\mathrm{max}}$ and  
$B^+_{\mathrm{max}}$ have continuous $G_K$-actions and Frobenius actions $\varphi$, here the topology on these rings are also $p$-adic topology. We have a natural continuous $G_K$-equivariant embedding $B_{\mathrm{max}}^+\hookrightarrow B_{\mathrm{dR}}^+$.
If we put $t:=\mathrm{log}([\varepsilon])$, then we can see that $t\in A_{\mathrm{max}}$, 
$\varphi(t)=pt, g(t)=\chi(g)t$ for any $g\in G_K$ and $\mathrm{Ker}(\theta)=tB_{\mathrm{dR}}^+
\subset B_{\mathrm{dR}}^+$ is the maximal ideal of $B^+_{\mathrm{dR}}$. If we put $B_{\mathrm{max}}:=B_{\mathrm{max}}^+
[\frac{1}{t}], B_{\mathrm{dR}}:=B_{\mathrm{dR}}^+[\frac{1}{t}]$, which are 
equipped 
with the inductive limit topology of $\frac{1}{t^n}B^+_{\mathrm{max}}$ and 
$\frac{1}{t^n}B^+_{\mathrm{dR}}$for any $n\in\mathbb{N}$, we have a natural continuous embedding 
$B_{\mathrm{max}}\hookrightarrow B_{\mathrm{dR}}$. We put $B_e:=B_{\mathrm{max}}^{\varphi =1}$, $(B_e\subseteq)\widetilde{B}^+_{\mathrm{rig}}:=\cap_{n=0}^{\infty}\varphi^{n}(B_{\mathrm{max}})
\subset B_{\mathrm{max}}$ (these are closed sub rings of $B_{\mathrm{max}}$), $\mathrm{Fil}^i B_{\mathrm{dR}}:=t^i B_{\mathrm{dR}}^+$ for any $i\in\mathbb{Z}$. We put $\mathrm{log}([\tilde{p}]):=\mathrm{log}(p)+
\sum_{n=1}^{\infty}\frac{(-1)^{(n-1)}}{n}(\frac{[\tilde{p}]}{p}-1)^{n}\in B^+_{\mathrm{dR}}$, here $\mathrm{log}:\mathbb{C}_p^{\times}\rightarrow \mathbb{C}_p$ is a branch of $\mathrm{log}$ which we fix in this article. Let $B_{\mathrm{log}}:=B_{\mathrm{max}}[\mathrm{log}([\tilde{p}])]\subseteq B_{\mathrm{dR}}$. There is a derivation $N:B_{\mathrm{log}}\rightarrow B_{\mathrm{log}}$ over $B_{\mathrm{max}}$ such that $N(\mathrm{log}([\tilde{p}])):=-1$. We have the following fundamental short exact sequence $\cite[\mathrm{Proposition}\,1.17]{BK}$
\begin{equation*}
0\rightarrow \mathbb{Q}_p\rightarrow B_e\oplus B_{\mathrm{dR}}^+\rightarrow B_{\mathrm{dR}}\rightarrow 0.
\end{equation*}
\begin{defn}
An $E$-representation of $G_K$ is a finite dimensional $E$-vector space $V$ with a continuous $E$-linear action of $G_K$. We call $E$-representation for simplicity when there will be no risk 
of confusion about $K$.
\end{defn}

Next we define $E$-$B$-pair of $G_K$, which is the $E$-coefficient version of $B$-pair. 
\begin{defn}\label{-3}
An $E$-$B$-pair of $G_K$ is a couple $W=(W_e, W^+_{\mathrm{dR}})$ such that
\begin{itemize}
\item[(1)]$W_e$ is a finite $B_e\otimes_{\mathbb{Q}_p}E$-module with a continuous 
semi-linear $G_K$-action which is free as $B_e$-module,
\item[(2)]$W_{\mathrm{dR}}^+\,\subseteq W_{\mathrm{dR}}:=B_{\mathrm{dR}}\otimes_{B_e}W_e$ is a $G_K$-stable $B_{\mathrm{dR}}^+\otimes_{\mathbb{Q}_p}E$-lattice, i.e. $W_{\mathrm{dR}}^+$ is a finitely generated $B^+_{\mathrm{dR}}\otimes_{\mathbb{Q}_p}E$-module which generates $W_{\mathrm{dR}}$ as a $B_{\mathrm{dR}}\otimes_{\mathbb{Q}_p}E$-module.
\end{itemize}
Here $G_K$ acts on $B_{?}\otimes_{\mathbb{Q}_p}E$ by $g(x\otimes y):=g(x)\otimes y$ 
for any $g\in G_K, x\in B_{?}, y\in E$ for $?\in\{e, \mathrm{dR}\}$. 
\end{defn}

We call $E$-$B$-pair for simplicity when there will be no risk 
of confusing about $K$. And we simply call a $B$-pair when $E=\mathbb{Q}_p$, then this definition 
is the same as that of Berger $\cite[\mathrm{Introduction}]{Be07}$.

\begin{rem}
Later we will prove that $W_e$ is also free over $B_e\otimes_{\mathbb{Q}_p}E$ 
(Lemma $\ref{0}$) and that $W^+_{\mathrm{dR}}$ is also free over $B^+_{\mathrm{dR}}\otimes_{\mathbb{Q}_p}E$ (Lemma $\ref{-2}$).
\end{rem}

\begin{defn}
Let $W_j:=(W_{e,j}, W^+_{\mathrm{dR},j})$ be $E$-$B$-pairs 
for $j=1,2$.
 Then a morphism of $E$-$B$-pairs $f:W_1\rightarrow W_2$ is defined as a $B_e\otimes_{\mathbb{Q}_p}E$-semi-linear $G_K$-equivariant morphism $f:W_{e,1}\rightarrow W_{e,2}$ such that $id_{B_{\mathrm{dR}}}\otimes_{B_e}f:B_{\mathrm{dR}}\otimes_{B_e}W_{e,1}
\rightarrow B_{\mathrm{dR}}\otimes_{B_e}W_{e,2}$ maps $W^+_{\mathrm{dR},1}$ to 
$W^+_{\mathrm{dR},2}$.
\end{defn}

\begin{rem}
Let $V$ be an $E$-representation of $G_K$. Then $W(V):=(B_e\otimes_{\mathbb{Q}_p}V, B^+_{\mathrm{dR}}\otimes_{\mathbb{Q}_p}V)$ is an $E$-$B$-pair of $G_K$. By the fundamental short exact sequence $0\rightarrow \mathbb{Q}_p\rightarrow 
B_e\oplus B^+_{\mathrm{dR}}\rightarrow B_{\mathrm{dR}}\rightarrow 0$, it is easy to see that the functor $V\mapsto W(V)$ is a fully faithful functor from the category of $E$-representations of $G_K$ to the category of $E$-$B$-pairs of $G_K$ ($\cite[\mathrm{Introduction}]
{Be07}$).
\end{rem}

Next, we prove a technically important lemma concerning 
to the B$\mathrm{\acute{e}}$zout property of $B_e\otimes_{\mathbb{Q}_p}E$. This 
lemma is a generalization of $\cite[\mathrm{Proposition}\, 1.1.9]{Be07}$ to any coefficent case.
\begin{lemma}\label{-1}
$B_e\otimes_{\mathbb{Q}_p}E$ is a B$\mathrm{\acute{e}}$zout domain i.e. a domain and every finitely generated ideal is generated by one element.
\end{lemma}
\begin{proof}
The proof is essentially same as that of $\cite[\mathrm{Proposition}\, 1.1.9]{Be07}$. First, we have a natural isomorphism of rings 
$B_e\otimes_{\mathbb{Q}_p}E\isom B_e\otimes_{\mathbb{Q}_p}
(E_0\otimes_{E_0}E)\isom B_{\mathrm{max}}^{\varphi^{f'}=1}\otimes_{E_0}E$.
(Here $E_0$ is the maximal unramified extension of $\mathbb{Q}_p$ in $E$ and 
$f'=[E_0:\mathbb{Q}_p]$.)
Because the natural map $B_{\mathrm{max}}^{\varphi^{f'}=1}\otimes_{E_0}E
\hookrightarrow B_{\mathrm{dR}}$ is injective by $\cite[\mathrm{Proposition}\,
7.14]{Co02}$, so 
$B_e\otimes_{\mathbb{Q}_p}E$ is a domain. 
Next, we show that for any $f,g\in B_{\mathrm{max}}^{\varphi^{f'}=1}\otimes_{E_0}E$,
 the ideal generated by $f$ and $g$ in $B_{\mathrm{max}}^{\varphi^{f'}=1}\otimes_{E_0}E$ is generated by one element.
For this, we first note that $\widetilde{B}^{\dagger}_{\mathrm{rig}}\otimes_{E_0}
E$ (for the definition of $\widetilde{B}^{\dagger}_{\mathrm{rig}}$, see $1.2$ of this paper.) is a B$\mathrm{\acute{e}}$zout domain by $\cite[\mathrm{Theorem}2.9.6]{Ke05}$ (In the definition of $\cite[2.1]{Ke05}$, if we take 
$K_0:=k_{E}$ the residue field of $E$, $\mathcal{O}:=\mathcal{O}_E$, $\sigma
:=\varphi^{f'}\otimes_{W(k_{E})}\mathrm{id}_{\mathcal{O}_E}$, then $\Gamma^{\mathrm{alg}}_{\mathrm{an},\mathrm{con}}\isom \widetilde{B}^{\dagger}_{\mathrm{rig}}\otimes_{E_0}E$.) 
Because we have the injection $B_{\mathrm{max}}^{\varphi^{f'}=1}\otimes_{E_0}E\hookrightarrow 
\widetilde{B}^{\dagger}_{\mathrm{rig}}[\frac{1}{t}]\otimes_{E_0}E$, for large $n\in
\mathbb{Z}_{\ge 0}$ we have $t^nf, t^ng\in \widetilde{B}^{\dagger}_{\mathrm{rig}}
\otimes_{E_0}E$. Because $\widetilde{B}^{\dagger}_{\mathrm{rig}}
\otimes_{E_0}E$ is B$\mathrm{\acute{e}}$zout, there exists an $h\in\widetilde{B}^{\dagger}_{\mathrm{rig}}\otimes_{E_0}E$ such that $f\widetilde{B}^{\dagger}_{\mathrm{rig}}
\otimes_{E_0}E+ g\widetilde{B}^{\dagger}_{\mathrm{rig}}
\otimes_{E_0} E= h\widetilde{B}^{\dagger}_{\mathrm{rig}}
\otimes_{E_0}E$. And $h\widetilde{B}^{\dagger}_{\mathrm{rig}}\otimes_{E_0}E$ is a $\sigma$-module over $\widetilde{B}^{\dagger}_{\mathrm{rig}}
\otimes_{E_0}E$ because $t^nf, t^ng$ are preserved by $\sigma$-action.
So by $\cite[\mathrm{Proposition}\,3.3.2]{Ke05}$, we can choose the generator $h$ such that $\sigma(h)=\pi_E^k h$ for some $k\in\mathbb{Z}$. By using the 
element $t_E\in B^+_{\mathrm{max}}\otimes_{E_0}E$ defined in $\cite[\mathrm{Proposition}\,8.10]{Co02}$ and by $\cite[\mathrm{Lemma}\,8.17]{Co02}$, we get 
$\frac{h}{t_E^k}\in (\widetilde{B}^{\dagger}_{\mathrm{rig}}[\frac{1}{t}])^{\varphi^{f'}=1}
\otimes_{\mathrm{E}_0}E =B_{\mathrm{max}}^{\varphi^{f'}=1}\otimes_{E_0}E$ 
(Here, for the last equality, we use the fact $(\widetilde{B}^{\dagger}_{\mathrm{rig}}[\frac{1}{t}])^{\varphi^{f'}=1}=B_{\mathrm{max}}^{\varphi^{f'}=1}$, see $\cite[\mathrm{Lemma}\,1.1.7]{Be07}$). 
Then we can show that the ideal generated by $f$ and $g$  in $B_{\mathrm{max}}^{\varphi^{f'}=1}\otimes_{E_0}E$ is generated by $\frac{h}{t_E^k}$ in the same way as in $\cite[\mathrm{Proposition}\,1.1.9]{Be07}$.

\end{proof}
From the above lemma, we get the following lemmas which are also
technically important.
\begin{lemma}\label{0}
Let $W_e$ be a finite $W_e\otimes_{\mathbb{Q}_p}E$ module 
which is free over $B_e$. Then $W_e$ is also 
free over $B_e\otimes_{\mathbb{Q}_p}E$.
In particular, for any $E$-$B$-pair $W:=(W_e,W^+_{\mathrm{dR}})$, 
$W_e$ is finite free over $B_e\otimes_{\mathbb{Q}_p}E$.
\end{lemma}
\begin{proof}
Because $W_e$ is finitely generated over $B_e\otimes_{\mathbb{Q}_p}E$ which 
is B$\mathrm{\acute{e}}$zout domain, by the remark after Lemma 2.4 $\cite{Ke04}$ it suffices to show that $W_e$ is 
 torsion free over $B_e\otimes_{\mathbb{Q}_p}E$. 
Let $\mathrm{Frac}B_e$ be the fraction field of $B_e$. 
Because $B_e\otimes_{\mathbb{Q}_p}E$ is a domain and $(\mathrm{Frac}B_e)
\otimes_{\mathbb{Q}_p}E$ is a localization of $B_e\otimes_{\mathbb{Q}_p}E$, 
so $(\mathrm{Frac}B_e)\otimes_{\mathbb{Q}_p}E$ is also a domain and 
the natural map $B_e\otimes_{\mathbb{Q}_p}E\hookrightarrow (\mathrm{Frac}B_e)
\otimes_{\mathbb{Q}_p}E$ is injective.
And because $(\mathrm{Frac}B_e)\otimes_{\mathbb{Q}_p}E$ is finite 
over the field $\mathrm{Frac}B_e$, so $(\mathrm{Frac}B_e)\otimes_{\mathbb{Q}_p}E$
 is also a field. Because $W_e$ is free over $B_e$ by assumption, the 
natural map $W_e\hookrightarrow \mathrm{Frac}B_e\otimes_{B_e}W_e$ is injective.
And of course $\mathrm{Frac}B_e\otimes_{B_e}W_e$ is torsion free over $(\mathrm{Frac}B_e)\otimes_{\mathbb{Q}_p}E$ which is a field. By these, we conclude that $W_e$ is torsion 
free over $B_e\otimes_{\mathbb{Q}_p}E$.
\end{proof}
\begin{lemma}\label{-2}
Let $W:=(W_e,W^+_{\mathrm{dR}})$ be an $E$-$B$-pair. Then 
$W^+_{\mathrm{dR}}$ is finite free over $B^+_{\mathrm{dR}}\otimes_{\mathbb{Q}_p}E$.
\end{lemma}
\begin{proof}
By lemma $\ref{0}$, $W_e$ is free over $B_e\otimes_{\mathbb{Q}_p}E$. So 
$W_{\mathrm{dR}}:=B_{\mathrm{dR}}\otimes_{B_e}W_e$ is also free over 
$B_{\mathrm{dR}}\otimes_{\mathbb{Q}_p}E$. Because $E\subseteq B^+_{\mathrm{dR}}$, we have natural isomorphisms $B^+_{\mathrm{dR}}\otimes_{\mathbb{Q}_p}E\isom 
\oplus_{\sigma:E\hookrightarrow B^+_{\mathrm{dR}}}B^+_{\mathrm{dR}} :a\otimes b \mapsto (a\sigma(b))_{\sigma}$ ($a\in B^+_{\mathrm{dR}}$, $b\in E$) and $B_{\mathrm{dR}}\otimes_{\mathbb{Q}_p}E\isom \oplus_{\sigma:E\hookrightarrow B^+_{\mathrm{dR}}}B_{\mathrm{dR}}$. By using these decompositions, we get decompositions $W^+_{\mathrm{dR}}\isom \oplus_{\sigma:E\hookrightarrow B^+_{\mathrm{dR}}}W^+_{\mathrm{dR},\sigma}$ and $W_{\mathrm{dR}}\isom \oplus_{\sigma:E\hookrightarrow B^+_{\mathrm{dR}}}W_{\mathrm{dR},\sigma}$ ($W^+_{\mathrm{dR},\sigma}$ and $W_{\mathrm{dR},\sigma}$ are the $\sigma$-components). Then, for each $\sigma$, $W^+_{\mathrm{dR},\sigma}$ is a $B^+_{\mathrm{dR}}$-lattice of $W_{\mathrm{dR},\sigma}$. So $W^+_{\mathrm{dR},\sigma}$ is finite free 
over $B^+_{\mathrm{dR}}$ of same rank as that of $W_{\mathrm{dR},\sigma}$ over $B_{\mathrm{dR}}$ because $B^+_{\mathrm{dR}}$ is a discrete valuation 
ring. And bacause $W_{\mathrm{dR}}$ is free over $B_{\mathrm{dR}}\otimes_{\mathbb{Q}_p}E$, so $W_{\mathrm{dR},\sigma}$ have same rank for any $\sigma$. 
So $W^+_{\mathrm{dR},\sigma}$ also have same rank for any $\sigma$. 
So $W^+_{\mathrm{dR}}$ is finite free over $B^+_{\mathrm{dR}}\otimes_{\mathbb{Q}_p}E$.
\end{proof}

By using these lemmas, we can define rank, tensor products and duals of $E$-$B$-pairs.
\begin{defn}
\begin{itemize}
\item[(1)]Let $W:=(W_e,W^+_{\mathrm{dR}})$ be an $E$-$B$-pair, then 
we define the rank of $W$ by $\mathrm{rank}(W):=\mathrm{rank}_{B_e\otimes_{\mathbb{Q}_p}E}(W_e)$.
\item[(2)]Let $W_1:=(W_{e,1},W^+_{\mathrm{dR},1})$ and $W_2:=
(W_{e,2},W^+_{\mathrm{dR},2})$ be $E$-$B$-pairs.
Then we define the tensor product of $W_1$ and $W_2$ by 
$W_1\otimes W_2:= (W_{e,1}\otimes_{B_e\otimes_{\mathbb{Q}_p}E}W_{e,2}, 
W^+_{\mathrm{dR},1}\otimes_{B^+_{\mathrm{dR}}\otimes_{\mathbb{Q}_p}E}
W^+_{\mathrm{dR},2})$.
\item[(3)]
Let $W:=(W_e,W^+_{\mathrm{dR}})$ be an $E$-$B$-pair. Then 
we define the dual of $W$ by $W^{\vee}:=(\mathrm{Hom}_{B_e}(W_e,B_e), 
\mathrm{Hom}_{B^+_{\mathrm{dR}}}(W^+_{\mathrm{dR}}, B^+_{\mathrm{dR}}))$.
Here, we define the $E$-action on $W^{\vee}$ by $af(x):=f(ax)$ for any 
$a\in E$, $f\in \mathrm{Hom}_{B_e}(W_e,B_e)$ (resp. $f\in\mathrm{Hom}_{B^+_{\mathrm{dR}}}(W^+_{\mathrm{dR}}, B^+_{\mathrm{dR}}))$), $x\in W_e$ (resp. $x\in W^+_{\mathrm{dR}}$).
\end{itemize}
\end{defn}

\begin{lemma}\label{1}
Let $W_1$, $W_2$ be finite free $B_e\otimes_{\mathbb{Q}_p}E$-modules with continuous 
semi-linear $G_K$-actions.
Let $f:W_1\rightarrow W_2$ be a $B_e\otimes_{\mathbb{Q}_p}E$-semi-linear $G_K$-morphism.
Then $\mathrm{Ker}(f)$, $\mathrm{Im}(f)$ and $\mathrm{Cok}(f)$ are all finite free over $B_e\otimes_{\mathbb{Q}_p}E$.
\end{lemma}
\begin{proof}
Because $\mathrm{Im}(f)$ is finite torsion free $B_e\otimes_{\mathbb{Q}_p}E$-module, it is free by the remark after Lemma 2.4 $\cite{Ke04}$. 
From this, we get a splitting (as $B_e\otimes_{\mathbb{Q}_p}E$-modules) 
of the short exact sequence  
$0\rightarrow \mathrm{Ker}(f)\rightarrow W_1\rightarrow \mathrm{Im}(f)\rightarrow 0$. So $\mathrm{Ker}(f)$ is finite over $B_e\otimes_{\mathbb{Q}_p}E$, so $\mathrm{Ker}(f)$ is finite torsion free over $B_e\otimes_{\mathbb{Q}_p}E$. 
So $\mathrm{Ker}(f)$ is also finite free.
By $\cite[\mathrm{Lemma}\, 2.1.4]{Be07}$, $\mathrm{Cok}(f)$ is free over $B_e$.
Then by  Lemma $\ref{0}$, $\mathrm{Cok}(f)$ is free over $B_e\otimes_{\mathbb{Q}_p}E$.
\end{proof}
The category of $E$-$B$-pairs is not an abelian category since cokernels of 
morphisms do not exist in general.
We define the exactness in the category of $B$-pairs.
\begin{defn}
Let $W_i:=(W_{e,i},W^+_{\mathrm{dR},i})$ be $E$-$B$-pairs of $G_K$ for $i=1,2,3.$
Let $f:W_1\rightarrow W_2, g:W_2\rightarrow W_3$ be morphisms of $E$-$B$-pairs.
Then we say that
\begin{equation*}
0\rightarrow W_1\rightarrow W_2\rightarrow W_3\rightarrow 0
\end{equation*}
is exact if the following two sequences are exact in usual sense
\begin{equation*}
0\rightarrow W_{e,1}\rightarrow W_{e,2}\rightarrow W_{e,3}\rightarrow 0,
\end{equation*}
\begin{equation*}
0\rightarrow W^+_{\mathrm{dR},1}\rightarrow W^+_{\mathrm{dR},2}\rightarrow W^+_{\mathrm{dR},3}\rightarrow 0.
\end{equation*}
\end{defn}

\begin{lemma}\label{2}
Let $W_1:=(W_{e,1},W^+_{\mathrm{dR},1})$, $W_2:=(W_{e,2},W^+_{\mathrm{dR},2})$ be $E$-$B$-pairs. Let $f:W_1\rightarrow W_2$ be a morphism of $E$-$B$-pairs.
We put $\mathrm{Ker}(f):=(\mathrm{Ker}(f_e:W_{e,1}\rightarrow W_{e,2}), \mathrm{Ker}(f_{\mathrm{dR}}:W^+_{\mathrm{dR},1}\rightarrow W^+_{\mathrm{dR},2}))$, $\mathrm{Im}(f)
:=(\mathrm{Im}(f_e:W_{e,1}\rightarrow W_{e,2}), \mathrm{Im}(f_{\mathrm{dR}}:W^+_{\mathrm{dR},1}\rightarrow W^+_{\mathrm{dR},2}))$.
Then $\mathrm{Ker}(f)$ and $\mathrm{Im}(f)$ are $E$-$B$-pairs.
\end{lemma}
\begin{proof}
$\mathrm{Ker}(f_e)$ and $\mathrm{Im}(f_e)$ are finite free $B_e\otimes_{\mathbb{Q}_p}E$-modules by Lemma $\ref{1}$. From this, we get the canonical isomorphisms 
$\mathrm{Ker}(id_{B_{\mathrm{dR}}}\otimes_{B_e}f_e)\isom B_{\mathrm{dR}}\otimes_{B_e}\mathrm{Ker}(f_e)$, 
$\mathrm{Im}(id_{\mathrm{dR}}\otimes_{B_e}f_e)\isom B_{\mathrm{dR}}\otimes_{B_e}\mathrm{Im}(f_e)$.
As for $B^+_{\mathrm{dR}}$-modules, $B_{\mathrm{dR}}$ is a flat $B^+_{\mathrm{dR}}$-module because $B^+_{\mathrm{dR}}$ is a principal ideal domain. So we get the canonical isomorphisms 
$\mathrm{Ker}(id_{B_{\mathrm{dR}}}\otimes_{B^+_{\mathrm{dR}}}f_{\mathrm{dR}})\isom B_{\mathrm{dR}}\otimes_{B^+_{\mathrm{dR}}}\mathrm{Ker}(f_{\mathrm{dR}})$, $\mathrm{Im}(id_{B_{\mathrm{dR}}}\otimes_{B^+_{\mathrm{dR}}}f_{\mathrm{dR}})\isom B_{\mathrm{dR}}\otimes_{B^+_{\mathrm{dR}}}\mathrm{Im}(f_{\mathrm{dR}})$.
So we get the natural isomorphisms 
$B_{\mathrm{dR}}\otimes_{B_e}\mathrm{Ker}(f_e)\isom B_{\mathrm{dR}}\otimes_{B^+_{\mathrm{dR}}}\mathrm{Ker}(f_{\mathrm{dR}})$ and $B_{\mathrm{dR}}\otimes_{B_e}\mathrm{Im}(f_e)\isom B_{\mathrm{dR}}\otimes_{B^+_{\mathrm{dR}}}\mathrm{Im}(f_{\mathrm{dR}})$, i.e. $\mathrm{Ker}(f)$ and $\mathrm{Im}(f)$ are $B$-pairs.
Because these maps are all $E$-linear, $\mathrm{Ker}(f)$ and $\mathrm{Im}(f)$ are $E$-$B$-pairs.
\end{proof}
\begin{defn}
Let $W_1\subset W_2$ be two $E$-$B$-pairs such that $W_1$ is a sub $E$-$B$-pair of $W_2$.
Then we say that $W_1$ is saturated in $W_2$ if $W^+_{\mathrm{dR},2}/W^+_{\mathrm{dR},1}$ is a free $B^+_{\mathrm{dR}}$-module. Then $W_2/W_1:=(W_{e,2}/W_{e,1}, W^+_{\mathrm{dR},2}/W^+_{\mathrm{dR},1})$ is an $E$-$B$-pair by Lemma $\ref{1}$. 
\end{defn}
\begin{lemma}{\label{3}}
Let $W_1\subset W_2$ be two $E$-$B$-pairs.
Then there exists unique $E$-$B$-pair $W_1^{\mathrm{sat}}:=(W_{e,1}^{\mathrm{sat}}, W_{\mathrm{dR},1}^{+,\mathrm{sat}})$ such that 
$W_1\subset W_1^{\mathrm{sat}}\subset W_2$, $W_{e,1}=W_{e,1}^{\mathrm{sat}}$ and $W_1^{\mathrm{sat}}$ is saturated in $W_2$.
We call $W_1^{\mathrm{sat}}$ the saturation of $W_1$ in $W_2$.
\end{lemma}
\begin{proof}
We put $W_{e,1}^{\mathrm{sat}}:=W_{e,1}$ and put $W_{\mathrm{dR},1}^{+,\mathrm{sat}}:=W_{\mathrm{dR},1}\cap W_{\mathrm{dR},2}^+$. Then it is easy to see that 
$W_1^{\mathrm{sat}}:=(W_{e,1}^{\mathrm{sat}}, W^{+,\mathrm{sat}}_{\mathrm{dR},1})$ is an $E$-$B$-pair satisfying all the desired conditions. Conversely, 
if $W'_1:=(W_{e,1}', W^{'+}_{\mathrm{dR},1})$ satisfies the same conditions, it is easy to see that $W^{'+}_{\mathrm{dR},1}$ satisfies $W^{'+}_{\mathrm{dR},1}=W_{\mathrm{dR},1}\cap W^+_{\mathrm{dR},2}$. Uniqueness of $W_1^{\mathrm{sat}}$ follows from this.
\end{proof}

Now we can define trianguline or split trianguline $E$-representations 
and trianguline  or split trianguline $E$-$B$-pairs.
\begin{defn}$\label{a}$
\begin{itemize}
\item[(1)]Let $W$ be an $E$-$B$-pair. We say that $W$ is a split 
trianguline $E$-$B$-pair if there is a filtration $0=W_0\subset W_1\subset \cdots \subset W_l=W$ by sub $E$-$B$-pairs such that for any $i$, $W_{i}$ is saturated in $W_{i+1}$ and the quotient $W_{i+1}/W_{i}$ is a rank one $E$-$B$-pair.
\item[(1)']Let $W$ be an $E$-$B$-pair. We say that $W$ is a trianguline 
$E$-$B$-pair if $W\otimes_E E':=(W_e\otimes_E E', W^+_{\mathrm{dR}}\otimes_E E')$ is split triangline for some finite extension $E'$ of $E$.
\item[(2)]Let $V$ be an $E$-representation. We say that $V$ is a split 
trianguline (resp. trianguline) $E$-representation if $W(V)$ is a split 
trianguline (resp. trianguline) $E$-$B$-pair.
\end{itemize}
\end{defn}

In this paper, we classify two dimensional split 
trianguline $E$-representations.

Next we recall the generalization of the usual $p$-adic Hodge theory to the
case of $B$-pairs following $\cite[2.3]{Be07}$. First we recall the definition of filtered $(\varphi,N,G_K)$-module over $K$.
\begin{defn}
Let $L$ be a finite Galois extension of $K$.
An $E$-filtered $(\varphi,N,\mathrm{Gal}(L/K))$-module over $K$ is a finite 
$L_0\otimes_{\mathbb{Q}_p}E$-module $D$ (where $L_0$ is the maximal unramified extension 
of $\mathbb{Q}_p$ in $L$) such that
\begin{itemize}
\item[(1)]$D$ has a Frobenius semi-linear operator $\varphi_D:D\hookrightarrow D$ such that $id_{L_0}\otimes \varphi_D:L_{0}\otimes_{\varphi,L_0}D\isom D$ is 
an $L_0\otimes_{\mathbb{Q}_p}E$-linear isomorphism. (Here $\varphi$ acts on 
$L_0\otimes_{\mathbb{Q}_p}E$ by $\varphi(x\otimes y)=\varphi(x)\otimes y$ for 
any $x\in L_0$, $y\in E$.)
\item[(2)]$N:D\rightarrow D$ is an $L_0\otimes_{\mathbb{Q}_p}E$-linear morphism such that $p\varphi N= N\varphi$.
\item[(3)]$D_L:=L\otimes_{L_0}D$ has a decreasing filtration by sub $L\otimes_{\mathbb{Q}_p}E$-modules $\mathrm{Fil}^iD_L$ for $i\in\mathbb{Z}$ such that $\mathrm{Fil}^{-i}D_L=D_L$ and $\mathrm{Fil}^iD_L=0$ for sufficiently large $i\ge 0$.
\item[(4)]$\mathrm{Gal}(L/K)$ acts $L_0\otimes_{\mathbb{Q}_p}E$ (or $L\otimes_{\mathbb{Q}_p}E$)-semi-linearly on $D$ (or $D_L$) such that $g\varphi=\varphi g$ and $g N=N g$ and $g(\mathrm{Fil}^iD_L)=\mathrm{Fil}^iD_L$ for any $g\in \mathrm{Gal}(L/K)$ and $i\in\mathbb{Z}$. (Here $\mathrm{Gal}(L/K)$ acts on 
$L\otimes_{\mathbb{Q}_p}E$ by $g(x\otimes y)=g(x)\otimes y$ for 
any $x\in L$, $y\in E$ and $g\in\mathrm{Gal}(L/K)$.)
\end{itemize}
We call $D$ an $E$-filtered $(\varphi,N,G_K)$-module if $D$ is an $E$-filtered $(\varphi, N,\mathrm{Gal}(L/K))$-module for some finite Galois extension 
$L$ of $K$.
\end{defn}

\begin{defn}
Let $W:=(W_e, W^+_{\mathrm{dR}})$ be an $E$-$B$-pair of $G_K$ and let 
$L$ be a finite Galois extension of $K$. Then we define
 $D^{L}_{\mathrm{cris}}(W):=(B_{\max}\otimes_{B_e}W_e)^{G_L}$, $D^{L}_{\mathrm{st}}(W)
 :=(B_{\mathrm{log}}\otimes_{B_e}W_e)^{G_L}$, $D^{L}_{\mathrm{dR}}:=(B_{\mathrm{dR}}\otimes
 _{B_e}W_e)^{G_L}$, i.e. fixed parts of $G_L$.
As in the case of usual $p$-adic representations, we can show that 
$\mathrm{dim}_{L'}(D^{L}_{?}(W))\le \mathrm{rank}_{B_e}(W_e)$ where $L'=L_0$ if $?\in
\{\mathrm{cris}, \mathrm{st}\}$ and $L'=L$ if $?=\mathrm{dR}$.
We say that $W$ is potentially cristalline (resp. potentially semi-stable, resp. de Rham) if 
$\mathrm{dim}_{L'}(D^{L}_{?}(W))=\mathrm{rank}_{B_e}(W_e)$ for $?=\mathrm{cris}$ (resp. $?=\mathrm{st}$, resp. $?=\mathrm{dR}$) for some  finite Galois extension $L$ of $K$. 
\end{defn}
For a potentially semi-stable $E$-$B$-pair $W$ and for sufficiently large $L$ such that $\mathrm{dim}_{L_0}D^L_{\mathrm{st}}(W)\allowbreak =\mathrm{rank}_{B_e}W_e$, 
we can equip $D^{L}_{\mathrm{st}}(W)$ 
with an $E$-filtered $(\varphi, N,\mathrm{Gal}(L/K))$-module structure as follows. The $(\varphi,N,\mathrm{Gal}(L/K))$-module strucure on $D^{L}_{\mathrm{st}}(W)$ is induced from the action of $\varphi$ and $N$ on $B_{\mathrm{log}}$ and from the action of $G_K$ on $B_{\mathrm{log}}\otimes_{B_e}W_e$. The filtration on $D^{L}_{\mathrm{dR}}(W)\isom L\otimes_{L_0}D^L_{\mathrm{st}}(W)$ is defined by $\mathrm{Fil}^i
D^L_{\mathrm{dR}}(W):=D^L_{\mathrm{dR}}(W)\cap t^i W^+_{\mathrm{dR}}\subseteq W_{\mathrm{dR}}$. 
So we get a functor $W\mapsto D^{L}_{\mathrm{st}}(W)$ from the category 
of potentially semi-stable $E$-$B$-pairs of $G_K$ which are semi-stable $E$-$B$-pairs 
of $G_L$ to the category of $E$-filtered $(\varphi, N,\mathrm{Gal}(L/K))$-modules over $K$. Berger generalized $p$-adic monodromy theorem and 
``weakly admissible implies admissible" theorem to the case of $B$-pairs.
\begin{thm}\label{4}
\begin{itemize}
\item[$(1)$]All de Rham $E$-$B$-pairs are potentially semi-stable.
\item[$(2)$]The functor $W\mapsto D^{L}_{\mathrm{st}}(W)$ realizes an equivalence
of categories between the category of potentially semi-stable $E$-$B$-pairs of $G_K$ which are semi-stable $E$-$B$-pairs 
of $G_L$ to the category of $E$-filtered $(\varphi, N,\mathrm{Gal}(L/K))$-modules over $K$.
\item[$(3)$]The functor $W\mapsto D^{L}_{\mathrm{cris}}(W)$ realizes an equivalence
of categories between the category of potentially cristalline $E$-$B$-pairs of $G_K$ which are cristalline $E$-$B$-pairs of $G_L$ to 
the category of $E$-filtered $(\varphi,\mathrm{Gal}(L/K))$-modules over $K$.
\end{itemize}
\end{thm}
\begin{proof}
$\cite[\mathrm{Proposition}\, 2.3.4]{Be07}$, $\cite[\mathrm{Theorem}\, 2.3.5]
{Be07}$.
\end{proof}
\begin{rem}
An inverse functor of $D^L_{\mathrm{st}}$ is defined as follows.
For an $E$-filtered $(\varphi, N,\allowbreak \mathrm{Gal}(L/K))$-module $D$ over $K$, put $W_e(D):=
(B_{\mathrm{st}}\otimes_{L_0}D)^{\varphi=1, N=0}$ and $W^+_{\mathrm{dR}}(D)
:=\mathrm{Fil}^{0}(B_{\mathrm{dR}}\otimes_{L}D_L)$. Then we can show that
$W(D):=(W_e(D), W^+_{\mathrm{dR}}(D))$ is an $E$-$B$-pair of $G_K$.
$D\mapsto W(D)$ is an inverse of $W\mapsto D^L_{\mathrm{st}}(W)$ 
($\cite[2.3]{Be07}$ ).
\end{rem}
\begin{rem}
The proof of $(2)$ and $(3)$ of the above theorem is much easier 
than that in the case of $p$-adic representations, because in this case 
there are no conditions about weakly admissibility of filtered $(\varphi,N,G_K)$-modules.
\end{rem}

\subsection{$(\varphi, \Gamma_K)$-modules over Robba ring.}

In the paper $\cite{Be07}$, Berger established the equivalence between the
category of $B$-pairs and the category of $(\varphi, \Gamma_K)$-modules over Robba ring $B^{\dagger}_{\mathrm{rig},K}$. In this subsection, we first recall the general facts about $\varphi$-modules over Robba rings following $\cite{Ke}$ and next recall the construction of Robba ring $B^{\dagger}_{\mathrm{rig},K}$ and the definition of $(\varphi,\Gamma_K)$-modules over $B^{\dagger}_{\mathrm{rig},K}$ following $\cite{Be02}$. 

In applications of $(\varphi,\Gamma_K)$-modules to $p$-adic Hodge theory, the notion of slope 
and the slope filtration theorem of Kedlaya are very important.
So we recall the general fact about slopes of $\varphi$-modules over general Robba ring and the slope filtration theorem of Kedlaya following $\cite[1]{Ke}$.
Let $L$ be a complete discrete valuation field of characteristic zero, $\mathcal{O}_L$ its integer ring, $k_L$ the residue field of $L$. 
$\pi_L$ is  a uniformizer of $L$. 
Then we define the Robba ring of $L$ by $\mathcal{R}_L
:=\{f(x):=\sum_{k\in\mathbb{Z}}a_n x^n| a_n\in L$, $f(x)$ converges for any $x\in L$ such that $r\le |x|<1$
 for some $r<1 \}$. We assume that there is an endomorphism $\phi_L:L\rightarrow L$ which is  a lifting of the $q=p^f$-th power Frobenius on $k_L$ for some $f\in\mathbb{N}$ and we assume that there is a $\phi_L$-semi-linear endomorphism $\phi:\mathcal{R}_L\rightarrow \mathcal{R}_L:
 \sum_{k\in\mathbb{Z}}a_k x^k\mapsto \sum_{k\in\mathbb{Z}}\phi_L(a_k)\phi(x)^k$
  such that $\phi(x)-x^q\in \mathcal{R}_L$ has all coefficients in $\pi_L\mathcal{O}_L$.
 \begin{defn}
 A $\phi$-module over $\mathcal{R}_L$ is a finite free $\mathcal{R}_L$-module $M$ equipped with $\phi$-semi-linear action $\phi_M:M\rightarrow M$ such that 
 $id_{\mathcal{R}_L}\otimes \phi_M:\mathcal{R}_L\otimes_{\phi,\mathcal{R}_L}M\rightarrow M$ is an $\mathcal{R}_L$-linear isomorphism.
 \end{defn}
The category of $\phi$-modules over $\mathcal{R}_L$ is not an abelian category
because the cokernel of a morphism is not a free $\mathcal{R}_L$-module in general. So it is important to know when the cokernel is free. 
\begin{defn}
Let $M_1\subset M$ be a sub $\phi$-module of a $\phi$-module $M$ over $\mathcal{R}_L$. Then we say that $M_1$ is saturated in $M$ if the quotient 
$M/M_1$ is a free $\mathcal{R}_L$ module.
Then $M/M_1$ is a $\phi$-module over $\mathcal{R}_L$.
\end{defn}
Put $\mathcal{R}_L^{\mathrm{bd}}:=\{f(x)\in \mathcal{R}_L| |f(x)| $ is bounded in $r\le |x|<1$ for some $r<1 \}=\{f(x)=\sum_{k\in\mathbb{Z}}a_k x^k|\{ a_k\}_{k\in\mathbb{Z}}\subset L $ is
 bounded and $f(x)$ converges for any $x$ such that $r\le |x|<1$ for some $r<1\}$. Put $\mathcal{R}_L^{\mathrm{int}}
 :=\{\sum_{k\in\mathbb{Z}}a_k x^k\in\mathcal{R}_L| a_k\in\mathcal{O}_L$ for any $k\in\mathbb{Z}\}$. Then $\mathcal{R}_L^{\mathrm{bd}}$ is a discrete valuation field with the integer ring $\mathcal{R}_L^{\mathrm{int}}$ and a valuation
  on $\mathcal{R}_L^{\mathrm{bd}}$ 
is defined by $w(f(x)):=\mathrm{inf}_{n\in\mathbb{Z}}\{v_L(a_n)\}$ for any $f(x)=
\sum_{n\in\mathbb{Z}}a_nx^n \in \mathcal{R}_L^{\mathrm{bd}}$, where $v_L$ is 
the valuation of $L$ such that $v_L(p)=1$. The residue field of $\mathcal{R}_L^{\mathrm{bd}}$
is $k_L((x))$. By the above assumption on $\phi$, $\phi(x)\in \mathcal{R}_L^{\mathrm{int}}$, $\phi$ preserves $\mathcal{R}_L^{\mathrm{bd}}$ and $\mathcal{R}_L^{\mathrm{int}}$ and $w(\phi(f))=w(f)$ for any $f\in\mathcal{R}_L^{\mathrm{bd}}$. Moreover it is known that $\mathcal{R}_L^{\times}=\mathcal{R}_L^{\mathrm{bd},\times}$ ($\cite[\mathrm{Example}\, 1.4.2]{Ke}$).

\begin{defn} 
For a $\phi$-module $M$ over $\mathcal{R}_L$ of rank $n$, the top exterior power $\wedge^n M$ has rank one over $\mathcal{R}_L$. Let $\mathbf{v}$ be a 
generator of $\wedge^n M$ and write $\phi(\mathbf{v})=r\mathbf{v}$ for some $r\in \mathcal{R}_L^{\times}=\mathcal{R}_L^{\mathrm{bd} \times}$. Define the degree of $M$ by setting $\mathrm{deg}(M):=w(r)$, define the slope of $M$ by setting $\mu(M):=\mathrm{deg}(M)/\mathrm{rank}(M)$. It is easy to check that $\mathrm{deg}(M)$ and $\mu(M)$ does not depend on the choice of the generator of $\wedge^n M$.
 \end{defn}
If $M_1$, $M_2$ are $\phi$-modules over $\mathcal{R}_L$, then we can see that
$\mu(M_1\otimes_{\mathcal{R}_L}M_2)=\mu(M_1)+\mu(M_2)$ ($\cite[\mathrm{Remark}
\,1.4.5]{Ke}$).
Next we define pure slope $\phi$-modules.
Before defining this, we need to recall some definitions.
\begin{defn}
Let $M$ be a $\phi$-module over $\mathcal{R}_L$ and let $a$ be a positive integer.
Then we define the $a$-pushforward $[a]_{\ast}M$ of $M$ to be the $\phi^a$-module $M$ over $\mathcal{R}_L$ such that $\phi^a$-semi-linear morphism is defined by $\phi_M^a:M\rightarrow M$.
\end{defn} 
 If $M$ is a $\phi$-module over $\mathcal{R}_L$, then it is easy to see that
$[a]_{\ast}M$ is a $\phi^a$-module over $\mathcal{R}_L$ with $\mathrm{deg}
([a]_{\ast}M)=a\mathrm{deg}(M)$, $\mu([a]_{\ast}M)=a\mu(M)$.
First we define pure slope zero $\phi$-modules, which are called $\acute{\mathrm{e}}$tale $\phi$-modules.
\begin{defn}
A $\phi$-module $M$ over $\mathcal{R}_L$ is said to be $\acute{\mathrm{e}}$tale if it can be 
obtained by base extension from a $\phi$-module over $\mathcal{R}_L^{\mathrm{int}}$, that is, $M$ admits a $\phi_M$-stable $\mathcal{R}_L^{\mathrm{int}}$-lattice $N$ such that $\phi_M$ induces an $\mathcal{R}_L^{\mathrm{int}}$-linear isomorphism $\mathcal{R}_L^{\mathrm{int}}\otimes_{\phi,\mathcal{R}_L^{\mathrm{int}}}N\isom N$.
\end{defn} 
From the above definition, if $M$ is an $\acute{\mathrm{e}}$tale $\phi$-module over $\mathcal{R}_L$, then it is easy to see that $\mathrm{deg}(M)=0$ and $\mu(M)=0$.
On this setting, we define pure slope $\phi$-modules as follows. 
(We put $e_L$ the absolute ramified index of $L$.)
\begin{defn}
Let $M$ be a $\phi$-module over $\mathcal{R}_L$ with slope
$s=c/de_L$ where $c,d$ are coprime integers with $d>0$. We say
 that $M$ is pure of slope $s$ if, for some $\phi^d$-module $N$ of rank one 
such that $\mathrm{deg}(N)=-c/e_L$, $([d]_{\ast}M)\otimes_{\mathcal{R}_L}N$ is an $\acute{\mathrm{e}}$tale $\phi^d$-module.
\end{defn}

\begin{rem}
In $\cite{Be07}$, the definition of pure slope $\phi$-modules 
over Robba ring is different from this definition ($\cite[1.2]{Be07}$).
But it is proved in the proof of $\cite[\mathrm{Theorem}1.7.1]{Ke}$ that two definitions are the same.
\end{rem}
Now we can state the slope filtration theorem of Kedlaya.
\begin{thm}$\label{5}$
Let $M$ be a $\phi$-module over $\mathcal{R}_L$.
Then $M$ admits a unique filtration 
$0=M_0\subset M_1\subset \cdots\subset M_l=M$ by sub $\phi$-modules over 
$\mathcal{R}_L$ such that $M_{i}\subset M_{i+1}$ is saturated and
 $M_{i+1}/M_{i}$  is pure of slope $s_{i+1}$ for any $0\le i\le l-1$ with
$s_1<s_2<\cdots<s_l$.
\end{thm}
\begin{proof}
$\cite[\mathrm{Theorem}1.7.1]{Ke}$.
\end{proof}
Next we recall the construction of Robba ring $B^{\dagger}_{\mathrm{rig},K}$ and some related rings following $\cite{Be02}$.
First, for rational numbers $0\le r\le s\le +\infty$, we put $\widetilde{A}_{[r,s]}:=\widetilde{\mathbb{A}}^{+}\{\frac{p}{[\bar{\pi}^r]}, \frac{[\bar{\pi}^s]}{p}\}$, the p-adic completion of $\widetilde{\mathbb{A}}^{+}[\frac{p}{[\bar{\pi}^r]}, \frac{[\bar{\pi}^s]}{p}]$ (here $\pi:=[\varepsilon]-1$, $\bar{\pi}:=\varepsilon -1$, so $[\bar{\pi}]=[\varepsilon-1]$). When $r=0$, we put $\frac{p}{[\bar{\pi}^r]}:=1$ and 
when $s=+\infty$, we put $\frac{[\bar{\pi}^s]}{p}:=1$. We put $\widetilde{B}_{[r,s]}:=\widetilde{A}_{[r,s]}[\frac{1}{p}]$. Then we have natural continuous 
$G_K$-action on $\widetilde{A}_{[r,s]}$ and on $\widetilde{B}_{[r,s]}$. Frobenius induces isomorphisms $\varphi:\widetilde{A}_{[r,s]}\isom\widetilde{A}_{[pr,ps]}$ and $\varphi:
\widetilde{B}_{[r,s]}\isom\widetilde{B}_{[pr,ps]}$, here we equippe these rings 
with $p$-adic topology. For $r\le r_0\le s_0\le s$, it is known that the natural map $\widetilde{A}_{[r,s]}\hookrightarrow \widetilde{A}_{[r_0,s_0]}$ is injective.
For $r >0$, we put $\widetilde{B}^{\dagger, r}:=\widetilde{B}_{[r,+\infty]}$, 
$\widetilde{B}^{\dagger}:=\cup_{r>0}\widetilde{B}^{\dagger,r}$,
$\widetilde{B}^{\dagger,r}_{\mathrm{rig}}:=\cap_{r\le s< +\infty}\widetilde{B}_{[r,s]}$ (equipped with Frechet topology defined by any $\widetilde{B}_{[r,s]}$) and $\widetilde{B}^{\dagger}_{\mathrm{rig}}:=\cup_{r>0}\widetilde{B}^{\dagger,r}_{\mathrm{rig}}$ (equipped with inductive limit topology). Then we have natural inclusions $\widetilde{B}^{\dagger,r}\subset \widetilde{B}^{\dagger,r}_{\mathrm{rig}}$, $\widetilde{B}^{\dagger}
\subset\widetilde{B}^{\dagger}_{\mathrm{rig}}$ and $\widetilde{B}^{\dagger}
\subset \widetilde{\mathbb{B}}$ (but $\widetilde{B}^{\dagger}_{\mathrm{rig}}
\not\subseteq \widetilde{\mathbb{B}}$). Frobenius induces isomorphisms 
$\varphi:\widetilde{B}^{\dagger,r}\isom\widetilde{B}^{\dagger,pr}$, $\widetilde{B}^{\dagger,r}_{\mathrm{rig}}\isom\widetilde{B}^{\dagger,pr}_{\mathrm{rig}}$,
$\widetilde{B}^{\dagger}\isom\widetilde{B}^{\dagger}$ and $\widetilde{B}^{\dagger}
_{\mathrm{rig}}\isom\widetilde{B}^{\dagger}_{\mathrm{rig}}$.
Moreover we can easily check that $A_{\mathrm{max}}=\widetilde{A}_{[0,(p-1)/p]}$, $B_{\mathrm{max}}=\widetilde{B}_{[0,(p-1)/p]}$ and $\widetilde{B}^+_{\mathrm{rig}}=\widetilde{B}_{[0,+\infty)}:=\cap_{0<s<+\infty}\widetilde{B}_{[0,s]}$.
The rings $\widetilde{B}^{\dagger}_{\mathrm{rig}}, \widetilde{B}^{\dagger}$ are respectively equal to $\Gamma^{\mathrm{alg}}_{\mathrm{an,con}},\Gamma^{\mathrm{alg}}_{\mathrm{con}}$ defined by Kedlaya ($\cite{Ke04}$).
We recall a relation between $\widetilde{B}^{\dagger}_{\mathrm{rig}}$ and $B_{\mathrm{dR}}$.
It is easy to see that $\widetilde{B}_{[(p-1)/p,(p-1)/p]}=\widetilde{\mathbb{A}}^+\{\frac{p}{[\tilde{p}]}, 
\frac{[\tilde{p}]}{p}\}[\frac{1}{p}]$. Then we have a natural continuous 
$G_K$-equivariant injective morphism 
$\widetilde{B}_{[(p-1)/p,(p-1)/p]}\hookrightarrow B^+_{\mathrm{dR}}:\sum_{k\ge 0}a_k(\frac{p}{[\tilde{p}]})^k +\sum_{l\ge 0}b_l(\frac{[\tilde{p}]}{p})^l \mapsto \sum_{k\ge 0}a_k(\frac{p}{[\tilde{p}]})^k +
\sum_{l\ge 0} b_l(\frac{[\tilde{p}]}{p})^l$. In particular, we have a natural continuous inclusion
$i_0:\widetilde{B}^{\dagger, (p-1)/p}_{\mathrm{rig}}=\cap_{(p-1)/p \le r<+\infty}\widetilde{B}_{[(p-1)/p, r]}\hookrightarrow \widetilde{B}_{[(p-1)/p,(p-1)/p]}\hookrightarrow B^+_{\mathrm{dR}}$.
For any $n\in\mathbb{N}$, we define a continuous $G_K$-equivariant injection $i_n:=i_0\circ\varphi^{-n}: \widetilde{B}^{\dagger, (p-1)p^{n-1}}_{\mathrm{rig}}\isom\widetilde{B}^{\dagger,(p-1)/p}_{\mathrm{rig}}\hookrightarrow B^+_{\mathrm{dR}}$.

Next we define $B^{\dagger}_{K,\mathrm{rig}}$ and $B^{\dagger}_{K}$.
We put $A_{K_0}:=\{\sum_{k\ge -\infty}^{+\infty} a_k\pi^k| a_k\in\mathcal{O}_{K_0}, a_k\rightarrow 0\, (k\rightarrow -\infty)\}$ and $B_{K_0}:=A_{K_0}[\frac{1}{p}]$ where $\pi:=[\varepsilon]-1$. $A_{K_0}$ is a complete discrete valuation ring such that $p$ is a prime element, the residue field is $E_{K_0}:=k((\varepsilon-1))$ and the fraction field is $B_{K_0}$. We can show that $A_{K_0}\subset \widetilde{\mathbb{A}}$ and $B_{K_0}\subset \widetilde{\mathbb{B}}$. $\varphi$ and the action of $G_K$ on $\widetilde{\mathbb{A}}$ preserves $A_{K_0}$ and $\varphi(\pi)=(\pi+1)^p-1$ and $g(\pi)=
(\pi+1)^{\chi(g)}-1$ for any $g\in G_K$. Let $\mathbb{A}$ be the $p$-adic completion of the maximal 
unramified extension of $A_{K_0}$ in $\widetilde{\mathbb{A}}$, $\mathbb{B}$ be its fraction field. Then $\varphi$ and the action of $G_K$ also preserve $\mathbb{A}$ and $\mathbb{B}$. We put $A_K:=\mathbb{A}^{H_K}, B_K:=\mathbb{B}^{H_K}, B^{\dagger}_K:=B_K\cap \widetilde{B}^{\dagger}, B^{\dagger,r}_K:=B_K\cap \widetilde{B}^{\dagger,r}$. By definition, these rings are equipped with natural continuous actions of $\varphi$ and $\Gamma_K$. 
If we put $E_K:=(E_{K_0}^{\mathrm{sep}})^{H_K}\subset \tilde{\mathbb{E}}$ (here $E_{K_0}^{\mathrm{sep}}$ is the separable closure of $E_{K_0}$ in $\tilde{\mathbb{E}}$), it is known that $E_K$ is a separable 
extension of $E_{K_0}$ of degree $[K_{\infty}:K_0(\zeta_{p^{\infty}})]$ by the theory of fields of norm. $A_K$ is a complete discrete valuation ring such that $p$ is a prime element, the residue field is $E_K$ and the fraction field is $B_K$. Let $K'_0$ be the maximal unramified extension of $K_0$ in $K_{\infty}$.
For sufficiently large $r>0$ such that a prime element $\pi_{E_K}\in E_K$ lifts 
to $\tilde{\pi}_{E_K}\in B^{\dagger, r}_K$, it is known that $B^{\dagger, r}_K=\{\sum_{k\ge -\infty}^{+\infty}a_k\tilde{\pi}_{E_K}^k| a_k\in K'_0$, $f(X):=\sum_{k\ge -\infty}^{+\infty} a_k X^k$ is a bounded function on $X\in\mathbb{C}_p$ such  that $p^{-1/re_K}\le |X| < 1\}$. We put $B^{\dagger,r}_{\mathrm{rig},K}$ the Frechet completion of $B^{\dagger,r}_{K}$ and $B^{\dagger}_{\mathrm{rig},K}:=\cup_{r\gg 0}B^{\dagger,r}_{\mathrm{rig},K}$. Then it is known that $B^{\dagger,r}_{\mathrm{rig},K}=\{\sum_{k\ge -\infty}^{+\infty} a_k\tilde{\pi}_{E_K}^k| a_k\in K'_0$, $f(X):=
\sum_{k\ge -\infty}^{+\infty}a_k X^k$ converges on $X\in \mathbb{C}_p$ such that $p^{-1/re_K}\le |X|<1\}$. So $B^{\dagger,r}_{\mathrm{rig},K}, B^{\dagger}_{\mathrm{rig},K}$ are the Robba rings with coefficients in $K'_0$. We can show that $B^{\dagger,r}_{\mathrm{rig},K}\subset \widetilde{B}^{\dagger,r}_{\mathrm{rig}}$, $B^{\dagger}_{\mathrm{rig},K}\subset \widetilde{B}^{\dagger}_{\mathrm{rig}}$ and the Frechet topology on $B^{\dagger,r}_{\mathrm{rig},K}$ is the induced topology  from that on $\widetilde{B}^{\dagger,r}_{\mathrm{rig}}$ and $\varphi$ induces inclusions
$B^{\dagger,r}_{\mathrm{rig},K}\hookrightarrow B^{\dagger,pr}_{\mathrm{rig},K}$, $B^{\dagger}_{\mathrm{rig},K}\hookrightarrow B^{\dagger}_{\mathrm{rig},K}$
and $\Gamma_K$ continuously acts on $B^{\dagger,r}_{\mathrm{rig},K}$, $B^{\dagger}_{\mathrm{rig},K}$.
By restricting $i_n:\widetilde{B}^{\dagger,(p-1)p^{n-1}}_{\mathrm{rig}}\hookrightarrow B^+_{\mathrm{dR}}$ to $B^{\dagger,(p-1)p^{n-1}}_{\mathrm{rig},K}$, we have a $G_K$-equivariant injection 
$i_n:B^{\dagger,(p-1)p^{n-1}}_{\mathrm{rig},K}\hookrightarrow B^+_{\mathrm{dR}}$.
\begin{defn}
An $E$-$(\varphi,\Gamma_K)$-module over $B^{\dagger}_{\mathrm{rig},K}$ is a finite $B^{\dagger}_{\mathrm{rig,K}}\otimes_{\mathbb{Q}_p}E$-module $D$ equipped with a Frobenius semi-linear action $\varphi_D$ and a continuous semi-linear action of $\Gamma_K$ such that $D$ is free as 
$B^{\dagger}_{\mathrm{rig},K}$-module, $id_{B^{\dagger}_{\mathrm{rig},K}}\otimes\varphi_D:B^{\dagger}_{\mathrm{rig},K}\otimes_{\varphi,B^{\dagger}_{\mathrm{rig},K}}D\rightarrow D$ is an isomorphism and that the actions of $\varphi$ and $\Gamma_K$ commute. Here $\varphi$ and $\Gamma_K$ act on $B^{\dagger}_{\mathrm{rig},K}\otimes_{\mathbb{Q}_p}E$ as $\varphi\otimes id$ , $\gamma\otimes id$ for any $\gamma\in\Gamma_K$.
\end{defn}
When $E=\mathbb{Q}_p$, we omit the notation $\mathbb{Q}_p$ in the above definition, i.e. we simply call $(\varphi,\Gamma_K)$-modules over $B^{\dagger}_{\mathrm{rig},K}$. 
 
Next we prove a lemma concerning 
the freeness over $B^{\dagger}_{\mathrm{rig},K}\otimes_{\mathbb{Q}_p}E$ of an $E$-$\varphi$-modules.
\begin{lemma}\label{7}
Let $D$ be an $E$- $\varphi$-module over $B^{\dagger}_{\mathrm{rig},K}$.
 Then $D$ is also free as $B^{\dagger}_{\mathrm{rig},K}\otimes_{\mathbb{Q}_p}E$-module.
\end{lemma}
\begin{proof}
Let $E'_0:=E\cap K'_0$, $f':=[E'_0:\mathbb{Q}_p]$ and $I:=\{\sigma:E'_0\hookrightarrow E\}$. Then we have a canonical decomposition $B^{\dagger}_{\mathrm{rig},K}\otimes_{\mathbb{Q}_p}E\isom \oplus_{\sigma\in I} B^{\dagger}_{\mathrm{rig},K}\otimes_{E'_0,\sigma}Ee_{\sigma}$, where $e_{\sigma}$ is the idempotent in $B^{\dagger}_{\mathrm{rig},K}\otimes_{\mathbb{Q}_p}E$ corresponding to $1$ in $B^{\dagger}_{\mathrm{rig},K}\otimes_{E'_0,\sigma}E$. Then any component $B^{\dagger}_{\mathrm{rig},K}\otimes_{E'_0,\sigma}Ee_{\sigma}$ is a Robba ring with coefficients in 
$E$, i.e. they are non canonically isomorphic to $\mathcal{R}_E$.
The action of $\varphi$ is given by $\varphi(a\otimes b e_\sigma)=\varphi(a)\otimes be_{\sigma \varphi^{-1}|_{E'_0}}$ for any $a\in B^{\dagger}_{\mathrm{rig},K}$ and $b\in E$ and $\sigma\in I$. Because $\varphi|_{E'_0}$ transitively acts on $I$, so $\varphi$ transitively acts 
on the components of this decomposition. Because $\sharp I=f'$, $\varphi^{f'}$ preserves the component $B^{\dagger}_{\mathrm{rig},K}\otimes_{E'_0,\sigma}Ee_{\sigma}$ for any $\sigma\in I$.
Let $D$ be an $E$-$\varphi$-module over $B^{\dagger}_{\mathrm{rig},K}$.
Then we also have a canonical decomposition $D\isom\oplus_{\sigma\in I}D\otimes_{B^{\dagger}_{\mathrm{rig},K}\otimes_{\mathbb{Q}_p}E}B^{\dagger}_{\mathrm{rig},K}\otimes_{E'_0,\sigma}Ee_{\sigma}\isom\oplus_{\sigma\in I}D\otimes_{E'_0,\sigma}Ee_{\sigma}$. Then, for any $\sigma\in I$, the component $D\otimes_{E'_0,\sigma}Ee_{\sigma}$ is a $B^{\dagger}_{\mathrm{rig},K}\otimes_{E'_0,\sigma}Ee_{\sigma}$-module which is finite torsion free as $B^{\dagger}_{\mathrm{rig},K}$-module. So it is free as $B^{\dagger}_{\mathrm{rig},K}$-module by $\cite[\mathrm{Proposition}\,2.5]{Ke04}$. Because $B^{\dagger}_{\mathrm{rig},K}\otimes_{E'_0,\sigma}Ee_{\sigma}$ is also a Robba ring, in particular, is a B$\acute{\mathrm{e}}$zout domain and because the natural map $B^{\dagger}_{\mathrm{rig},K}\hookrightarrow B^{\dagger}_{\mathrm{rig},K}\otimes_{E'_0,\sigma}E$ is injective, then it is easy to see that $D\otimes_{E'_0,\sigma}Ee_{\sigma}$ is finite torsion free over $B^{\dagger}_{\mathrm{rig},K}\otimes_{E'_0,\sigma}Ee_{\sigma}$. So it is free over $B^{\dagger}_{\mathrm{rig},K}\otimes_{E'_0,\sigma}Ee_{\sigma}$ for any $\sigma\in I$ by $\cite[\mathrm{Proposition}\,2.5]{Ke04}$.
Then we can take a basis $\{v_1e_{\sigma},\cdots, v_ke_{\sigma}\}$ of $D\otimes_{E'_0,\sigma}Ee_{\sigma}$ over $B^{\dagger}_{\mathrm{rig},K}\otimes_{E'_0,\sigma}Ee_{\sigma}$ for any fixed $\sigma\in I$. Then $\{\varphi^{i}(v_1)e_{\sigma\varphi^{-i}|_{E'_0}},\cdots,\varphi^{i}(v_k)e_{\sigma\varphi^{-i}|_{E'_0}}\}$ is a basis of $D\otimes_{E'_0,\sigma\varphi^{-i}|_{E'_0}}Ee_{\sigma\varphi^{-i}|_{E'_0}}$ for any $0\le i \le f'-1$ by the Frobenius structure on $D$. Then 
$\{\sum_{i=0}^{f'-1}\varphi^{i}(v_1)e_{\sigma\varphi^{-i}|_{E'_0}},\cdots, 
\sum_{i=0}^{f'-1}\varphi^{i}(v_k)e_{\sigma\varphi^{-i}|_{E'_0}}\}$ is a basis of $D$ over $B^{\dagger}_{\mathrm{rig},K}\otimes_{\mathbb{Q}_p}E$. 
So $D$ is a free $B^{\dagger}_{\mathrm{rig},K}\otimes_{\mathbb{Q}_p}E$-module.
\end{proof}
By this lemma, we can define the rank and the tensor products and the duals of $E$-$(\varphi,\Gamma_K)$-modules as follows.
\begin{defn}
\begin{itemize}
\item[(1)]
Let $D$ be an $E$-$(\varphi,\Gamma_K)$-module over $B^{\dagger}_{\mathrm{rig},K}$. We define the rank of $D$ by $\mathrm{rank}(D):=\mathrm{rank}_{B^{\dagger}_{\mathrm{rig},K}\otimes_{\mathbb{Q}_p}E}(D)$, i.e. the rank of $D$ as free $B^{\dagger}_{\mathrm{rig},K}\otimes_{\mathbb{Q}_p}E$-module.
\item[(2)]
Let $D_1$, $D_2$ be $E$-$(\varphi,\Gamma_K)$-modules over $B^{\dagger}_{\mathrm{rig},K}$.
Then we define the tensor product of $D_1$ and $D_2$ by $D_1\otimes D_2:=D_1\otimes_{B^{\dagger}_{\mathrm{rig},K}\otimes_{\mathbb{Q}_p}E}D_2$, which is a free $B^{\dagger}_{\mathrm{rig},K}\otimes_{\mathbb{Q}_p}E$-module by Lemma $\ref{7}$, on which $\varphi$ and $\Gamma_K$ act 
 by $\varphi_{D_1\otimes D_2}:=\varphi_{D_1}\otimes_{B^{\dagger}_{\mathrm{rig},K}\otimes_{\mathbb{Q}_p}E}\varphi_{D_2}$ and $\gamma\otimes_{B^{\dagger}_{\mathrm{rig},K}\otimes_{\mathbb{Q}_p}E}\gamma$ for any $\gamma\in\Gamma_K$.
\item[(3)]
Let $D$ be an $E$-($\varphi,\Gamma)$-module over $B^{\dagger}_{\mathrm{rig},K}$. Then we define the dual of $D$ by $D^{\vee}:=\mathrm{Hom}_{B^{\dagger}_{\mathrm{rig},K}}(D, B^{\dagger}_{\mathrm{rig},K})$, here $E$-action on $D^{\vee}$ is 
defined by $af(x):=f(ax)$ for any $a\in E$, $f\in D^{\vee}$, $x\in D$.
\end{itemize}
 \end{defn}
Next we prove the slope filtration theorem for $E$-$(\varphi,\Gamma_K)$-modules over $B^{\dagger}_{\mathrm{rig},K}$.
\begin{thm}$\label{8}$
Let $D$ be an $E$-$(\varphi,\Gamma_K)$-module over $B^{\dagger}_{\mathrm{rig},K}$. Then $D$ admits unique filtration $0=D_0\subset D_1\subset 
\cdots \subset D_l=D$ such that $D_{i}$ is an $E$-$(\varphi, \Gamma_K)$-module over $B^{\dagger}_{\mathrm{rig},K}$, $D_{i}$ is saturated in $D_{i+1}$ 
and the quotient $D_{i+1}/D_{i}$ is pure of slope $s_{i+1}$ for any $0\le i \le l-1$ with $s_1<s_2<\cdots<s_l$.
Here the slope of an $E$-$(\varphi,\Gamma_K)$-module $M$ is the slope of $M$ 
as $\varphi$-module over $B^{\dagger}_{\mathrm{rig},K}$.
\end{thm}
\begin{proof}
Let $D$ be an $E$-$(\varphi,\Gamma_K)$-module over $B^{\dagger}_{\mathrm{rig},K}$. Then, by the slope filtration theorem (Theorem $\ref{5}$), $D$ admits unique filtration $0=D_0\subset D_1\subset \cdots \subset D_l=D$ by sub $\varphi$-modules over $B^{\dagger}_{\mathrm{rig},K}$ such that
$D_{i}$ is saturated in $D_{i+1}$ and the quotient $D_{i+1}/D_{i}$ is pure of slope $s_{i+1}$ for any $0\le i\le l-1$ with $s_1<s_2<\cdots<s_l$. Then the action of $E$ preserves $D_i$ for any $i$ by uniqueness of filtration. So $D_i$ and $D_{i+1}/D_{i}$ are $E$-$\varphi$-modules over $B^{\dagger}_{\mathrm{rig},K}$ for any $i$. For any $\gamma\in\Gamma_K$, $0\subset \gamma(D_1)\subset\cdots\subset \gamma(D_l)=D$ also satisfies the same conditions as $0\subset D_1\subset \cdots \subset D$ by the commutativity of $\Gamma_K$ and $\varphi$.  
So we get $\gamma(D_i)=D_i$ for any $i$ and $\gamma\in\Gamma_K$ by uniqueness 
of filtration. Hence $D_i$ and $D_{i+1}/D_{i}$ are $E$-$(\varphi,\Gamma_K)$-modules over $B^{\dagger}_{\mathrm{rig}.K}$ for any $i$.
\end{proof}
In particular, it follows from the above theorem that when $D$ is rank one,
 then $D$ is pure of slope $s$ for some $s\in\mathbb{Q}$. 
Concerning this slope $s$, we know more precise information as follows,
 which we need to classify rank one $E$-$B$-pairs.
\begin{lemma}\label{9}
Let $D$ be a rank one $E$-$(\varphi,\Gamma_K)$-module. Then 
$D$ is pure and the slope $\mu(D)$ is contained in $\frac{1}{fe_E}\mathbb{Z}$.
\end{lemma}
\begin{proof}
The claim that $D$ is pure follows from the above theorem.
We prove that $\mu(D)$ is contained in $\frac{1}{fe_E}\mathbb{Z}$ by using the $\Gamma_K$ 
structure on $D$.
First we consider the following short exact sequence of finite groups 
\begin{equation*}
1\rightarrow \mathrm{Gal}(E'_0/K_0)\rightarrow \mathrm{Gal}(E'_0/\mathbb{Q}_p)
\rightarrow \mathrm{Gal}(K_0/\mathbb{Q}_p)\rightarrow 1.
\end{equation*}
Here we put $E'_0:=E\cap K'_0\supseteq K_0$. We put $f'':=[E'_0:K_0]$ and $f':=ff''=[E'_0:\mathbb{Q}_p]$. Then we have a natural surjection
$\Gamma_K\twoheadrightarrow \mathrm{Gal}(K'_0K/K)\isom\mathrm{Gal}(K'_0/K_0)
\twoheadrightarrow\mathrm{Gal}(E'_0/K_0)$. So, if we consider $\varphi$ (which acts on $B^{\dagger}
_{\mathrm{rig},K}$) as an element in $\mathrm{Gal}(E'_0/\mathbb{Q}_p)$, there is a $\gamma\in\Gamma_K$ such that $\varphi^{f}=\gamma$ in $\mathrm{Gal}(E'_0/K_0)$.
On this setting, we consider a rank one $E$-$(\varphi,\Gamma_K)$-module $D$ over $B^{\dagger}
_{\mathrm{rig},K}$. Put $I:=\{\sigma:E'_0\hookrightarrow E\}\isom
\mathrm{Gal}(E'_0/\mathbb{Q}_p)$ and $D_{\sigma}:=D\otimes_{B^{\dagger}_{\mathrm{rig},K}\otimes_{\mathbb{Q}_p}E}B^{\dagger}_{\mathrm{rig},K}\otimes_{E'_0,\sigma}E$ for any $\sigma\in I$.
Because $D$ is rank one, $D_{\sigma}$ is a rank one free $B^{\dagger}_{\mathrm{rig},K}\otimes_{E'_0,\sigma}E$-module. Take a base $e_{\sigma}$ of $D_{\sigma}$ for any $\sigma\in I$. Then $\varphi$ sends $D_{\sigma}$ to 
$D_{\sigma\varphi^{-1}}$ and $\gamma\in \Gamma_K$ sends $D_{\sigma}$ to $D_{\sigma\gamma^{-1}}$. We calculate the slope of $D$ as follows.
Let $e_1$ be a base of $D_{id}$ corresponding to $id\in\mathrm{Gal}(E'_0/\mathbb{Q}_p)$.
Then we have $\varphi^{f}(e_1)=\alpha e_{\varphi^{-f}}$ for some $\alpha\in B^{\dagger}_{\mathrm{rig},K}
\otimes_{E'_0,\varphi^{-f}}E$. Because $\varphi^{f}=\gamma\in \mathrm{Gal}(E'_0/K_0)$ and 
$\gamma$ induces an isomorphism $\gamma:D_{id}\isom D_{\gamma^{-1}}$, there exists some $\beta
\in B^{\dagger}_{\mathrm{rig},K}\otimes_{E'_0,id}E$ such that $\gamma(\beta e_1)=\alpha
e_{\gamma^{-1}}$. Because the actions of $\varphi$ and $\gamma$ commute, we have $\varphi^{2f}(e_1)=\varphi^f(\alpha e_{\varphi^{-f}})=\varphi^f(\alpha e_{\gamma^{-1}})=\varphi^{f}(\gamma(\beta e_1))=\gamma(\varphi^{f}(\beta)\varphi^{f}(e_1))=\gamma(\varphi^{f}(\beta)\gamma(\beta e_1))=\gamma(\varphi^{f}(\beta))\gamma^{2}(\beta)\gamma^2(e_1)$. Repeating this procedure, we get $\varphi^{ff''}(e_1)=\gamma(\varphi^{f(f''-1)}(\beta))\gamma^{2}(\varphi^{f(f''-2)}(\beta))\cdots \gamma^{f''}(\beta)\gamma^{f''}(e_1):=\tilde{\beta}ce_1$ (Here we put $\tilde{\beta}:=\gamma(\varphi^{f(f''-1)}(\beta))\gamma^{2}(\varphi^{f(f''-2)}(\beta))\cdots \gamma^{f''}(\beta)$ and $\gamma^{f''}(e_1):=ce_1$ for some $c\in B^{\dagger}_{\mathrm{rig},K}\otimes_{E'_0, id}E$ because $\gamma^{f''}=id\in \mathrm{Gal}(E'_0/\mathbb{Q}_p)$). Then $D_{id}$ is a rank one $\varphi^{ff''}$-module over 
$B^{\dagger}_{\mathrm{rig},K}\otimes_{E'_0,id}E$ of slope $w_1(\tilde{\beta})+w_1(c)$, where $w_1$ is the valuation of $(B^{\dagger}_{\mathrm{rig},K}\otimes_{E'_0,id}E)^{\mathrm{bd}}$ 
which is the natural extension of the valuation of $(B^{\dagger}_{\mathrm{rig},K})^{\mathrm{bd}}
=B^{\dagger}_{K}$. Because $\varphi$ and $\gamma$ does not change valuation, we have
$w_1(\tilde{\beta})=f''w_1(\beta)\in \frac{f''}{e_E}\mathbb{Z}$. Because $\Gamma_K$ acts continuously on $D$ and $\Gamma_K\subseteq \mathbb{Z}_p^{\times}$, we have $w_1(c)=0$. So the slope of $D_{id}$ as 
$\varphi^{ff''}$-module over $B^{\dagger}_{\mathrm{rig},K}\otimes_{E'_0,id}E$ is contained in $\frac{f''}{e_E}\mathbb{Z}$. We can easily see that 
the slope of $D_{id}$ as $\varphi^{ff''}$-module over $B^{\dagger}_{\mathrm{rig},K}\otimes_{E'_0,id}E$ is same as the slope of $D_{id}$ as $\varphi^{ff''}$-module over $B^{\dagger}_{\mathrm{rig},K}$.
So $D_{id}$ also has slope $w_1(\tilde{\beta})\in \frac{f''}{e_E}\mathbb{Z}$ as $\varphi^{ff''}$-module over  $B^{\dagger}_{\mathrm{rig}
,K}$. Using the Frobenius structure on $D$, 
we can show that $D_{\sigma}$ also has
slope $w_1(\tilde{\beta})$ as $\varphi^{ff''}$-module over $B^{\dagger}_{\mathrm{rig},K}$ for any $\sigma\in I$.
So $D$ also has slope $w_1(\tilde{\beta})$ as $\varphi^{ff''}$-module over $B^{\dagger}_{\mathrm 
{rig},K}$. Because $\mu([ff'']_{\ast}D)=ff''\mu(D)$ by $\cite[\mathrm{Remark}
\, 1.4.5]{Ke}$, $D$ has slope $\frac{w_1(\tilde{\beta})}{ff''}\in \frac{f''}{ff''e_E}\mathbb{Z}=\frac{1}{fe_E}\mathbb{Z}$ as $\varphi$-module over $B^{\dagger}_{\mathrm{rig},K}$. We have finished the proof of this lemma. 
\end{proof}
Next we prove a lemma concerning a relation between tensor products and slopes.
\begin{lemma}$\label{10}$
Let $D_1$, $D_2$ be $E$-$(\varphi, \Gamma_K)$-modules over $B^{\dagger}_{\mathrm{rig},K}$
 which are pure of slope $s_1$, $s_2$ respectively. 
Then $D_1\otimes D_2$ is pure of slope $s_1+s_2$.
\end{lemma}
\begin{proof}
If we decompose $B^{\dagger}_{\mathrm{rig} ,K}\otimes_{\mathbb{Q}_p}E\isom
\oplus_{\sigma:E'_0\hookrightarrow E}B^{\dagger}_{\mathrm{rig},K}\otimes_{E'_0,\sigma}Ee_{\sigma}$, we can decompose $D_1$, $D_2$ into $\sigma$ components,
$D_i\isom\oplus_{\sigma:E'_0\hookrightarrow E}D_{i,\sigma}$ for $i=1,2$, where $D_{i,\sigma}$ is the $B^{\dagger}_{\mathrm{rig},K}\otimes_{E'_0,\sigma}E$ component of 
$D_i$. Then, by the proof of Lemma $\ref{9}$ and by $\cite[\mathrm{Lemma}\, 1.6.3]{Ke}$, we can see that $D_{i,\sigma}$ is a $\varphi^{ff''}$-module over $B^{\dagger}_{\mathrm{rig},K}\otimes_{E'_0,\sigma}E$ which is pure of slope 
$ff''s_i$ for any $\sigma$ and $i=1,2$. So $D_{1,\sigma}\otimes_{B^{\dagger}_{\mathrm{rig},K}\otimes_{E'_0,\sigma}E}D_{2,\sigma}$ is a $\varphi^{ff''}$-module over $B^{\dagger}_{\mathrm{rig},K}\otimes_{E'_0,\sigma}E$ which is pure of slope $ff''(s_1+s_2)$ by $\cite[\mathrm{Corollary}\, 1.6.4]{Ke}$. Because $D_1\otimes D_2\isom \oplus_{\sigma:E'_0\hookrightarrow E}D_{1,\sigma}\otimes_{B^{\dagger}_{\mathrm{rig},K}\otimes_{E'_0,\sigma}E}D_{2,\sigma}$, we can show that $D_1\otimes D_2$ is a $\varphi$-module over $B^{\dagger}_{\mathrm{rig},K}$ which is pure of slope $s_1+s_2$ in the same way as the proof of Lemma $\ref{9}$ and by using $\cite[\mathrm{Lemma}\,1.6.3]{Ke}$.
\end{proof}

\subsection{Equivalence between $B$-pairs and $(\varphi,\Gamma_K)$-modules.}
In this subsection, we recall a result of Berger on the equivalence between the category of $B$-pairs and the category of $(\varphi,\Gamma_K)$-modules over
$B^{\dagger}_{\mathrm{rig},K}$.
First we recall the construction of a functor from $(\varphi,\Gamma_K)$-modules to $B$-pairs ($\cite[2.2]{Be07}$).
Let $D$ be a $(\varphi,\Gamma_K)$-module 
of rank $d$ over $B^{\dagger}_{\mathrm{rig},K}$.
Then Berger showed that $W_e(D):=(\widetilde{B}^{\dagger}_{\mathrm{rig}}[\frac{1}{t}]\otimes_{B^{\dagger}_{\mathrm{rig},K}}D)^{\varphi=1}$ is a free $B_e$-module of rank $d$ ($\cite[\mathrm{Proposition}\, 2.2.6]{Be07}$). It is equipped with a continuous semi-linear $G_K$-action. 
On the other hand, for sufficiently large $r_0>0$, we can take unique $\Gamma_K$-stable finite free $B^{\dagger,r}_{\mathrm{rig},K}$-submodule $D^{r}\subset D$ such that $B^{\dagger}_{\mathrm{rig},K}\otimes_{B^{\dagger,r}_{\mathrm{rig},K}}D^r=D$ and $id_{B^{\dagger,pr}_{\mathrm{rig},K}}\otimes\varphi_D:B^{\dagger,pr}_{\mathrm{rig},K}\otimes_{\varphi,B^{\dagger,r}_{\mathrm{rig},K}}D^r\isom D^{pr}$ for any $r\ge r_0$ ($\cite[\mathrm{Theorem}\, 1.3.3]{Be04}$).
Then Berger showed that the continuous $G_K$-module $W^+_{\mathrm{dR}}(D):=B^+_{\mathrm{dR}}\otimes_{i_n,B^{\dagger,(p-1)p^{n-1}}_{\mathrm{rig},K}}D^{(p-1)p^{n-1}}$ over $B^+_{\mathrm{dR}}$ is independent  of any $n$ such that $(p-1)p^{n-1}\ge r_0$ and showed that there is a canonical $G_K$-isomorphism
$B_{\mathrm{dR}}\otimes_{B_e}W_e(D)\isom B_{\mathrm{dR}}\otimes_{B^+_{\mathrm{dR}}}W^+_{\mathrm{dR}}(D)$ ($\cite[\mathrm{Proposition}\, 2.2.6]{Be07}$).
We put $W(D):=(W_e(D), W^+_{\mathrm{dR}}(D))$. This is a $B$-pair of rank $d$. This defines a functor from the category of $(\varphi,\Gamma_K)$-modules over $B^{\dagger}_{\mathrm{rig},K}$ to the category of $B$-pairs of $G_K$.
\begin{rem}
We can also define an inverse functor from the category of $B$-pairs 
to the category of $(\varphi,\Gamma_K )$-modules ($\cite[2.2]{Be07}$), but 
the definition of this is very difficult.
In this paper, we do not need this construction. So we omit the definition 
of this functor. We denote this inverse functor by $W\mapsto D(W)$.
\end{rem}
\begin{thm}$\label{6}$
The functor $D\mapsto W(D)$ is an exact functor and this gives an equivalence
 of categories between the category of $E$-$(\varphi,\Gamma_K)$-modules over $B^{\dagger}_{\mathrm{rig},K}$ and the category of $E$-$B$-pairs of $G_K$.
\end{thm}
\begin{proof}
For $\mathbb{Q}_p$-coefficient case, it was proved by Berger ($\cite[\mathrm{Theorem}\, 2.2.7]{Be07}$).
Let $D$ be an $E$-$(\varphi,\Gamma_K)$-module over $B^{\dagger}_{\mathrm{rig},K}$, $a\in E$. Then the multiplication by $a$ gives an endomorphism of $D$ as 
$\mathbb{Q}_p$-$(\varphi,\Gamma_K)$-modules. From the functoriality, these 
multiplications give an $E$-action on $W(D)$. So $W(D)$ is an $E$-$B$-pair.
For any $E$-$B$-pair $W$, we can define in similar way the $E$-action 
on $D(W)$. So we get the desired equivalence.
\end{proof}

If we restrict the functor $D\mapsto W(D)$ to $\acute{\mathrm{e}}$tale $E$-$(\varphi,\Gamma_K)$-modules, we get the followng theorem.
\begin{thm}$\label{16}$
The functor $D\mapsto W(D)$ gives an equivalence of categories between the category of $\acute{\mathrm{e}}$tale $E$-$(\varphi,\Gamma_K)$-modules and the category of $E$-$B$-pairs of the form $W(V)$ for some $E$-representation $V$.
\end{thm}
\begin{proof}
$\cite[\mathrm{Proposition}\, 2.2.9]{Be07}$. 
\end{proof}

\begin{defn}
Let $W$ be an $E$-$B$-pair and $s\in\mathbb{Q}$.
Then we say that $W$ is pure of slope $s$ if $D(W)$ is pure of slope $s$.
\end{defn}

\begin{rem}
In fact, we can define directly $B$-pairs with pure slope 
without using $(\varphi,\Gamma_K)$-modules, but we omit this definition in 
this article (see $\cite[3.2]{Be07}$).
\end{rem}

The next theorem is the $E$-$B$-pair version of slope filtration 
theorem.

\begin{thm}$\label{11}$
Let $W$ be an $E$-$B$-pair. Then there exists unique filtration $0=W_0\subset W_1\subset \cdots \subset W_l=W$ by sub $E$-$B$-pairs of $W$ such that $W_{i}$ is  saturated in $W_{i+1}$ and the quotient $W_{i+1}/W_{i}$ is pure of slope $s_i$ 
for any $i$ with $s_1<s_2<\cdots <s_l$.
\end{thm}
\begin{proof}
This follows from Theorem $\ref{8}$ and Theorem $\ref{6}$.
\end{proof}
\begin{lemma}$\label{12}$
Let $W$ be a rank one $E$-$B$-pair.
Then $W$ is pure and its slope is contained in $\frac{1}{fe_E}\mathbb{Z}$.
\end{lemma}
\begin{proof}
This follows from Lemma $\ref{9}$.
\end{proof}

\subsection{Classification of rank one $E$-$B$-pairs.}
In this subsection, we classify rank one $E$-$B$-pairs.
By lemma $\ref{12}$, any rank one $E$-$B$-pairs are pure and their slopes are 
contained in $\frac{1}{fe_E}\mathbb{Z}$. We fix a prime element $\pi_E$ of $E$.
 First we construct a special rank one $E$-$B$-pair $W_0$ which is pure of slope $\frac{1}{fe_E}$.
To construct this, first we define a rank one $E$-filtered $\varphi$-module $D_0$ over $K$, 
which corresponds to $W_0$ by $D^K_{\mathrm{cris}}$. (Here, for any $E$-filtered ($\varphi,N$)-module $D$, we can define 
the rank of $D$ by $\mathrm{rank}(D):=\mathrm{rank}_{K_0\otimes_{\mathbb{Q}_p}E}(D)$ because we can show that $D$ is a free $K_0\otimes_{\mathbb{Q}_p}E$-module in the same way as in the case of $E$-($\varphi,\Gamma_K$)-modules over $B^{\dagger}_{\mathrm{rig},K}$.)
Let $D_{0}:=K_0\otimes_{\mathbb{Q}_p}Ee\isom\oplus_{\sigma:K_0\hookrightarrow E}Ee_{\sigma}$ be a free rank one 
$K_0\otimes_{\mathbb{Q}_p}E$-module. 
We define a filtered $\varphi$-module structure 
on $D_0$ as follows. Define a $K_0\otimes_{\mathbb{Q}_p}E$-semi-linear action of $\varphi$ on $D_0$ by 
$\varphi(e_{id}):=e_{\varphi^{-1}}$, 
$\varphi(e_{\varphi^{-1}}):=e_{\varphi^{-2}}$, $\cdots$, 
$\varphi(e_{\varphi^{-(f-2)}}):=e_{\varphi^{-(f-1)}}$, $\varphi(e_{\varphi^{-(f-1)}}):=
\pi_E e_{id}$, i,e a $\varphi$-module such that $\varphi^f=\pi_E$. 
We define a decreasing filtration on $D_{0, K}:=K\otimes_{K_0}D_0$ by $\mathrm{Fil}^0D_{0, K}:=D_{0, K}$, $\mathrm{Fil}^1D_{0, K}:=0$. Then $D_0$ is a filtered $\varphi$-module over $K$ with a natural $E$-action. 
For any $i\ge 0$, then  $D_0^{\otimes i}$ is the rank one $E$-filtered $\varphi$-module such that $\varphi^f=\pi_E^{i}$ and $\mathrm{Fil}^0D_{0,K}^{\otimes i}=
D_{0,K}^{\otimes i}$ and $\mathrm{Fil}^1D_{0,K}^{\otimes i}=0$. If we put $D_0^{\otimes -1}:=D_0^{\vee}$ the dual of 
$D_0$ defined in the same way as in the case of $E$-$B$-pairs,
then we can see that $D_0^{\otimes -1}$ is an $E$-filtered $\varphi$-module such that $\varphi^f=\pi_E^{-1}$ and $\mathrm{Fil}^0D_{0,K}^{\otimes -1}=
D_{0,K}^{\otimes -1}$ and $\mathrm{Fil}^1D_{0,K}^{\otimes -1}=0$. Moreover, if we define the tensor products of $E$-filtered $(\varphi,N)$-modules in the same way as in the case of $E$-($\varphi,\Gamma_K$)-modules, then we can see that 
$D_0^{\otimes -1}\otimes D_0$ is the trivial $E$-filtered $\varphi$-module. So $D_0^{\otimes i}$ is well defined for any $i\in\mathbb{Z}$
 and satisfies $D_0^{\otimes i}\otimes D_0^{\otimes j}\isom D_0^{\otimes i+j}$ for any $i,j\in\mathbb{Z}$.
 
 By ``weakly admissible imply admissible" theorem for $B$-pairs (Theorem 
 $\ref{4}$ (3)), 
 $W_0:=W(D_0)=(W_e(D_0), W^+_{\mathrm{dR}}(D_0))$ is the rank one cristalline $E$-$B$-pair 
  such that $D^K_{\mathrm{cris}}(W_0)\isom D_0$. 
  \begin{lemma}$\label{13.5}$
  $W_0$ is pure of slope $\frac{1}{fe_E}$.
  \end{lemma}
  \begin{proof}
  Because $W_0$ is rank one, $W_0$ is pure by Lemma $\ref{12}$.
  So it suffices to show that the slope of $W_0$ is equal to $\frac{1}{fe_E}$.
  In $\cite[2.2]{Be04}$, Berger defined a functor $D\mapsto \mathcal{M}(D)$ 
  from the category of filtered $(\varphi, N, G_K)$-modules 
  to the category of (locally trivial) $(\varphi, \Gamma_K)$-modules over $B^{\dagger}_{\mathrm{rig}, K}$. In $\cite[\mathrm{Proposition}\, 2.3.4]{Be07}$, he showed that 
  $\mathcal{M}(D^L_{\mathrm{st}}(W))=D(W)$ for any potentially semi-stable 
  $B$-pair $W$ and for sufficient large $L$. By this compatibility, we have $D(W_0)\isom \mathcal{M}(D_0)$.
  On the other hand, for a rank one filtered $(\varphi, N)$-module $D$, 
  the slope of $\mathcal{M}(D)$ is equal to $t_N(D)-t_H(D)$ by $\cite[\mathrm{Theorem}\, 4.2.1]{Be04}$.
  Here, for a rank one filterd $(\varphi,N)$-module $D:=K_0e$, we define $t_N(D):=\mathrm{val}_p(\varphi(e)/e)$ 
 and define $t_H(D)$ as unique integer $k$ such that $\mathrm{Fil}^k(D_K)/\mathrm{Fil}^{k+1}(D_K)\ne 0$.
 Because $\mu(\wedge^i M)=i\mu(M)$ for a $\varphi$-module $M$ over Robba ring, 
 we have $\mu(D(W_0))=\mu(\mathcal{M}(D_0))=\frac{1}{[E:\mathbb{Q}_p]}\mu(\wedge^{[E:\mathbb{Q}_p]}\mathcal{M}(D_0))=\frac{1}{[E:\mathbb{Q}_p]}(t_N(\wedge^{[E:\mathbb{Q}_p]}D_0)-t_H(\wedge^{[E:\mathbb{Q}_p]}D_0))$. By definition of $D_0$, it is easy to see that $t_{N}(\wedge^{[E:\mathbb{Q}_p]}D_0)=\frac{[E:K_0]}{e_E}$, 
$t_H(\wedge^{[E:\mathbb{Q}_p]}D_0)=0$. So we get $\mu(\mathcal{M}(D_0))=\frac{[E:K_0]}{[E:\mathbb{Q}_p]e_E}=\frac{1}{fe_E}$. So the slope of $W_0$ is $\frac{1}{fe_E}$ by definition of the slopes of $E$-$B$-pairs.
  \end{proof}
By using $W_0$, we can classify all the rank one $E$-$B$-pairs as follows.
Let $\delta:K^{\times}\rightarrow E^{\times}$ be a continuous character where 
$K$ and $E$ are equipped with $p$-adic topology. We put $\delta(\pi_K):=u\pi_E^{i}$ for $u\in \mathcal{O}_E^{\times}$ and $i\in \mathbb{Z}$. We put $\delta_0:K^{\times}\rightarrow \mathcal{O}_E^{\times}$ the unitary 
continuous character such that $\delta_0|_{\mathcal{O}_K^{\times}}=\delta|_{\mathcal{O}_K^{\times}}$ and $\delta_0(\pi_K):=u$. Then by local class field theory, this extends uniquely 
to a continuous character $\tilde{\delta}_0:G_K\rightarrow \mathcal{O}_E^{\times}$ such 
that $\tilde{\delta}_0\circ\mathrm{rec}_K=\delta_0$, here $\mathrm{rec}_K:
K^{\times}\rightarrow G_K^{\mathrm{ab}}$ is the reciprocity map as 
in Notation. Then we put 
$W(\delta):=W(E(\tilde{\delta}_0))\otimes W_0^{\otimes i}$.
By the definitions of tensor products and duals, we can see that these are 
compatible with tensor products and with duals, i.e. $W(\delta_1)\otimes W(\delta_2)\isom W(\delta_1\delta_2)$ and $W(\delta_1)^{\vee}\isom W(\delta_1^{-1})$ for any characters $\delta_1,\delta_2:K^{\times}\rightarrow E^{\times}$. In particular we have $W_0^{\otimes i}\isom W(D_0^{\otimes i})$ for any $i\in\mathbb{Z}$, so $W_0^{\otimes i}$ is pure of slope $\frac{i}{fe_E}$ by Lemma $\ref{10}$ and Lemma $\ref{13.5}$ for any $i\in\mathbb{Z}$.
Then we have the following lemma.
\begin{lemma}$\label{14}$
Let $\delta:K^{\times}\rightarrow E^{\times}$ be a continuous character.
Then the $E$-$B$-pair $W(\delta)$ does not depend on the choice of a 
uniformizer $\pi_E$ of $E$ and on the choice of a uniformizer 
$\pi_K$ of $K$.
\end{lemma}
\begin{proof}
Because $\delta_0$ does not 
depend on the choice of a uniformizer $\pi_K$ of $K$, so 
$W(\delta)$ does not depend on $\pi_K$. 
Let $\pi_E^{'}=v\pi_E$ be another prime of $E$ ($v\in\mathcal{O}_E^{\times}$).
Put $D_0^{'}$ a rank one $E$-filtered $\varphi$-module over $K_0$ defined
 by replacing $\pi_E$ with $\pi_E^{'}$ in the definition of $D_0$. Put 
 $W'_0:=W(D'_0)$ the corresponding $E$-$B$-pair. If we write $\delta(\pi_K)=u\pi_E^i=uv^{-i}{\pi'}^i_E$, then 
 $\delta_0|_{\mathcal{O}_K^{\times}}=\delta'_0|_{\mathcal{O}_K^{\times}}=\delta|_{\mathcal{O}_K^{\times}}$ and $\delta_0(\pi_K)=u$ and $\delta'_0(\pi_K)=uv^{-i}$. Then it suffices 
 to show that $W(E(\tilde{\delta}_0))\otimes W_0^{\otimes i}\isom W(E(\tilde{\delta}'_0))\otimes
 W_0^{' \otimes i}$. Because $\tilde{\delta_0}/\tilde{\delta'_0}$ is a unramified character, it suffices to show that $D^K_{\mathrm{cris}}(E(\tilde{\delta}_0/\tilde{\delta}'_0))\isom 
 D_0^{' \otimes i}\otimes D_0^{\otimes -i}$ as $E$-filtered $\varphi$-modules.
 By calculation, $D^K_{\mathrm{cris}}(E(\tilde{\delta}_0/\tilde{\delta}'_0))=K_0\otimes_{\mathbb{Q}_p}
 Ee$ on which $\varphi$ acts by $\varphi^f(e)=\tilde{\delta_0}/\tilde{\delta'_0}(\pi_K)e=(u/uv^{-i})e
 =v^ie$. By definition, $\varphi$ acts on $D_0^{' \otimes i}\otimes D_0^{\otimes -i}=K_0\otimes_{\mathbb{Q}_p}Ee'$ by $\varphi^f(e')=(\pi'_E/\pi_E)^i e'=v^i e'$. So these are ismorphic 
 as $\varphi$-module. Concerning the filtrations, $\mathrm{Fil}^0D^K_{\mathrm{dR}}(E(\tilde{\delta}_0/\tilde{\delta}'_0))=D^K_{\mathrm{dR}}(E(\tilde{\delta}_0/\tilde{\delta}'_0))$ and $\mathrm{Fil}^1D^K_{\mathrm{dR}}(E(\tilde{\delta}_0/\tilde{\delta}'_0))=0$ because $\tilde{\delta}_0/\tilde{\delta}'_0$ is a unramified 
 character. On the other hand, $\mathrm{Fil}^0D_{0,K}^{' \otimes i}\otimes D_{0,K}^{\otimes -i}=D_{0,K}^{' \otimes i}\otimes D_{0,K}^{\otimes -i}$ and $\mathrm{Fil}^1D_{0,K}^{' \otimes i}\otimes D_{0,K}^{\otimes -i}=0$ by definition of $D_0$ and $D'_0$. So, these are isomorphic as $E$-filtered $\varphi$-modules. 
\end{proof}
\begin{rem}\label{-4}
We can also show in the same way that $W(\delta)$ defined here is isomorphic to $W(\delta)$ defined before Theorem 0.1 in Introduction.
\end{rem}

The classification theorem of rank one $E$-$B$-pairs is the following.
\begin{thm}$\label{15}$
Let $W$ be a rank one $E$-$B$-pair.
Then there exists unique continuous character $\delta:K^{\times}\rightarrow E^{\times}$
such that $W\isom W(\delta)$.
 \end{thm}
\begin{proof}
Let $W$ be a rank one $E$-$B$-pair. Then, by Lemma $\ref{12}$, $W$ is pure of slope $\frac{i}{fe_E}$ for unique $i\in\mathbb{Z}$. Because $W_0$ is pure of slope $\frac{1}{fe_E}$ by lemma $\ref{13.5}$, $W\otimes W_0^{\otimes -i}$ is pure of slope zero by Lemma $\ref{10}$. So, by Theorem $\ref{16}$, there exists unique continuous character $\tilde{\delta}_0:G_K\rightarrow \mathcal{O}_E^{\times}$ such that $W\otimes W_0^{\otimes -i}\isom W(E(\tilde{\delta}_0))$. So $W\isom
  W_0^{\otimes i}\otimes W(E(\tilde{\delta}_0))$. If we define a continuous character
  $\delta:K^{\times}\rightarrow E^{\times}$ such that $\delta(\pi_K):=\tilde{\delta}_0\circ\mathrm{rec}_K(\pi_K)\pi_E^{i}$, $\delta|_{\mathcal{O}_K^{\times}}=
  \tilde{\delta}_0\circ\mathrm{rec}_K|_{\mathcal{O}_K^{\times}}$, then we have 
  $W\isom W(\delta)$ by the definition of $W(\delta)$. The uniqueness of $\delta$ follows from uniqueness of $i$ and $\tilde{\delta}_0$ above.
  \end{proof}

Next we recall some facts about Sen's theory for $B$-pairs to define Hodge-Tate weight of $B$-pairs.
Let $U$ be a $\mathbb{C}_p$-representation of $G_K$, i.e., $U$ is 
a finite dimensional $\mathbb{C}_p$-vector space equipped with a continuous 
semi-linear $G_K$-action. 
Then the union $U^{H_K}_{\mathrm{fini}}$ of finite dimensional sub-$K_{\infty}$-vector spaces of $U^{H_K}$ which are stable by $\Gamma_K$ is the largest 
subspace with this property and satisfies $\mathbb{C}_p\otimes_{K_{\infty}}U^{H_K}_{\mathrm{fini}}\isom U$ (cf. $\cite[2.3]{Be07}$ and the reference there).
Then $U^{H_K}_{\mathrm{fini}}$ is equipped with a $K_{\infty}$-linear 
operator $\nabla_U:=\mathrm{log}(\gamma)/\mathrm{log}(\chi(\gamma))$ for any 
$\gamma\in\Gamma_K$ which is sufficiently close to $1$.
\begin{defn}
Let $W$ be a $B$-pair, then $W^+_{\mathrm{dR}}/tW^+_{\mathrm{dR}}$ is a 
$\mathbb{C}_p$-representation of $G_K$. We put $D_{\mathrm{Sen}}(W):=
(W^+_{\mathrm{dR}}/tW^+_{\mathrm{dR}})^{H_K}_{\mathrm{fini}}$ and put 
$\Theta_{\mathrm{Sen},W}:=\nabla_{W^+_{\mathrm{dR}}/tW^+_{\mathrm{dR}}}$.
Then we say that $W$ is Hodge-Tate if $\Theta_{\mathrm{Sen},W}$ is diagonalizable and all the eigenvalues are contained in $\mathbb{Z}$.
Then we call these eigenvalues the Hodge-Tate weights of $W$.
\end{defn} 

By using this definition, we calculate the Hodge-Tate weights of 
rank one $E$-$B$-pair.
Let $W(\delta)$ be a rank one $E$-$B$-pair for some continuous character 
$\delta:K^{\times}\rightarrow E^{\times}$. If we put $W(\delta)\isom 
W(E(\tilde{\delta}_0))\otimes W_0^{\otimes i}$ as above, then we have $W(\delta)^+_{\mathrm{dR}}\isom B^+_{\mathrm{dR}}\otimes_{\mathbb{Q}_p}E(\tilde{\delta}_0)$ because $W^+_{0,\mathrm{dR}}\isom B^+_{\mathrm{dR}}\otimes_{\mathbb{Q}_p}E$ by the definition of $W_0$. So $W^+_{\mathrm{dR}}(\delta)$ comes from an $E$-representation $E(\tilde{\delta}_0)$. In particular $W^+_{\mathrm{dR}}(\delta)/tW^+_{\mathrm{dR}}(\delta)\isom \mathbb{C}_p\otimes_{\mathbb{Q}_p}E(\tilde{\delta}_0)$ also comes from $E(\tilde{\delta}_0)$. Then, by $\cite[\mathrm{Proposition}\,4.1.2]{Be-Co}$ and the remark after this propositon, $D_{\mathrm{Sen}}(W(\delta))$ is a free $K_{\infty}\otimes_{\mathbb{Q}_p}E$-module of rank one and the operator 
$\Theta_{\mathrm{Sen},W(\delta)}$ acts by multiplication by $w(\delta)$ on $D_{\mathrm{Sen}}(W(\delta))$ for some $w(\delta)\in K\otimes_{\mathbb{Q}_p}E$. If we decompose $w(\delta)$ into $\sigma$-components $K\otimes_{\mathbb{Q}_p}E\isom\oplus_{\sigma:K\hookrightarrow E}Ee_{\sigma}:w(\delta) \mapsto (w(\delta)_{\sigma}e_{\sigma})$. Then $W(\delta)$ is Hodge-Tate $B$-pair if and only if $w(\delta)_{\sigma}$ is contained in $\mathbb{Z}$ for any $\sigma$. 

\begin{defn}\label{-5}
In the above situation, we say the set of numbers $\{w(\delta)_{\sigma}\}_{\sigma}$ the generalized 
Hodge-Tate weight of $W(\delta)$.
\end{defn}

\section{Calculations of cohomologies of $E$-$B$-pairs.}
For classifying two dimensional split trianguline $E$-representations, we need to calculate the extension groups $\mathrm{Ext}^1(W_2, W_1)$ of $E$-$B$-pairs of $W_1$ by $W_2$ for 
rank one $E$-$B$-pairs $W_1$, $W_2$. 
In this section, we define Galois cohomology $\mathrm{H}^i(G_K,W)$ for any $E$-$B$-pair $W$ and interprete $\mathrm{H}^1(G_K,W)$ as 
the extension group $\mathrm{Ext}^1(B_E, W)$. (
Here $B_E:=(B_e\otimes_{\mathbb{Q}_p}E, B^+_{\mathrm{dR}}\otimes_{\mathbb{Q}_p}E)$ is the trivial $E$-$B$-pair.) 
Next we review Liu's results ($\cite{Li}$) which are generalizations 
of Tate duality and Euler-Poincar$\mathrm{\acute{e}}$ characteristic formula to the case of 
($\varphi,\Gamma_K)$-modules over $B^{\dagger}_{\mathrm{rig},K}$.
Finally we calculate some Galois cohomologies of $E$-$B$-pairs which we need for classification of two dimensional split trianguline $E$-representations.

\subsection{Definition of Galois cohomology of $E$-$B$-pairs.}
Let $W:=(W_e, W^+_{\mathrm{dR}})$ be an $E$-$B$-pair. We define 
the complex of $G_K$-modules $C^{\bullet}(W)$ by $\delta_0:C^0(W):=W_e\oplus W^+_{\mathrm{dR}}\rightarrow C^1(W):=W_{\mathrm{dR}}:(x,y)\mapsto x-y$, $C^{i}(W)=0$ for $i\not= 0,1$.

\begin{defn}
Let $W:=(W_e, W^+_{\mathrm{dR}})$ be an $E$-$B$-pair.
Then we define the Galois cohomology of $W$ by 
$\mathrm{H}^i(G_K, W):=\mathrm{H}^i(G_K, C^{\bullet}(W))$.
Here the right hand side is the usual Galois cohomology of a complex
of continuous $G_K$-modules which are computed 
by using continuous cochain complexes.
\end{defn}
By definition, we have the following long exact sequence
\begin{equation*}
\cdots\rightarrow\mathrm{H}^i(G_K,W)\rightarrow\mathrm{H}^i(G_K,W_e)\oplus
\mathrm{H}^i(G_K,W^+_{\mathrm{dR}})\rightarrow\mathrm{H}^i(G_K,W_{\mathrm{dR}})\rightarrow\mathrm{H}^{i+1}(G_K,W)\rightarrow\cdots .
\end{equation*}
From this, we get the following.
\begin{itemize}
\item[(1)]$\mathrm{H}^0(G_K, W)=(W_e\cap W_{\mathrm{dR}}^+)^{G_K}$.
\item[(2)]There exists the exact sequence of $E$-vector spaces
\begin{equation*}
0\rightarrow W_{\mathrm{dR}}^{G_K}/(W_e^{G_K}+W^{+ G_K}_{\mathrm{dR}})\rightarrow\mathrm{H}^1(G_K, W)
\end{equation*}
\begin{equation*}
\rightarrow\mathrm{Ker}(\mathrm{H}^1(G_K, W_e)\oplus\mathrm{H}^1(G_K, W^+_{\mathrm{dR}})\rightarrow\mathrm{H}^1(G_K, W_{\mathrm{dR}}))\rightarrow 0 .
\end{equation*}
\end{itemize}
Next, for any $E$-$B$-pair, we construct a canonical isomorphism\begin{equation*}
 \mathrm{H}^1(G_K,W)\isom
\mathrm{Ext}^1(B_E, W),
\end{equation*}
where the right hand side is the extension group 
 in the category of $E$-$B$-pairs.
For any continuous $G_K$-module $V$, we functorially define the diagram 
$\mathrm{C}^0(V)\xrightarrow[\delta_0]{}\mathrm{C}^1(V)\xrightarrow[\delta_1]{}
\mathrm{C}^2(V)$ as follows:
\begin{itemize}
\item[(1)]
$\mathrm{C}^0(V):=V$, 
\item[]$\mathrm{C}^i(V):=\{f:G_K^{\times i}\rightarrow V \,\mathrm{a}\,\, \mathrm{continuous}\,\,\mathrm{ function}\}$ for $i=1,2$.
\item[(2)]$\delta_0:\mathrm{C}^0(V)\rightarrow \mathrm{C}^1(V):x\mapsto (g\mapsto gx-x)$,
\item[]$\delta_1:\mathrm{C}^1(V)\rightarrow \mathrm{C}^2(V):f\mapsto((g_1,g_2)\mapsto f(g_1g_2)-f(g_1)-g_1f(g_2))$. 
\end{itemize}
Then we have $\mathrm{H}^0(G_K, V)\isom\mathrm{Ker}(\delta_0)$, $\mathrm{H}^1(G_K, V)\isom \mathrm{Ker}(\delta_1)/\mathrm{Im}(\delta_0)$. 
So, for any $E$-$B$-pair $W$, we have $\mathrm{H}^1(G_K, W)\isom \mathrm{Ker}(\tilde{\delta}_1)
/\mathrm{Im}(\tilde{\delta_0})$, here $\tilde{\delta}_0$, $\tilde{\delta}_1$ are defined by 
\begin{itemize}
\item[]$\tilde{\delta}_0:C^0(W_e)\oplus C^0(W^+_{\mathrm{dR}})
\rightarrow C^1(W_e)\oplus C^1(W^+_{\mathrm{dR}})\oplus C^0(W_{\mathrm{dR}}):(x,y)\mapsto
(\delta_0(x), \delta_0(y), x-y)$,
\item[]$\tilde{\delta}_1:C^1(W_e)\oplus C^1(W^+_{\mathrm{dR}})
\oplus C^0(W_{\mathrm{dR}})\rightarrow C^2(W_e)\oplus C^2(W^+_{\mathrm{dR}})\oplus
C^1(W_{\mathrm{dR}}):(f_1,f_2,x)\mapsto(\delta_1(f_1),\delta_1(f_2), f_1-f_2-\delta_0(x))$.
\end{itemize}
By using this expression, we define a map $\mathrm{H}^1(G_K, W)\rightarrow 
\mathrm{Ext}^1(B_E, W)$ as follows. Let $(f_1,f_2,\alpha)\in\mathrm{Ker}(\tilde{\delta}_1)$. Then we define an $E$-$B$-pair $X:=(X_e,X^+_{\mathrm{dR}},\iota)$ as follows:
\begin{itemize}
\item[(1)]$X_e:=W_e\oplus (B_e\otimes_{\mathbb{Q}_p}E)e_{\mathrm{cris}}$ on which $G_K$ acts by $g(x, ae_{\mathrm{cris}}):=(gx+gaf_1(g), gae_{\mathrm{cris}})$
for any $x\in W_e$, $a\in B_e\otimes_{\mathbb{Q}_p}E$ and $g\in G_K$.
\item[(2)]$X^+_{\mathrm{dR}}:=W^+_{\mathrm{dR}}\oplus (B^+_{\mathrm{dR}}\otimes_{\mathbb{Q}_p}E)e_{\mathrm{dR}}$ on which $G_K$ acts by 
$g(y, be_{\mathrm{dR}}):=(gy+gbf_2(g), gbe_{\mathrm{dR}})$ for any $y\in W^+_{\mathrm{dR}}$, $b\in B^+_{\mathrm{dR}}\otimes_{\mathbb{Q}_p}E$ and $g\in G_K$. 
\item[(3)]$\iota:B_{\mathrm{dR}}\otimes_{B_e}X_e\isom B_{\mathrm{dR}}
 \otimes_{B^+_{\mathrm{dR}}}X^+_{\mathrm{dR}}:(x, ae_{\mathrm{cris}})\mapsto 
 (\iota_{W}(x)+a\alpha, ae_{\mathrm{dR}})$ for any $x\in B_{\mathrm{dR}}\otimes_{B_e}W_e$ and $a\in B_{\mathrm{dR}}\otimes_{\mathbb{Q}_p}E$. Here, $\iota_W:B_{\mathrm{dR}}\otimes_{B_e}W_e\isom B_{\mathrm{dR}}\otimes_{B^+_{\mathrm{dR}}}W^+_{\mathrm{dR}}$ is the given isomorphism in the definition of $B$-pair $W=(W_e,W^+_{\mathrm{dR}},
 \iota_W)$.
 \end{itemize}
(Here we see an $E$-$B$-pair $W$ as a triple $W:=(W_e,W^+_{\mathrm{dR}},\iota_W)$
 such that $W_e$ (resp. $W^+_{\mathrm{dR}}$) is finite free over 
$B_e\otimes_{\mathbb{Q}_p}E$ (resp. $B^+_{\mathrm{dR}}\otimes_{\mathbb{Q}_p}E$) with continuous semi-linear $G_K$-actions and $\iota_W:B_{\mathrm{dR}}\otimes_{B_e}W_e\isom B_{\mathrm{dR}}\otimes_{B^+_{\mathrm{dR}}}W^+_{\mathrm{dR}}$ is a $G_K$-equivariant $B_{\mathrm{dR}}\otimes_{\mathbb{Q}_p}E$-semi-linear isomorphism. Then this definition is equivalent to Definition $\ref{-3}$.)
 Because $(f_1, f_2, \alpha)\in \mathrm{Ker}(\tilde{\delta}_1)$, we can easily 
 see that $X$ is a well-defined $E$-$B$-pair which sits in the following short exact sequence 
 \begin{equation*}
 0\rightarrow W \rightarrow X\rightarrow B_E \rightarrow 0.
 \end{equation*}
 So we can see the isomorphism class $[X]$ of the extension $X$ as an element in $\mathrm{Ext}^1(B_E, W)$. By construction, 
 we can easily see that this defines an $E$-linear morphism $\mathrm{Ker}(\tilde{\delta}_1)\rightarrow \mathrm{Ext}^1(B_E, W)$. Then, by standard argument, we see that this morphism factors through 
 $\mathrm{Ker}(\tilde{\delta}_1)/\mathrm{Im}(\tilde{\delta}_0)\rightarrow \mathrm{Ext}^1(B_E, W)$ and that this is in fact $E$-linear isomorphism.
 So we get the following proposition.
\begin{prop}$\label{20}$
Let $W$ be an $E$-$B$-pair.
 Then we have the following functorial $E$-linear isomorphisms.
\begin{itemize}
\item[$\mathrm{(1)}$] $\mathrm{H}^0(G_K,W)\isom (W_e\cap W^+_{\mathrm{dR}})^{G_K}\isom
\mathrm{Hom}(B_{E}, W)$, here $\mathrm{Hom}(W_1, W_2)$ is the $E$-vector space 
of morphisms in the category of $E$-$B$-pairs from $W_1$ to $W_2$ for any $E$-$B$-pairs $W_1$, $W_2$.
\item[$\mathrm{(2)}$]$\mathrm{H}^1(G_K,W)\isom \mathrm{Ext}^1(B_E, W)$.
\end{itemize}
\end{prop}
\begin{proof}
We have already proved (2). For (1), it is easy to see that $(W_e\cap W^+_{\mathrm{dR}})^{G_K}=
\mathrm{Hom}(B_{\mathbb{Q}_p}, W)$ for any $\mathbb{Q}_p$-$B$-pair $W$. Because $\mathrm{Hom}(B_{\mathbb{Q}_p}, W)=\mathrm{Hom}(B_{E}, W)$ for any $E$-$B$-pair $W$ (here on the left hand side, we see $W$ as $\mathbb{Q}_p$-$B$-pair), so 
we get (1).
\end{proof}
\begin{rem}
In the case of $W=W(V)$ for an $E$-representation $V$ of $G_K$,
we have $\mathrm{H}^i(G_K, V)\isom \mathrm{H}^1(G_K, W(V))$: indeed, 
by the fundamental short exact sequence $0\rightarrow \mathbb{Q}_p
\rightarrow B_e\oplus B^+_{\mathrm{dR}}\rightarrow B_{\mathrm{dR}}\rightarrow 0$, we have a quasi-isomorphism
 $V[0]\isom C^{\bullet}(W(V))$.
 \end{rem}

In the application to classification of two dimensional 
potentially semi-stable 
split trianguline $E$-representations, we need to know 
when an extension is cristalline, semi-stable or de Rham $E$-$B$-pair.
So, as in the case of usual Galois cohomology of $p$-adic representations, 
we define Bloch-Kato's cohomologies $\mathrm{H}^1_e(G_K, W)$, $\mathrm{H}^1_f(G_K, W)$, and $\mathrm{H}^1_g(G_K, W)$ ($\cite[3]{BK}$) as follows.
\begin{defn}$\label{d}$
Let $W:=(W_e, W^+_{\mathrm{dR}})$ be an $E$-$B$-pair.
Then we define 
\begin{itemize}
\item[]$\mathrm{H}^1_e(G_K, W):=\mathrm{Ker}(\mathrm{H}^1(G_K,W)\rightarrow \mathrm{H}^1(G_K, W_e))$,
\item[] $\mathrm{H}^1_f(G_K,W):=\mathrm{Ker}
(\mathrm{H}^1(G_K,W)\rightarrow \mathrm{H}^1(G_K, B_{\mathrm{cris}}\otimes_{B_e}W_e))$, 
\item[]$\mathrm{H}^1_g(G_K, W):=\mathrm{Ker}(\mathrm{H}^1(G_K,W)\rightarrow
\mathrm{H}^1(G_K, W_{\mathrm{dR}}))$.
\end{itemize}
Here the above maps are induced by the natural maps $C^{\bullet}(W)\rightarrow W_e
\rightarrow B_{\mathrm{cris}}\otimes_{B_e}W_e\rightarrow B_{\mathrm{dR}}\otimes_{B_e}W_e\isom W_{\mathrm{dR}}$, so we have natural injections $\mathrm{H}^1_e(G_K, W)\subseteq \mathrm{H}^1_f(G_K, W) \subseteq \mathrm{H}^1_g(G_K, W)$.
\end{defn}
\begin{rem}$\label{b}$
If $W$ is a cristalline (resp. de Rham) $E$-$B$-pair, then $[X]\in \mathrm{Ext}^1(B_E, W)\isom \mathrm{H}^1(G_K, W)$ is in $\mathrm{H}^1_f(G_K, W)$ (resp. in $\mathrm{H}^1_g(G_K, W)$) if and only if $X$ is a cristalline (resp. de Rham) $E$-$B$-pair.
\end{rem}
As in the case of usual $p$-adic representations, we have a 
dimension formula of $\mathrm{H}^1_f(G_K, W)$ and $\mathrm{H}^1_e(G_K, W)$.
Before stating this, we prove the following lemma.
\begin{lemma}$\label{21}$
Let $W$ be a de Rham $E$-$B$-pair.
Then the canonical map $\mathrm{H}^1(G_K, W^+_{\mathrm{dR}})\rightarrow 
\mathrm{H}^1(G_K, W_{\mathrm{dR}})$ is injective.
\end{lemma}
\begin{proof}
The proof is the same as that in the case of $p$-adic representation 
($\cite[\mathrm{Lemma}\, 3.8.1]{BK}$), but we give 
the proof for the convenience of readers.
Consider the following short exact sequence
\begin{equation*}
0\rightarrow W_{\mathrm{dR}}^{+ G_K}\rightarrow W_{\mathrm{dR}}^{G_K}\rightarrow (W_{\mathrm{dR}}/W^+_{\mathrm{dR}})^{G_K}.
\end{equation*}
From this we have 
\begin{align*}
\mathrm{dim}_E(D_{\mathrm{dR}}(W))&\le 
\mathrm{dim}_E(W_{\mathrm{dR}}^{+ G_K})+\mathrm{dim}_E((W_{\mathrm{dR}}
/W^+_{\mathrm{dR}})^{G_K}) \\
 &\le \sum_{i\in\mathbb{Z}}\mathrm{dim}_E(\mathbb{C}_p(i)\otimes_{\mathbb{C}_p}(W^+_{\mathrm{dR}}/tW^+_{\mathrm{dR}}))^{G_K} \\
 &=\mathrm{rank}(W)=
\mathrm{dim}_E(D_{\mathrm{dR}}(W)).
\end{align*}
 Here we use the fact that 
de Rham $B$-pair is Hodge-Tate and $W$ is de Rham. So the map $W_{\mathrm{dR}}^{G_K}\rightarrow (W_{\mathrm{dR}}/W_{\mathrm{dR}}^+)^{G_K}$ is surjective. Hence 
the natural map $\mathrm{H}^1(G_K, W^+_{\mathrm{dR}})\rightarrow\mathrm{H}^1(G_K, W_{\mathrm{dR}})$ is injective.
\end{proof}
The dimension formulas are as follows.
\begin{prop}$\label{22}$
Let $W$ be a de Rham $E$-$B$-pair. Then we have 
\begin{itemize}
\item[]$\mathrm{dim}_E(\mathrm{H}^1_f(G_K, W))=\mathrm{dim}_E(D^K_{\mathrm{dR}}(W)/\mathrm{Fil}^0D^K_{\mathrm{dR}}(W))+\mathrm{dim}_E(\mathrm{H}^0(G_K, W))$,
\item[]$\mathrm{dim}_E(\mathrm{H}^1_f(G_K, W)/\mathrm{H}^1_e(G_K, W))=\mathrm{dim}_E(D^K_{\mathrm{cris}}(W)/(1-\varphi)D^K_{\mathrm{cris}}(W))$.
\end{itemize}
\end{prop}
\begin{proof}
This proof is essentially the same as that of $\cite[\mathrm{Corollary}\, 3.8.4]{BK}$.
First we consider the following short exact sequence
\begin{equation*}
0\rightarrow B_e\rightarrow B_{\mathrm{cris}}\xrightarrow[]{1-\varphi} B_{\mathrm{cris}}\rightarrow 0.
\end{equation*}
Tensoring by $W_e$ over $B_e$, we can easily see that $C^{\bullet}(W)$ is naturally quasi-isomorphic to $C^{'\bullet}(W):=(W_{\mathrm{cris}}\oplus W^+_{\mathrm{dR}}
\rightarrow W_{\mathrm{cris}}\oplus W_{\mathrm{dR}}:(x,y)\mapsto ((1-\varphi)x,x-y))$ (Here we put $W_{\mathrm{cris}}:=B_{\mathrm{cris}}\otimes_{B_e}W_e$). So, by definition of $\mathrm{H}^1_f$ and by the above lemma, we have the 
following two exact sequences
\begin{equation*}
0\rightarrow \mathrm{H}^0(G_K, W)\rightarrow D^K_{\mathrm{cris}}(W)\oplus
\mathrm{Fil}^0D^K_{\mathrm{dR}}(W)\rightarrow
D^K_{\mathrm{cris}}(W)\oplus D^K_{\mathrm{dR}}(W)\rightarrow \mathrm{H}^1_f(G_K, W)
\rightarrow 0,
\end{equation*}
\begin{equation*}
0\rightarrow \mathrm{H}^0(G_K, W)\rightarrow D^K_{\mathrm{cris}}(W)^{\varphi=1}
\oplus \mathrm{Fil}^0D^K_{\mathrm{dR}}(W)\rightarrow D^K_{\mathrm{dR}}(W)
\rightarrow \mathrm{H}^1_e(G_K, W)\rightarrow 0.
\end{equation*}
From these, we get the desired formulas.
\end{proof}

\subsection{Euler-Poincar$\acute{\mathrm{e}}$ characteristic formula 
and Tate local duality for 
$B$-pairs.}
In this subsection, we review Liu's results genaralizing some fundamental 
results of Tate on Galois cohomology of $p$-adic representations, following 
$\cite{Ke07}$ and $\cite{Li}$. 
Liu constructed a category of $B$-quotients, which is a minimal abelian category 
in which the category of $B$-pairs is contained as full subcategory.
Then he defined the Galois cohomology of $B$-quotients as the univerasl 
$\delta$-functor of $\mathrm{H}_{\mathrm{Liu}}^0(W):=\mathrm{Hom}(B_{\mathbb{Q}_p}, W)$ for any $B$-quotient $W$ (here $\mathrm{Hom}$ is the set of morphisms 
of $B$-quotients). In this paper, we denote this cohomology as $\mathrm{H}^i_{\mathrm{Liu}}(G_K, W)$ for 
a $B$-quotient $W$. Then it is shown in $\cite{Ke07}$ that there exists the isomorphism $H^1_{\mathrm{Liu}}(G_K, W)\isom \mathrm{Ext}^1(B_{\mathbb{Q}_p}, W)$ for any $B$-pair. So, for any $B$-pair $W$, we have isomorphisms $\mathrm{H}^1(G_K 
, W)\isom\mathrm{H}^1_{\mathrm{Liu}}(G_K, W)$ and $\mathrm{H}^0(G_K, W)\isom 
\mathrm{H}^0_{\mathrm{Liu}}(G_K, W)$. For any $\mathbb{Q}_p$-representation $V$ of 
$G_K$, it is shown in $\cite{Ke07}$ that there is a functorial isomorphism $\mathrm{H}^i(G_K, V)
\isom \mathrm{H}^i_{\mathrm{Liu}}(G_K, W(V))$ for any $i\in \mathbb{N}$. In this paper, 
we review the results only for $B$-pairs.
\begin{thm}$\label{23}$
Let $W$ be a $B$-pair. Then 
\begin{itemize}
\item[$\mathrm{(1)}$]For i=0,1,2, $\mathrm{H}^i_{\mathrm{Liu}}(G_K, W)$ is finite dimensional over $\mathbb{Q}_p$ and $\mathrm{H}^i_{\mathrm{Liu}}(G_K, W)=0$ for $i\not= 0,1,2$.
\item[$\mathrm{(2)}$] We have $\sum_{i=0}^{2}(-1)^i\mathrm{dim}_{\mathbb{Q}_p}\mathrm{H}^i_{\mathrm{Liu}}(G_K, W)=-[K:\mathbb{Q}_p]\mathrm{rank}_{B_e}(W_e)$.
\item[$\mathrm{(3)}$]There is a natural perfect pairing
\begin{equation*}
\mathrm{H}^i_{\mathrm{Liu}}(G_K, W)\times\mathrm{H}^{2-i}_{\mathrm{Liu}}(G_K, W^{\vee}(\chi))\rightarrow 
\mathrm{H}^2_{\mathrm{Liu}}(G_K, W\otimes W^{\vee}(\chi))\rightarrow 
\mathrm{H}^2_{\mathrm{Liu}}(G_K, W(\mathbb{Q}_p(\chi))\isom \mathbb{Q}_p.
\end{equation*}
Here $\chi:G_K\rightarrow \mathbb{Z}_p^{\times}\hookrightarrow \mathcal{O}_E^{\times}$ is the $p$-adic cyclotomic character and the last isomorphism is the one which appears in usual Tate duality.
\end{itemize}
\end{thm}
\begin{proof}
$\cite[\mathrm{Theoem}\, 8.1]{Ke07}$, $\cite[\mathrm{Theorem}4.3.]{Li}$, $\cite[\mathrm{Theorem}4.7.]{Li}$.
\end{proof}
\begin{rem}
The above perfect pairing is defined by using $(\varphi,\Gamma_K)$-modules $\cite[4.2]{Li}$) and the equivalence between $B$-pairs and $(\varphi,\Gamma_K)$-modules ($\cite[p.8]{Ke07}$). On the other hand, for $\mathbb{Q}_p$-$B$-pair $W$, 
we can define a pairing 
\begin{itemize}
 \item[]$\mathrm{H}^1(G_K, W)\times\mathrm{H}^1(G_K, W^{\vee}(\chi))\rightarrow 
\mathrm{H}^2(G_K, W\otimes W^{\vee}(\chi))\rightarrow 
\mathrm{H}^2(G_K, W(\mathbb{Q}_p(\chi))\isom \mathbb{Q}_p$,
\end{itemize}
in the usual way by using cocycles. If we identify $\mathrm{H}^1(G_K, W^{\vee}(\chi))\isom \mathrm{H}^1_{\mathrm{Liu}}(G_K, W^{\vee}(\chi))$ as above, then 
we can check that the both pairings are same by using the fact that 
$\mathrm{H}^i_{\mathrm{Liu}}(G_K, W)$ is the universal $\delta$ functor ($\cite
[\mathrm{Theorem}\, 8.1]{Ke07}$ )
and $\mathrm{H}^i(G_K, W)$ is a $\delta$-functor.
\end{rem}
Let $W$ be an $E$-$B$-pair. Then by (2) above and the remark before the above theorem, we have 
\begin{align*}
\mathrm{dim}_{\mathbb{Q}_p}\mathrm{H}^1(G_K, W) &=\mathrm{dim}_{\mathbb{Q}_p}\mathrm{H}^1_{\mathrm{Liu}}(G_K, W) \\
&=[K:\mathbb{Q}_p][E:\mathbb{Q}_p]\mathrm{rank}(W)+ \mathrm{dim}_{\mathbb{Q}_p}
\mathrm{H}^0_{\mathrm{Liu}}(G_K, W)+\mathrm{dim}_{\mathbb{Q}_p}
\mathrm{H}^2_{\mathrm{Liu}}(G_K, W) \\
 &=[K:\mathbb{Q}_p][E:\mathbb{Q}_p]\mathrm{rank}(W)+ \mathrm{dim}_{\mathbb{Q}_p}
\mathrm{H}^0_{\mathrm{Liu}}(G_K, W)+\mathrm{dim}_{\mathbb{Q}_p}
\mathrm{H}^0_{\mathrm{Liu}}(G_K, W^{\vee}(\chi)) \\
 &=[K:\mathbb{Q}_p][E:\mathbb{Q}_p]\mathrm{rank}(W)+ \mathrm{dim}_{\mathbb{Q}_p}
\mathrm{H}^0(G_K, W)+\mathrm{dim}_{\mathbb{Q}_p}
\mathrm{H}^0(G_K, W^{\vee}(\chi)). 
\end{align*}
Dividing this equality by $[E:\mathbb{Q}_p]$, we get the following dimension formula. 
\begin{prop}$\label{24}$
Let $W$ be an $E$-$B$-pair.
Then we have 
\begin{itemize}
\item[]$\mathrm{dim}_{E}\mathrm{H}^1(G_K, W)
=[K:\mathbb{Q}_p]\mathrm{rank}(W)+ \mathrm{dim}_{E}
\mathrm{H}^0(G_K, W)+\mathrm{dim}_{E}
\mathrm{H}^0(G_K, W^{\vee}(\chi))$.
\end{itemize}
\end{prop}
As in the case of $p$-adic representations, we have dualities between $\mathrm{H}^1_e$, $\mathrm{H}^1_f$ and $\mathrm{H}^1_g$.
\begin{prop}$\label{25}$
Let $W$ be a $B$-pair.
In the above perfect pairing in Theorem $\ref{23}$ $(3)$, $\mathrm{H}^1_e(G_K, W)$ and $\mathrm{H}^1_g(G_K, W^{\vee}(\chi))$ are the exact annihilators of each 
other. The same statement holds with $e$ replaced by $g$ and $g$ by $e$ 
and also when $e$ and $g$ are both replaced by $f$.
\end{prop}
\begin{proof}
This proof is also essentially same as that in the case of $p$-adic 
representations (cf. $\cite[\mathrm{Proposition}\, 3.8]{BK}$). 
\end{proof}

\subsection{Calculations of cohomologies of rank one $E$-$B$-pairs.}
To classify rank two split trianguline $E$-$B$-pairs, we need to calculate the dimension of 
 $\mathrm{Ext}^1(W(\delta_2), \allowbreak W(\delta_1))$ for 
 any continuous characters $\delta_1,\delta_2:K^{\times}\rightarrow E^{\times}$ .
 Twisting by $W(\delta_2^{-1})$, it is isomorphic to $\mathrm{Ext}^1(B_E, W(\delta_1/\delta_2)))\isom \mathrm{H}^1(G_K, W(\delta_1/\delta_2))$. 
 We will calculate this in the following ways.
 \begin{lemma}$\label{26}$
 For any embedding $\sigma:K\hookrightarrow E$, we define a continuous 
character  $\sigma(x):K^{\times}\rightarrow E^{\times}:y\mapsto \sigma(y)$.
 Then for any $\{k_{\sigma}\}_{\sigma}$ $(k_{\sigma}\in\mathbb{Z}$ for any $\sigma)$, we have an isomorphism of $E$-$B$-pairs
\begin{itemize}
\item[] $W(\prod_{\sigma:K\hookrightarrow E}\sigma(x)^{k_{\sigma}})\isom (B_e\otimes_{\mathbb{Q}_p}E, \oplus_{\sigma:K\hookrightarrow E}t^{k_{\sigma}}B^+_{\mathrm{dR}}\otimes_{K,\sigma}E)$.
\end{itemize}
 \end{lemma}
 \begin{proof}
 First we prove that there is an isomorphism 
 \begin{itemize}
 \item[]$W(id(x))\isom (B_e\otimes_{\mathbb{Q}_p}E, tB^+_{\mathrm{dR}}\otimes_{
 K,id}E\oplus \oplus_{\sigma:K\hookrightarrow E,\not= id}B^+_{\mathrm{dR}}\otimes_{K,\sigma}E)$, 
 \end{itemize}
 here $id:K\hookrightarrow E$ is the given embedding.
 We decompose $id(x):K^{\times}\rightarrow E^{\times}$ into $id(x):=\delta_0\delta_1$, here 
 $\delta_0|_{\mathcal{O}_K^{\times}}:=id|_{\mathcal{O}_K^{\times}}, \delta_0(\pi_K):=1$ and 
 $\delta_1|_{\mathcal{O}_K^{\times}}$ is trivial and $\delta_1(\pi_K):=\pi_K$. Then, 
 by definition of the reciprocity map $\mathrm{rec_K}:K^{\times}\rightarrow 
 G_K^{\mathrm{ab}}$ in Notation, we have 
$\tilde{\delta_0}=\chi_{\mathrm{LT}}$, i.e. the Lubin-Tate character associated to $\pi_K$. So we have 
 $W(\delta_0)=W(E(\chi_{\mathrm{LT}}))$. Then $D^K_{\mathrm{cris}}(W(\delta_0)):=K_0\otimes
 _{\mathbb{Q}_p}Ee$ is the filtered $\varphi$-module such that $\varphi$ acts by $\varphi^f(e)=\pi_K^{-1}e$ and the filtration on 
 $K\otimes_{K_0}D^K_{\mathrm{cris}}(W(\delta_0))=K\otimes_{\mathbb{Q}_p}Ee\isom
 \oplus_{\sigma:K \hookrightarrow E}Ee_{\sigma}$ is defined by 
 $\mathrm{Fil}^{-1}=K\otimes_{K_0}D^K_{\mathrm{cris}}(W(\delta_0))$ , $\mathrm{Fil}^0=
 \oplus_{\sigma\not= id}Ee_{\sigma}$ and $\mathrm{Fil}^1=0$ by 
 $\cite{Co02}$ Proposition 9.10 and Lemma 9.18.
 On the other hand, $D^K_{\mathrm{cris}}(W(\delta_1)):=K_0\otimes_{\mathbb{Q}_p}Ee'$ is the 
 filtered $\varphi$-module such that $\varphi^f(e')=\pi_K e'$ and $\mathrm{Fil}^0(K\otimes_{K_0}
 D^K_{\mathrm{cris}}(W(\delta_1)))=K\otimes_{K_0}D^K_{\mathrm{cris}}(W(\delta_1)), \mathrm{Fil}^1=0$.
 So we have $D^K_{\mathrm{cris}}(id(x))\isom D^K_{\mathrm{cris}}(\delta_0)\otimes 
 D^K_{\mathrm{cris}}(\delta_1):=K_0\otimes_{\mathbb{Q}_p}Ee''$ on which $\varphi$ acts by 
 $\varphi^f(e'')=e''$ and the filtration on $K\otimes_{K_0}D^K_{\mathrm{cris}}(id(x))\isom
  \oplus_{\sigma:K \hookrightarrow E}Ee''_{\sigma}$ is given by $\mathrm{Fil}^{-1}=K\otimes_{K_0}D^K_{\mathrm{cris}}(id(x))$, $\mathrm{Fil}^{0}=\oplus_{\sigma\not= id}Ee''_{\sigma}$ and $\mathrm{Fil}^1=0$. 
 Then, by the definition of $D^K_{\mathrm{cris}}$ for $E$-$B$-pairs, it is easy to see that 
  $D^K_{\mathrm{cris}}((B_e\otimes_{\mathbb{Q}_p}E, tB^+_{\mathrm{dR}}\otimes_{K,id}E
  \oplus \oplus_{\sigma\not= id} B^+_{\mathrm{dR}}\otimes_{K,\sigma}E))\isom 
  D^K_{\mathrm{cris}}(id(x))$ as filtered $\varphi$-modules. So, in this case, we have proved the lemma.
  In the case where $\prod_{\sigma:K\hookrightarrow E}\sigma(x)^{k_{\sigma}}=\sigma(x)$ for some $\sigma$, we can prove the lemma in the same way. By tensoring these, we can prove the lemma for any $\prod_{\sigma:K\hookrightarrow E}\sigma(x)^{k_{\sigma}}$.
  \end{proof}
  Let $N_{K/\mathbb{Q}_p}(x):K^{\times}\rightarrow \mathbb{Q}_p^{\times}$ be the norm map and let $|-|:\mathbb{Q}_p^{\times}\rightarrow E^{\times}$ be the absolute value character such that 
  $|p|:=\frac{1}{p}, |u|:=1$ for any $u\in \mathbb{Z}_p^{\times}$. 
  
  \begin{lemma}$\label{27}$
  There is an isomorphism $W(E(\chi))\isom W(N_{K/\mathbb{Q}_p}(x) |N_{K/\mathbb{Q}_p}(x)|)$.
  \end{lemma}
  \begin{proof}
When $K=\mathbb{Q}_p$, $\chi:G_{\mathbb{Q}_p}\rightarrow \mathbb{Z}_p^{\times}$ satisfies $x|x|=\chi\circ\mathrm{rec}_{\mathbb{Q}_p}$ by local class field 
theory. So in general case, $\chi:G_K\rightarrow \mathcal{O}_E^{\times}$ corresponds to $N_{K/\mathbb{Q}_p}(x) |N_{K/\mathbb{Q}_p}(x)|$ by local class field theory.
  \end{proof}

The next proposition is a generalization of $\cite[\mathrm{Proposition}\,3.1]
{Co07a}$.
 \begin{prop}$\label{28}$
 Let $\delta:K^{\times}\rightarrow E^{\times}$ be a continuous character.
 Then $\mathrm{H}^0(G_K, W(\delta))\isom E$ if and only if 
 $\delta=\prod_{\sigma:K\hookrightarrow E}\sigma(x)^{k_{\sigma}}$ such that 
 $k_{\sigma}\le0$ for any $\sigma$.
 Otherwise $\mathrm{H}^0(G_K, W(\delta))=0$.
 \end{prop}
 \begin{proof}
 If we assume that $\mathrm{H}^0(G_K, W(\delta))\not= 0$, then there is a non zero morphism 
 $f:B_E\rightarrow W(\delta)$ of $E$-$B$-pairs because $\mathrm{H}^0(G_K, W(\delta))
 =\mathrm{Hom}(B_E, W(\delta))$. Then by Lemma $\ref{2}$, $\mathrm{Ker}(f)$
 and $\mathrm{Im}(f)$ are also $E$-$B$-pairs. Because $B_E$ is a rank one $E$-$B$-pair,
 one of $\mathrm{Ker}(f)$ and $\mathrm{Im}(f)$ must be zero. Because $f\not= 0$ we have $\mathrm{Im}(f)
 \not= 0$ so we have $\mathrm{Ker}(f)=0$, i.e. $f$ must be injective.
 So we have an injection $f:B_e\otimes_{\mathbb{Q}_p}E\hookrightarrow
 W(\delta)_e$ of free $B_e\otimes_{\mathbb{Q}_p}E$-modules of the same rank. By Lemma $\ref{1}$, the cokernel is also a free $B_e\otimes_{\mathbb{Q}_p}E$-module. So the cokernel must be zero,  so $f:B_e\otimes_{\mathbb{Q}_p}E\isom W(\delta)_e$ must be an isomorphism.
 Then $W(\delta)^+_{\mathrm{dR}}$ is a $G_K$-stable $B^+_{\mathrm{dR}}$-lattice in $B_{\mathrm{dR}}\otimes_{\mathbb{Q}_p}E$ with $E$-action which contains $B^+_{\mathrm{dR}}
 \otimes_{\mathbb{Q}_p}E=\oplus_{\sigma:K\hookrightarrow E}B^+_{\mathrm{dR}}\otimes_{K,\sigma}E$. So $W(\delta)^+_{\mathrm{dR}}$ must be of the form $W(\delta)^+_{\mathrm{dR}}=
 \oplus_{\sigma:K\hookrightarrow E}t^{k_{\sigma}}B^+_{\mathrm{dR}}\otimes_{K,\sigma}E$ for 
 some $k_{\sigma}\le 0$ for any $\sigma$. So $W(\delta)=W(\prod_{\sigma:K\rightarrow E}\sigma(x)^{k_{\sigma}})$ by Lemma $\ref{26}$. In this case, it is clear that $\mathrm{Hom}
 (B_E, W(\delta))=E$, so we have proved the proposition.
 \end{proof}
 
By this proposition and Liu's Euler-Poincar$\acute{\mathrm{e}}$ formula, we can calculate 
the dimension of $\mathrm{Ext}^1(W(\delta_2), W(\delta_1))$ as follows. This is a generalization of $\cite[\mathrm{Theorem}3.9]{Co07a}$. 

\begin{prop}$\label{29}$
Let $\delta_1,\delta_2:K^{\times}\rightarrow E^{\times}$ be continuous characters. Then  $\mathrm{dim}_E\mathrm{Ext}^1(W(\delta_2), \allowbreak W(\delta_1))$ is equal to
\begin{itemize}
\item[$\mathrm{(1)}$]
$[K:\mathbb{Q}_p]+1$ when $\delta_1/\delta_2= \prod_{\sigma:K\hookrightarrow E}\sigma(x)^{k_{\sigma}}$ such that 
$k_{\sigma}\le 0$ for any $\sigma$,
\item[$\mathrm{(2)}$]
$[K:\mathbb{Q}_p]+1$ when $\delta_1/\delta_2= |N_{K/\mathbb{Q}_p}(x)|
\prod_{\sigma:K\hookrightarrow E}\sigma(x)^{k_{\sigma}}$ such that $k_{\sigma}\ge 1$ for any $\sigma$,
\item[$\mathrm{(3)}$]
$[K:\mathbb{Q}_p]$ otherwise.
\end{itemize}
\end{prop}
 \begin{proof}
By Lemma $\ref{27}$ and Proposition $\ref{28}$, we have
\begin{itemize}
\item[]$\mathrm{H}^0(G_K, W(\delta))=E$ if and only if 
$\delta=\prod_{\sigma:K\hookrightarrow E}\sigma(x)^{k_{\sigma}}$ such that 
$k_{\sigma}\le 0$ for any $\sigma$,
\item[]$\mathrm{H}^0(G_K, W(\delta)^{\vee}(\chi))\isom
\mathrm{H}^0(G_K, W(\delta^{-1}N_{K/\mathbb{Q}_p}(x)|N_{K/\mathbb{Q}_p}(x)|))= E$ if and only if $\delta= N_{K/\mathbb{Q}_p}(x)|N_{K/\mathbb{Q}_p}(x)|
\prod_{\sigma:K\hookrightarrow E}\sigma(x)^{k_{\sigma}}$ such that $k_{\sigma}\ge 0$ for any $\sigma$.
\end{itemize}
So we have the desired result by Corollary $\ref{24}$.
\end{proof}

To classify two dimensional split trianguline $E$-representations, we need to know which extension class in $\mathrm{Ext}^1
(W(\delta_2), W(\delta_1))$ is of the form $[W(V)]$ for some two dimensional $E$-representation $V$.
For this problem, the next exact sequence is very important, which is the $B$-pair analogue of the map $i_k$ defined in $\cite[3.7]{Co07a}$.
Let $\delta:K^{\times}\rightarrow E^{\times}$ be a continuous character and let $\{k_{\sigma}\}_{\sigma:K\hookrightarrow E}$ 
be $k_{\sigma}\in\mathbb{Z}_{\ge 0}$ for any $\sigma$.  We consider the natural inclusion of $E$-$B$-pairs $W(\delta)\hookrightarrow W(\delta)\otimes W(\prod_{\sigma:K\hookrightarrow E}\sigma(x)^{-k_{\sigma}})=W(\delta\prod_{\sigma:K\hookrightarrow E}\sigma(x)^{-k_{\sigma}})$ which is obtained by tensoring with $W(\delta)$ the natural inclusion $B_E\hookrightarrow W(\prod_{\sigma:K\hookrightarrow E}\sigma(x)^{-k_{\sigma}})$. Then, by Lemma $\ref{26}$, we have 
\begin{itemize}
\item[]$W^+_{\mathrm{dR}}(\delta\prod_{\sigma:K\hookrightarrow E}\sigma(x)^{-k_{\sigma}})=W^+_{\mathrm{dR}}(\delta)\otimes_{B^+_{\mathrm{dR}}\otimes_{\mathbb{Q}_p}E}(\oplus_{\sigma:K\hookrightarrow E} t^{-k_{\sigma}}B^+_{\mathrm{dR}}\otimes_{K,\sigma}E)$. 
\end{itemize}
As for $B_e$-part, 
we can prove in the same way as the proof of Proposition $\ref{28}$ that $W_e(\delta)\isom W_e(\delta\prod_{\sigma:K\hookrightarrow E}\sigma(x)^{-k_{\sigma}})$. So we have the following short exact 
sequence of complexes of $G_K$-modules
\begin{equation*}
0\rightarrow C^{\bullet}(W(\delta))\rightarrow C^{\bullet}(W(\delta\prod_{\sigma:K\hookrightarrow E}\sigma(x)^{-k_{\sigma}}))\rightarrow \oplus_{\sigma:K\hookrightarrow E}
t^{-k_{\sigma}}W^+_{\mathrm{dR}}(\delta)_{\sigma}/W^+_{\mathrm{dR}}(\delta)_{\sigma}[0]\rightarrow 0,
\end{equation*}
 where $W^+_{\mathrm{dR}}(\delta)_{\sigma}:=W^+_{\mathrm{dR}}(\delta)\otimes_{B^+_{\mathrm{dR}}\otimes_{\mathbb{Q}_p}E}(B^+_{\mathrm{dR}}\otimes_{K,\sigma}E)$ for any $\sigma$ where the last tensor product is taken by the projection to the $\sigma$-component of $B^+_{\mathrm{dR}}\otimes_{\mathbb{Q}_p}E\isom\oplus_{\sigma:K\hookrightarrow E}B^+_{\mathrm{dR}}\otimes_{K,\sigma}E$.
From this short exact sequence, we get the following long exact sequence
\begin{equation*}
\cdots \rightarrow \mathrm{H}^0(G_K, W(\delta\prod_{\sigma:K\hookrightarrow E}\sigma(x)^{
-k_{\sigma}}))\rightarrow \oplus_{\sigma:K\hookrightarrow E}\mathrm{H}^0(G_K, 
t^{-k_{\sigma}}W^+_{\mathrm{dR}}(\delta)_{\sigma}/W^+_{\mathrm{dR}}(\delta)_{\sigma})
\end{equation*}
\begin{equation*}
\rightarrow \mathrm{H}^1(G_K, W(\delta))\rightarrow \mathrm{H}^1(G_K, 
W(\delta\prod_{\sigma:K\hookrightarrow E}\sigma(x)^{
-k_{\sigma}}))\rightarrow\cdots.
\end{equation*}

The next lemma is a generalization of $\cite[\mathrm{Proposition}3.18]{Co07a}$.

\begin{lemma}$\label{30}$
Let $\{w(\delta)_{\sigma}\}_{\sigma}$ be 
the generalized Hodge-Tate weight of $W(\delta)$ defined in $\ref{-5}$.
Then there is an isomorphism of $E$-vector spaces 
\begin{itemize}
\item[]$\oplus_{\sigma:K\hookrightarrow E}\mathrm{H}^0(G_K, 
t^{-k_{\sigma}}W^+_{\mathrm{dR}}(\delta)_{\sigma}/W^+_{\mathrm{dR}}
(\delta)_{\sigma})\isom\oplus_{\sigma, w(\delta)_{\sigma}\in \{1,2,\cdots,k_{\sigma}\}}Ee_{\sigma}$.
\end{itemize}
Here $Ee_{\sigma}$ is a one dimensional $E$-vector space with base $e_{\sigma}$.
\end{lemma}
\begin{proof}
It suffices to show that $\mathrm{H}^0(G_K, 
t^{-k_{\sigma}}W^+_{\mathrm{dR}}(\delta)_{\sigma}/W^+_{\mathrm{dR}}
(\delta)_{\sigma})\isom E$ if and only if $w(\delta)_{\sigma}\in \{1,2,\cdots,k_{\sigma}\}$ and $\mathrm{H}^0(G_K, 
t^{-k_{\sigma}}W^+_{\mathrm{dR}}(\delta)_{\sigma}/W^+_{\mathrm{dR}}
(\delta)_{\sigma})=0$ otherwise. But, by definition of generalized Hodge-Tate 
weight, for any $i\in\mathbb{Z}$ and $\sigma$, we have $\mathrm{H}^0(G_K,t^{-i}W^+_{\mathrm{dR}}(\delta)_{\sigma}/t^{-i+1}W^+_{\mathrm{dR}}(\delta)_{\sigma})=E$ if and only if $w(\delta)_{\sigma}=i$ and we have $\mathrm{H}^0(G_K,t^{-i}W^+_{\mathrm{dR}}(\delta)_{\sigma}/t^{-i+1}W^+_{\mathrm{dR}}(\delta)_{\sigma})=0$ otherwise. From this, the result follows.
\end{proof}

By definition, for any continuous character $\delta:K^{\times}\rightarrow E^{\times}$,
we have the following short exact sequence
\begin{equation*}
0\rightarrow D^K_{\mathrm{dR}}(W(\delta))/(D^K_{\mathrm{cris}}(W(\delta))^{\varphi=1}+
\mathrm{Fil}^0D^K_{\mathrm{dR}}(W(\delta))) \rightarrow \mathrm{H}^1(G_K, W(\delta))
\end{equation*}
\begin{equation*}
\rightarrow \mathrm{Ker}(\mathrm{H}^1(G_K, W_e(\delta))\oplus \mathrm{H}^1(G_K,W^+_{\mathrm{dR}}(\delta))\rightarrow \mathrm{H}^1(G_K, W_{\mathrm{dR}}(\delta)))\rightarrow 0.
\end{equation*}
\begin{lemma}$\label{31}$
Let $\delta:K\rightarrow E^{\times}$ be a continuous character.
Then the $E$-vector space $\allowbreak D^K_{\mathrm{dR}}(W(\delta))/(D^K_{\mathrm{cris}}(W(\delta))^{\varphi=1}+\mathrm{Fil}^0D^K_{\mathrm{dR}}(W(\delta)))$ is isomorphic to 
\begin{itemize}
\item[$\mathrm{(1)}$]$\oplus_{\sigma, w(\delta)_{\sigma}\in\mathbb{Z}_{\ge 1}}Ee_{\sigma}/\Delta(E)$ when $\delta=\prod_{\sigma:K\hookrightarrow E}\sigma(x)^{k_{\sigma}}$ for some $\{k_{\sigma}\}_{\sigma}$ such that $k_{\sigma}\in \mathbb{Z}$ for any $\sigma$, 
\item[$\mathrm{(2)}$] $\oplus_{\sigma, w(\delta)_{\sigma}\in\mathbb{Z}_{\ge 1}}Ee_{\sigma}$ otherwise.
\end{itemize}
Here $\Delta:E\rightarrow \oplus_{\sigma. w(\delta)_{\sigma}\in\mathbb{Z}_{\ge 1}}Ee_{\sigma}$ is the diagonal map.
\end{lemma}
\begin{proof}
By definition of $W(\delta)$, it is easy to see that $W_{\mathrm{dR}}(\delta)\isom B_{\mathrm{dR}}
\otimes_{\mathbb{Q}_p}E(\delta_0)$ for some continuous character $\delta_0:G_K\rightarrow 
\mathcal{O}_E^{\times}$. If we decompose $B_{\mathrm{dR}}\otimes_{\mathbb{Q}_p}E(\delta_0)\isom
\oplus_{\sigma:K\hookrightarrow E}B_{\mathrm{dR}}\otimes_{K,\sigma}E(\delta_0)$, we have $(B_{\mathrm{dR}}\otimes_{K,\sigma}E(\delta_0))^{G_K}=0$ if and only if $w(\delta)_{\sigma}\notin\mathbb{Z}$. If $w(\delta)_{\sigma}\in\mathbb{Z}$, then we have 
$(t^{-w(\delta)_{\sigma}}B_{\mathrm{dR}}^+\otimes_{K,\sigma}E(\delta_0))^{G_K}\isom E$ 
and $(t^{-w(\delta)_{\sigma}+1}B_{\mathrm{dR}}^+\otimes_{K,\sigma}E(\delta_0))^{G_K}=0$. 
As for $W_e(\delta)$, we have $D^K_{\mathrm{cris}}(W(\delta))^{\varphi=1} =W_e(\delta)^{G_K}\not= 0$
if and only if $W_e(\delta)\isom B_e\otimes_{\mathbb{Q}_p}E$. Then $W(\delta)\isom (B_e\otimes_{\mathbb{Q}_p}E, \oplus_{\sigma:K\hookrightarrow E}t^{k_{\sigma}}B^+_{\mathrm{dR}}\otimes_{K,\sigma}E)$ for some $\{k_{\sigma}\}_{\sigma}$ such that $k_{\sigma}\in \mathbb{Z}$ for any $\sigma$.
Then, by Lemma $\ref{26}$, we have $\delta=\prod_{\sigma:K\hookrightarrow E}\sigma(x)^{k_{\sigma}}$. 
In this case, we have $D^K_{\mathrm{cris}}(W(\delta))^{\varphi=1} \isom E$.
Combining these computations, we get the lemma.
\end{proof}
Next we compare the following two maps, one is the boundary map 
\begin{equation*}\partial:\oplus_{\sigma:K\hookrightarrow E}\mathrm{H}^0(G_K,t^{-k_{\sigma}}W^+_{\mathrm{dR}}(\delta)_{\sigma}/W^+_{\mathrm{dR}}(\delta)_{\sigma})\rightarrow \mathrm{H}^1(G_K, W(\delta))
\end{equation*}
 of the exact seuence before Lemma $\ref{30}$ and the other 
is the natural map 
\begin{align*}
\iota:D^K_{\mathrm{dR}}(W(\delta))/\mathrm{Fil}^0D^K_{\mathrm{dR}}(W(\delta))&\rightarrow D^K_{\mathrm{dR}}(W(\delta))/(D^K_{\mathrm{cris}}(W(\delta))^{\varphi=1}+
\mathrm{Fil}^0D^K_{\mathrm{dR}}(W(\delta))) \\
 &\rightarrow \mathrm{H}^1(G_K, W(\delta)).
\end{align*}
From the proof of Lemma $\ref{31}$, we can see that the natural map 
$f:D^K_{\mathrm{dR}}(W(\delta))/\mathrm{Fil}^0D^K_{\mathrm{dR}}(W(\delta))
\allowbreak\isom (W_{\mathrm{dR}}(\delta)/W^+_{\mathrm{dR}}(\delta))^{G_K}$ is an isomorphism. On the other hand, we have a natural map 
\begin{align*}
 j:\oplus_{\sigma:K\hookrightarrow E}\mathrm{H}^0(G_K,t^{-k_{\sigma}}W^+_{\mathrm{dR}}(\delta)_{\sigma}/W^+_{\mathrm{dR}}(\delta)_{\sigma})&=\mathrm{H}^0(G_K, W^+_{\mathrm{dR}}(\delta\prod_{\sigma:K\hookrightarrow E}\sigma(x)^{-k_{\sigma}})/W^+_{\mathrm{dR}}(\delta)) \\  
 &\hookrightarrow (W_{\mathrm{dR}}(\delta)/W^+_{\mathrm{dR}}(\delta))^{G_K}.
 \end{align*}
\begin{lemma}$\label{32}$
With the above notations, we have 
$\iota\circ f^{-1}\circ j=(-1)\partial$.
\end{lemma}
\begin{proof}
This follows easily from diagram chase.
\end{proof}

\section{Classification of two dimensional split 
trianguline $E$-representations.}
In this section, we classify two dimensional split 
trianguline $E$-representations. As in $\cite[0.2]{Co07a}$, we 
would like to explicitly determine the parameter spaces of two dimensional 
split trianguline $E$-representations.

\subsection{Parameter spaces of two dimensional split trianguline 
$E$-represent$\-\allowbreak\,$ations.}
Let $W$ be a rank two split trianguline $E$-$B$-pair such that 
$[W]\in \mathrm{Ext}^1(W(\delta_2), W(\delta_1))$ for 
some continuous characters $\delta_1, \delta_2:K^{\times}\rightarrow E^{\times}$.
First we study a necessary condition on $(\delta_1,\delta_2)$ in order   
for $W$ to be $\mathrm{\acute{e}}$tale, i.e. of the form $W=W(V)$ for some $E$-representation $V$.
\begin{lemma}$\label{33}$
If $[W]\in \mathrm{Ext}^1(W(\delta_2), W(\delta_1))$ is $\mathrm{\acute{e}}$tale, 
then the pair $(\delta_1,\delta_2)$ satisfies $\mathrm{val}_p(\delta_1(\pi_K))+\mathrm{val}_p(\delta_2(\pi_K))=0$ and $\mathrm{val}_p(\delta_1(\pi_K))\ge 0$. 
Here $\mathrm{val}_p$ is the valuation of $E$ such that $\mathrm{val}_p(p)=1$.
\end{lemma}
\begin{proof}
Let $W$ be of the form $W=W(V)$ for some $E$-representation $V$.
By definition, $W(V)$ sits in the following short exact sequence of 
$E$-$B$-pairs
\begin{equation*}
0\rightarrow W(\delta_1)\rightarrow W(V) \rightarrow W(\delta_2)\rightarrow 0.
\end{equation*}
Because $W(V)$ is pure of slope 0, $\mathrm{det}W(V)\isom W(\delta_1\delta_2)$ is also pure of slope 0, so we have a condition $\frac{\mathrm{val}_p(\delta_1(\pi_K))+\mathrm{val}_p(\delta_2(\pi_K))}{f}=0$ because 
the slope of $W(\delta)$ is $\frac{\mathrm{val}_p(\delta(\pi_K))}{f}$. 
By the slope filtration theorem for $E$-$B$-pairs (Theorem $\ref{11}$), $W(\delta_1)$ must be pure of slope $\ge0$. This implies that $\frac{\mathrm{val}_p(\delta_1(\pi_K))}{f}\ge 0$. We have proved 
the lemma.
\end{proof}

To describe the classification of two dimensional split trianguline 
$E$-representations, let us introduce several notations.
First let us put
\begin{itemize}
\item[$\bullet$]$S^+:=\{(\delta_1,\delta_2)| \delta_1,\delta_2:K^{\times}\rightarrow E^{\times}$ continuous characters such that $\mathrm{val}_p(\delta_1(\pi_K))+\mathrm{val}_p(\delta_2(\pi_K))=0, \mathrm{val}_p(\delta_1(\pi_K))\ge 0\}$.
\end{itemize}
By Lemma $\ref{33}$, for classifying two dimensional split trianguline 
$E$-representations, it suffices only to consider split trianguline $E$-$B$-pairs $W$ such that $[W]\in \mathrm{Ext}^1(W(\delta_2), W(\delta_1))$ for 
some $(\delta_1,\delta_2)\in S^+$.

First we classify two dimensional split trianguline $E$-representations $V$ 
such that $[W(V)]=0\in \mathrm{Ext}^1(W(\delta_2), W(\delta_1))$ for 
some $(\delta_1,\delta_2)$, i.e. $W(V)=W(\delta_1)\oplus W(\delta_2)$.
\begin{lemma}$\label{33.5}$
Let $(\delta_1,\delta_2)\in S^+$.
Then $W(\delta_1)\oplus W(\delta_2)$ is of the form $W(V)$ for some two dimensional $E$-representation if and only if both $W(\delta_1)$ and $W(\delta_2)$ are pure of slope zero, i.e. $W(\delta_1)\isom W(V(E(\tilde{\delta}_1)))$ and 
$W(\delta_2)\isom W(V(E(\tilde{\delta}_2)))$ for some continuous characters 
$\tilde{\delta}_1, \tilde{\delta}_2:G_K\rightarrow \mathcal{O}^{\times}_E $.
 In this case, $V\isom E(\tilde{\delta}_1)\oplus E(\tilde{\delta}_2)$.
\end{lemma}
\begin{proof}
``If" part is trivial.
Let us assume that $W(\delta_1)\oplus W(\delta_2)\isom W(V)$ for some two dimensional $E$-representation $V$. So $W(\delta_1)\oplus W(\delta_2)$ is pure of slope zero. If $W(\delta_1)$ is not pure of slope zero, then this is pure of slope $u>0$ by the assumption that $(\delta_1,\delta_2)\in S^+$. Then $W(\delta_2)$ is pure of slope $-u<0$ because $(\delta_1,\delta_2)\in S^+$. Then we have 
the exact sequence
\begin{equation*}
0\rightarrow W(\delta_2)\rightarrow W(\delta_1)\oplus W(\delta_2)\rightarrow W(\delta_1)\rightarrow 0.
\end{equation*}
Then the slope filtration theorem for $E$-$B$-pairs (Theorem $\ref{11}$) 
implies that $W(\delta_1)\oplus W(\delta_2)$ is not pure of slope zero.
This is a contradiction.
\end{proof}
This lemma says that two dimensional split trianguline $E$-representation $V$ 
such that $[W(V)]=0\in \mathrm{Ext}^1(W(\delta_2), W(\delta_1))$ corresponds 
to two dimensional $E$-representations $V$ such that $V=E(\tilde{\delta}_1)\oplus E(\tilde{\delta}_2)$ for some continuous characters 
$\tilde{\delta}_1, \tilde{\delta}_2:G_K\rightarrow \mathcal{O}^{\times}_E $.
So, in this case, we finish the classification.

From now on, we only consider split trianguline $E$-$B$-pairs $W$ such that 
$[W]\not= 0\in \mathrm{Ext}^1(W(\delta_2), W(\delta_1))$.
For any $(\delta_1,\delta_2)\in S^+$, let us put
\begin{itemize}
\item[$\bullet$]$S(\delta_1,\delta_2):=\mathbb{P}_E(\mathrm{H}^1(G_K, W(\delta_1/\delta_2)))$.
\end{itemize}
Here, for any finite dimensional $E$-vector space $M$, we denote 
\begin{itemize}
\item[]$\mathbb{P}_E(M):=\{[v]|v\in M-\{0\}, [v]=[v'] \iff v'=av$ for some $a\in E^{\times}\}$.
\end{itemize}
Next, for any $(\delta_1,\delta_2)\in S^+$,  we define a subset 
$S'(\delta_1,\delta_2)\subseteq S(\delta_1,\delta_2)$ by 
\begin{itemize}
\item[$\bullet$]$S'(\delta_1,\delta_2):=\mathbb{P}_E(D^K_{\mathrm{dR}}(W(\delta_1/\delta_2))/
(D^K_{\mathrm{cris}}(W(\delta_1/\delta_2))^{\varphi=1}+\mathrm{Fil}^0D^K_{\mathrm{dR}}(W(\delta_1/\delta_2))))$,
\end{itemize}
here we see $S'(\delta_1,\delta_2)$ as a subset of $S(\delta_1,\delta_2)$ by the following 
natural inclusion
\begin{itemize}
\item[]$D^K_{\mathrm{dR}}(W(\delta_1/\delta_2))/
(D^K_{\mathrm{cris}}(W(\delta_1/\delta_2))^{\varphi=1}+\mathrm{Fil}^0D^K_{\mathrm{dR}}(W(\delta_1/\delta_2)))\hookrightarrow \mathrm{H}^1(G_K, W(\delta_1/\delta_2))$.
\end{itemize}
Then, by Lemma $\ref{31}$, $S'(\delta_1,\delta_2)$ is non canonically isomorphic to 
\begin{itemize}
\item[$\mathrm{(1)}$]$S'(\delta_1,\delta_2)\isom\mathbb{P}_E(\oplus_{\sigma, w(\delta_1/\delta_2)_{\sigma}\ge 1}Ee_{\sigma}/\Delta(E))$ when $\delta_1/\delta_2=\prod_{\sigma:K\hookrightarrow E}\sigma(x)^{k_{\sigma}}$ for some $\{k_{\sigma}\}_{\sigma}$ such that $k_{\sigma}\in \mathbb{Z}$ for any $\sigma$, 
\item[$\mathrm{(2)}$]$S'(\delta_1,\delta_2)\isom\mathbb{P}_E(\oplus_{\sigma, w(\delta_1/\delta_2)_{\sigma}\in\mathbb{Z}_{\ge 1}}Ee_{\sigma})$ otherwise,
\end{itemize}
here these isomorphisms depend on the choice of $E$-linear isomorphisms in Lemma $\ref{31}$.

Finally we define the subset ${S'}^{\acute{\mathrm{e}}t}(\delta_1,\delta_2)$ of $S'(\delta_1,\delta_2)$ as follows.
When $\delta_1/\delta_2\not= \prod
_{\sigma:K\hookrightarrow E}\sigma(x)^{k_{\sigma}}$ for any $\{k_{\sigma}\}$ such that $k_{\sigma}\in\mathbb{Z}$ for any $\sigma$, we fix an isomorphism 
$S'(\delta_1,\delta_2)\isom\mathbb{P}_E(\oplus_{\sigma, w(\delta_1/\delta_2)_{\sigma}\in\mathbb{Z}_{\ge 1}}Ee_{\sigma})$ as above and identify these two spaces. Then let us put
\begin{itemize}
\item[$\bullet$]${S'}^{\acute{\mathrm{e}}t}(\delta_1,\delta_2):=\{[(a_{\sigma}e_{\sigma})_{\sigma}]\in \mathbb{P}_E(\oplus_{\sigma, w(\delta_1/\delta_2)_{\sigma}\in\mathbb{Z}_{\ge 1}}Ee_{\sigma})| (\sum_{\sigma, a_{\sigma}\not= 0}w(\delta_1/\delta_2)_{\sigma})+ e_K \mathrm{val}_p
(\delta_2(\pi_K))\allowbreak\ge 0\}$.
\end{itemize} 
When $\delta_1/\delta_2= \prod
_{\sigma:K\hookrightarrow E}\sigma(x)^{k_{\sigma}}$ for some $\{k_{\sigma}\}_{\sigma}$ such that $k_{\sigma}\in\mathbb{Z}$ for any $\sigma$, we fix an isomorphism $S'(\delta_1,\delta_2)\isom \mathbb{P}_E(\oplus_{\sigma, w(\delta_1/\delta_2)_{\sigma}\ge 1}Ee_{\sigma}/\Delta(E))$ as above and identify these two spaces. Then let us put
\begin{itemize}
\item[$\bullet$]${S'}^{\acute{\mathrm{e}}t}(\delta_1,\delta_2):=
\{[\overline{(a_{\sigma}e_{\sigma})_{\sigma}}]\in \mathbb{P}_E(\oplus_{\sigma. w(\delta_1/\delta_2)_{\sigma}\ge 1}Ee_{\sigma}/\Delta(E))| (\sum_{\sigma, a_{\sigma}\not= 0}w(\delta_1/\delta_2)_{\sigma})+ e_K \mathrm{val}_p(\allowbreak\delta_2(\pi_K))\ge 0$ for any lifting $(a_{\sigma}e_{\sigma})_{\sigma}\in
\oplus_{\sigma, w(\delta_1/\delta_2)_{\sigma}\in\mathbb{Z}_{\ge 1}}Ee_{\sigma}$ of $[\overline{(a_{\sigma}e_{\sigma})_{\sigma}}]\}$. 
\end{itemize}
We put ${S'}^{non-\acute{\mathrm{e}}t}(\delta_1,\delta_2):=S'(\delta_1,\delta_2)\setminus
 {S'}^{\acute{\mathrm{e}}t}(\delta_1,\delta_2)$. Then we can easily see that the subsets ${S'}^{\acute{\mathrm{e}}t}(\delta_1,\delta_2)$, 
 ${S'}^{non-\acute{\mathrm{e}}t}(\delta_1,\delta_2)\subseteq S(\delta_1,\delta_2)$ do not depend on the choice 
 of an isomorphism $S'(\delta_1,\delta_2)\isom\mathbb{P}_E(\oplus_{\sigma, w(\delta_1/\delta_2)_{\sigma}\ge 1}Ee_{\sigma}/\Delta(E))$ or $S'(\delta_1,\delta_2)\isom\mathbb{P}_E(\oplus_{\sigma, w(\delta_1/\delta_2)_{\sigma}\in\mathbb{Z}_{\ge 1}}Ee_{\sigma})$.
 \begin{rem}
 When $K=\mathbb{Q}_p$, $S(\delta_1,\delta_2)$ is one point or $\mathbb{P}_E(E)$ and $S'(\delta_1,\delta_2)$ is empty or one point or $\mathbb{P}_E(E)$ and ${S'}^{\acute{\mathrm{e}}t}(\delta_1,\delta_2)$ is empty 
 or one point or $\mathbb{P}_E(E)$.
 If we compare this parameter space and Colmez's one ($\cite[4.3]{Co07a}$), 
 then we have $\sqcup_{(\delta_1,\delta_2)\in S^+}(S(\delta_1,\delta_2)\setminus S'(\delta_1,\delta_2))=\mathcal{S}_+^{\mathrm{ng}}\sqcup S_+^{\mathrm{st}}$, $\sqcup_{(\delta_1,\delta_2)\in
  S^+}{S'}^{\acute{\mathrm{e}}t}(\delta_1,\delta_2)=\mathcal{S}_+^{\mathrm{cris}}\sqcup \mathcal{S}_+^{\mathrm{ord}}$ and $\sqcup_{(\delta_1,\delta_2)\in S^+}{S'}^{non-\acute{\mathrm{e}}t}(\delta_1,\delta_2)=\mathcal{S}^{\mathrm{ncl}}$.
  \end{rem} 
For any $s\in S(\delta_1,\delta_2)$, we denote $W(s)$ an extension of $W(\delta_1)$ by $W(\delta_2)$ defined by $s$. The isomorphism class of $W(s)$ as $E$-$B$-pair depends only on the class $s$. By definition, $W(s)$ is a split trianguline $E$-$B$-pair which sits in a following non-split short exact sequence of $E$-$B$-pairs
\begin{equation*}
0\rightarrow W(\delta_1)\rightarrow W(s)\rightarrow W(\delta_2)\rightarrow 0.
\end{equation*}

We determine when $W(s)$ is of the form $W(V(s))$ for some $E$-representation $V(s)$.
The following theorem is a generalization of $\cite[\mathrm{Proposition}\,4.7]{Co07a}$ and is the most important step of the classification.
 \begin{thm}$\label{34}$
 Let $s\in S(\delta_1,\delta_2)$ for $(\delta_1,\delta_2)\in S^+$.
 Then the following conditions are equivalent.
 \begin{itemize}
 \item[$\mathrm{(1)}$]$W(s)$ is pure of slope zero, i.e. 
 $W(s)\isom W(V(s))$ for a two dimensional split trianguline 
$E$-representation $V(s)$.
 \item[$\mathrm{(2)}$] $s\notin {S'}^{non-\acute{\mathrm{e}}t}(\delta_1,\delta_2)$.
 \end{itemize}
 \end{thm}
 \begin{proof}
 First we prove the theorem in the case $\delta_1/\delta_2\not=\prod
_{\sigma:K\hookrightarrow E}\sigma(x)^{k_{\sigma}}$ for any $\{k_{\sigma}\}_{\sigma}$ such that $k_{\sigma}\in\mathbb{Z}$ for any $\sigma$.
Let us assume that $W(s)$ is not pure of slope zero. Then, by the slope filtration theorem for 
$E$-$B$-pairs (Theorem $\ref{11}$), there 
exist rank one $E$-$B$-pairs $W(\delta_3)$, $W(\delta_4)$ such that $W(\delta_3)$ is 
pure of slope $u<0$, $W(\delta_4)$ is pure of slope $-u>0$ and $W(s)$ sits in the 
following short exact sequence
\begin{equation*}
0\rightarrow W(\delta_3)\rightarrow W(s)\rightarrow W(\delta_4)\rightarrow 0.
\end{equation*}
On the other hand, by definition we have a following short exact sequence
\begin{equation}
0 \rightarrow W(\delta_1)\rightarrow W(s)\rightarrow W(\delta_2)\rightarrow 0.
\end{equation}
We consider the restriction to $W(\delta_3)$ of the projection $W(s)\rightarrow W(\delta_2)$,
which we denote $f:W(\delta_3)\rightarrow W(\delta_2)$. We claim that 
$f$ is an inclusion. Because $\mathrm{Ker}(f)$ and $\mathrm{Im}(f)$ are $E$-$B$-pairs by Lemma $\ref{2}$ and because $W(\delta_3)$ is rank one,
we have exactly one of the following two cases, $\mathrm{Ker}(f)\isom W(\delta_3)$ or $W(\delta_3)\isom \mathrm{Im}(f)$. If $\mathrm{Ker}(f)\isom W(\delta_3)$ then the inclusion $W(\delta_3)\hookrightarrow W(s)$ factors through a non zero map 
$f':W(\delta_3)\rightarrow W(\delta_1)$. But because the slope of $W(\delta_1)$ is strictly bigger 
than the slope of $W(\delta_3)$ 
by assumption, so we have $\mathrm{Hom}(W(\delta_3), W(\delta_1))=0$ 
by Proposition $\ref{28}$ or by $\cite[\mathrm{Lemma}\,6.2]{Ke07}$.
This is a contradiction. So we have $W(\delta_3)\isom \mathrm{Im}(f)$, i.e. $f$ is an inclusion. 
Then, by Proposition $\ref{28}$, we have $\delta_3=
\delta_2\prod_{\sigma:K\hookrightarrow E}\sigma(x)^{k_{\sigma}}$ for some $\{k_{\sigma}\}_{\sigma}$ 
such that $k_{\sigma}\in\mathbb{Z}_{\ge 0}$ for any $\sigma$. Then, because $\mathrm{det}(W(s))
=W(\delta_1\delta_2)\isom W(\delta_3\delta_4)$, we have $\delta_4=\delta_1\prod_{\sigma:K\hookrightarrow E}\sigma(x)^{-k_{\sigma}}$. 
The injectivity of $f$ means that the exact sequence $(1)$ splits after pulling back by $f:W(\delta_3)\hookrightarrow 
W(\delta_2)$, i.e. we have 
\begin{itemize}
\item[]$s\in \mathrm{Ker}(\mathrm{H}^1(G_K, W(\delta_1/\delta_2))\rightarrow 
\mathrm{H}^1(G_K, W(\prod_{\sigma:K\hookrightarrow E}\sigma(x)^{-k_{\sigma}}\delta_1/\delta_2)))$. 
\end{itemize}
This kernel is isomorphic to 
\begin{itemize}
\item[]$\mathrm{Im}(\partial:\mathrm{H}^0(G_K, \oplus_{\sigma:K\hookrightarrow E}t^{-k_{\sigma}}W^+_{\mathrm{dR}}(\delta_1/\delta_2)_{\sigma}/W^+_{\mathrm{dR}}(\delta_1/\delta_2)_{\sigma})\hookrightarrow 
\mathrm{H}^1(G_K, W(\delta_1/\delta_2)))$.
\end{itemize}
So, by Lemma $\ref{32}$, we have $s=[(a_{\sigma}e_{\sigma})_{\sigma}]\in \mathbb{P}_E(\oplus_{\sigma, w(\delta_1/\delta_2)_{\sigma}\in\mathbb{Z}_{\ge1}}Ee_{\sigma})=S'(\delta_1,\delta_2)$ such that any $\sigma$ with $a_{\sigma}\not=0$ satisfies $w(\delta_1/\delta_2)_{\sigma}\in \{1,2,\cdots, k_{\sigma}\}$. Then, $(\sum_{\sigma,a_{\sigma}\not= 0}w(\delta_1/\delta_2)_{\sigma})+e_K \mathrm{val}_p(\delta_2(\pi_K))\le (\sum_{\sigma}k_{\sigma})+e_K\mathrm{val}_p(\delta_2(\pi_K))=[K:\mathbb{Q}_p]$(slope of $W(\delta_3))<0$ by assumption. So $s\in {S'}^{non-\acute{\mathrm{e}}t}(\delta_1,\delta_2)$.

Next we assume that $s=[(a_{\sigma}e_{\sigma})_{\sigma}]\in {S'}^{non-\acute{\mathrm{e}}t}(\delta_1,\delta_2)$.
Then, by Lemma $\ref{32}$,  we can see that $s$ is contained in the image of  
\begin{itemize}
\item[]$\partial:\oplus_{\sigma,a_{\sigma}\not= 0}\mathrm{H}^0(G_K,t^{-w(\delta_1/\delta_2)_{\sigma}}W^+_{\mathrm{dR}}(\delta_1/\delta_2)_{\sigma}/W^+_{\mathrm{dR}}(\delta_1/\delta_2)_{\sigma})\rightarrow \mathrm{H}^1(G_K, W(\delta_1/\delta_2))$. 
\end{itemize}
So we have 
\begin{itemize}
\item[]$s\in
\mathrm{Ker}(\mathrm{H}^1(G_K,W(\delta_1/\delta_2))\rightarrow \mathrm{H}^1
(G_K,W(\prod_{\sigma,a_{\sigma}\not= 0}\sigma(x)^{-w(\delta_1/\delta_2)_{\sigma}}\delta_1/\delta_2))$. 
\end{itemize}
So there is an injection $g:W(\prod_{\sigma,a_{\sigma}\not= 0}\sigma(x)^{w(\delta_1/\delta_2)_{\sigma}}\delta_2)\hookrightarrow W(s)$ such that $f 'g$ is the natural inclusion $W(\prod_{\sigma,a_{\sigma}\not= 0}\sigma(x)^{w(\delta_1/\delta_2)_{\sigma}}\delta_2)\hookrightarrow W(\delta_2)$ (here $f':W(s)\rightarrow W(\delta_2)$ is the projection). If we take the saturation $W(\delta_3)$ of this 
inclusion $g$ (Lemma $\ref{3}$) and write the cokernel by $W(\delta_4)$, we have the following short exact sequence
\begin{equation*}
0\rightarrow W(\delta_3)\rightarrow W(s)\rightarrow W(\delta_4)\rightarrow 0.
\end{equation*}
Then, by Proposition $\ref{28}$, we have (the slope of $W(\delta_3))\le$ (the slope of $W(\prod_{\sigma,a_{\sigma}\not= 0}\sigma(x)^{w(\delta_1/\delta_2)_{\sigma}}\allowbreak\delta_2))=\frac{1}{[K:\mathbb{Q}_p]}(\sum_{\sigma,a_{\sigma}\not= 0}w(\delta_1/\delta_2)_{\sigma})
+\frac{\mathrm{val}_p(\delta_2(\pi_K))}{f}<0$ because $s\in S^{' non-et}(\delta_1,\delta_2)$. So $W(s)$ is not pure of slope $0$ by the slope filtration theorem.
So we have finished the proof when $\delta_1/\delta_2\not=\prod
_{\sigma:K\hookrightarrow E}\sigma(x)^{k_{\sigma}}$ for any $\{k_{\sigma}\}$ such that $k_{\sigma}\in\mathbb{Z}$ for any $\sigma$.

Next we prove the theorem in the case $\delta_1/\delta_2=\prod
_{\sigma:K\hookrightarrow E}\sigma(x)^{k_{\sigma}}$ for some $\{k_{\sigma}\}_{\sigma}$ such that $k_{\sigma}\in\mathbb{Z}$ for any $\sigma$.
First let us assume that $W(s)$ is not pure of slope 0.
Then, as in the first case, there are rank one $E$-$B$-pairs $W(\delta_3), W(\delta_4)$ such that 
the slope of $W(\delta_3):=u<0$ and $W(s)$ sits in the following short exact sequence
\begin{equation*}
0\rightarrow W(\delta_3)\rightarrow W(s)\rightarrow W(\delta_4) \rightarrow 0.
\end{equation*}
Then we can prove as in the first case that the induced map $W(\delta_3)\hookrightarrow W(\delta_2)$ is injective and $\delta_3=\delta_2\prod_{\sigma:K\hookrightarrow E}\sigma(x)^{k_{\sigma}}$ for some $\{k_{\sigma}\}_{\sigma}$ such that 
$k_{\sigma}\in\mathbb{Z}_{\ge 0}$ for any $\sigma$. So we have  
\begin{align*}
s\in& \,\mathrm{Ker}(\mathrm{H}^1(G_K,W(\delta_1/\delta_2))\rightarrow \mathrm{H}^1(G_K, W(\prod_{\sigma:K\hookrightarrow E}\sigma(x)^{-k_{\sigma}} \delta_1/\delta_2)))  \\
 =&\,\mathrm{Im}(\partial:\mathrm{H}^0(G_K, \oplus_{\sigma:K\hookrightarrow E}t^{-k_{\sigma}}W^+_{\mathrm{dR}}(\delta_1/\delta_2)_{\sigma}/W^+_{\mathrm{dR}}(\delta_1/\delta_2)_{\sigma})\rightarrow \mathrm{H}^1(G_K, W(\delta_1/\delta_2))). 
\end{align*}
By Lemma $\ref{32}$, we can see that  $s$ corresponds to 
$s=[\overline{(a_{\sigma}e_{\sigma})_{\sigma}}]\in \mathbb{P}_E(\oplus_{\sigma, w(\delta_1/\delta_2)_{\sigma}\in\mathbb{Z}_{\ge 1}}Ee_{\sigma}/\Delta(E))\allowbreak\isom S'(\delta_1,\delta_2)$ such that there is a lifting $[(a_{\sigma}e_{\sigma})_{\sigma}]$ of $[\overline{(a_{\sigma}e_{\sigma})_{\sigma}}]$ such that $a_{\sigma}=0$ for any $\sigma$ satisfying $w(\delta_1/\delta_2)_{\sigma}\notin\{1,2,\cdots, k_{\sigma}\}$. So $(\sum_{\sigma, a_{\sigma}\not= 0}w(\delta_1/\delta_2)_{\sigma})+ e_K \mathrm{val}_p(\delta_2(\pi_K))\le (\sum_{\sigma:K\hookrightarrow E}k_{\sigma}) + e_K \mathrm{val}_p(\delta_2(\pi_K))=[K:\mathbb{Q}_p]u< 0$.
So $s$ is contained in ${S'}^{non-\acute{\mathrm{e}}t}(\delta_1,\delta_2)$.

Next let us assume that $s\in {S'}^{non-\acute{\mathrm{e}}t}(\delta_1,\delta_2)$.
Then, by definition, $s=[\overline{(a_{\sigma}e_{\sigma})_{\sigma}}]\in \mathbb{P}_E(\oplus_{\sigma, w(\delta_1/\delta_2)_{\sigma}\in \mathbb{Z}_{\ge 1}}Ee_{\sigma}/\Delta(E))\isom S'(\delta_1,\delta_2)$ such that some lift $[(a_{\sigma}e_{\sigma})_{\sigma}]$ of $ [\overline{(a_{\sigma}e_{\sigma})_{\sigma}}]$ satisfies $(\sum_{\sigma,a_{\sigma}\not= 0}w(\delta_1/\delta_2)_{\sigma})+ e_K \mathrm{val}_p(\delta_2(\pi_K)) < 0$.
Then, by Lemma $\ref{32}$, we have 
\begin{itemize}
\item[]$s\in \mathrm{Ker}(\mathrm{H}^1(G_K, W(\delta_1/\delta_2))\rightarrow 
\mathrm{H}^1(G_K, W(\prod_{\sigma,a_{\sigma}\not= 0}\sigma(x)^{-w(\delta_1/\delta_2)_{\sigma}}\delta_1/\delta_2)))$. 
\end{itemize}
In particular, as in the proof of the first case, we can see that the natural inclusion 
$W(\prod_{\sigma,a_{\sigma}\not= 0}\allowbreak\sigma(x)^{w(\delta_1/\delta_2)_{\sigma}}\delta_2)\hookrightarrow W(\delta_2)$ factors through an inclusion $g:W(\prod_{\sigma,a_{\sigma}\not= 0}\sigma(x)^{w(\delta_1/\delta_2)_{\sigma}}\delta_2)\hookrightarrow W(s)$. If we take the saturation $W(\delta_3)$ of $g$ and write the cokenel by $W(\delta_4)$, then we have (the slope of $W(\delta_3))\le$ (the slope of $W(\prod_{\sigma,a_{\sigma}\not= 0}\sigma(x)^{w(\delta_1/\delta_2)_{\sigma}}\delta_2))=
\frac{1}{[K:\mathbb{Q}_p]}((\sum_{\sigma,a_{\sigma}\not= 0}w(\delta_1/\delta_2)_{\sigma})+ e_K 
\mathrm{val}_p(\delta_2(\pi_K)))<0$, where the first inequality follows from Proposition $\ref{28}$. So $W(s)$ is not pure of slope $0$ by 
the slope filtration theorem.
So we have finished the proof of this theorem in all the cases.
\end{proof}

By this theorem, we can determine all the two dimensional split 
trianguline $E$-representations. For any $s\in S(\delta_1,\delta_2)\setminus {S'}^{non-\acute{\mathrm{e}}t}(\delta_1,\delta_2)$, we write $V(s)$ the two dimensional split trianguline 
$E$-representation such that $W(V(s))= W(s)$.

\subsection{Irreducibility of $V(s)$.}

Next we determine when $V(s)$ is irreducible as $E$-representation of $G_K$.
For this, we put 
\begin{itemize}
\item[$\bullet$]$S^+_{0}:=\{(\delta_1,\delta_2)\in S^+|\mathrm{val}_p(\delta_1(\pi_K))
=\mathrm{val}_p(\delta_2(\pi_K))=0\}$,
\item[$\bullet$]$S^{+}_{*}:=S^+\setminus S^{+}_{0}$.
\end{itemize}
For any $(\delta_1,\delta_2)\in S^+_{*}$, we put 
\begin{itemize}
\item[(1)]${S'}^{\mathrm{ord}}(\delta_1,\delta_2):=\{[\overline{(a_{\sigma}e_{\sigma})_{\sigma}}]\in S^{' et}(\delta_1,\delta_2)|(\sum_{\sigma,a_{\sigma}\not= 0}
w(\delta_1/\delta_2)_{\sigma})+e_K\mathrm{val}_p(\delta_2(\pi_K))=0 $ for some lifting $[(a_{\sigma}e_{\sigma})_{\sigma}]$ of $[\overline{(a_{\sigma}e_{\sigma})_{\sigma}}] \}$ when $\delta_1/\delta_2= \prod_{\sigma:K\hookrightarrow E}\sigma(x)^{k_{\sigma}}$ for some $\{k_{\sigma}\}_{\sigma}$ such that $k_{\sigma}\in\mathbb{Z}$ for any $\sigma$,
\item[(2)]${S'}^{\mathrm{ord}}(\delta_1,\delta_2):=\{[a_{\sigma}e_{\sigma}]\in S^{' et}(\delta_1,\delta_2)|(\sum_{\sigma,a_{\sigma}\not= 0}w(\delta_1/\delta_2)_{\sigma})+e_K\mathrm{val}_p(\delta_2(\pi_K))=0\} $ otherwise. 
\end{itemize}
We can easily see that the subset ${S'}^{\mathrm{ord}}(\delta_1,\delta_2)$ does not depend on the choice of an isomorphism $S'(\delta_1,\delta_2)\isom\mathbb{P}_E(\oplus_{\sigma, w(\delta_1/\delta_2)_{\sigma}\ge 1}Ee_{\sigma}/\Delta(E))$ or $S'(\delta_1,\delta_2)\isom\mathbb{P}_E(\oplus_{\sigma, w(\delta_1/\delta_2)_{\sigma}\in\mathbb{Z}_{\ge 1}}Ee_{\sigma})$.
\begin{rem}
When $K=\mathbb{Q}_p$, $\mathcal{S}_+^{\mathrm{ord}}$ in $\cite[4.3]{Co07a}$ is equal to $\sqcup_{(\delta_1,\delta_2)\in S^+_{*}}S^{'\mathrm{ord}}(\delta_1,\delta_2)$ and $\mathcal{S}_0$ in $\cite[4.3]{Co07a}$ is equal to 
$\sqcup_{(\delta_1,\delta_2)\in S^+_0}S(\delta_1,\delta_2)$
\end{rem}
 The following proposition is a generalization of $\cite[\mathrm{Proposition}\,5.7,\, 5.8]{Co07a}$.
\begin{prop}$\label{35}$
Let $(\delta_1,\delta_2)\in S^+$ and $s\in S(\delta_1,\delta_2)\setminus {S'}^{non-\acute{\mathrm{e}}t}(\delta_1,\delta_2)$. Then 
the following conditions are equivalent.
\begin{itemize}
\item[$\mathrm{(1)}$]$V(s)$ is irreducible.
\item[$\mathrm{(2)}$]$(\delta_1,\delta_2)\notin S^+_{0}$ and $s\notin {S'}^{\mathrm{ord}}(\delta_1,\delta_2)$.
\end{itemize}
\end{prop}
\begin{proof}
First let us assume that $V(s)$ is reducible and $(\delta_1,\delta_2)\notin S^+_{0}$.
Because $V(s)$ is reducible, there exist two continuous characters $\delta_3,\delta_4:G_K\rightarrow \mathcal{O}_{E}^{\times}$ such that $W(V(s))$ sits in the following short exact sequence
\begin{equation*}
0\rightarrow W(E(\delta_3))\rightarrow W(V(s))\rightarrow W(E(\delta_4))\rightarrow 0.
\end{equation*}
Then the fact that $\mathrm{Hom}(W(E(\delta_3)), W(\delta_1))=0$ (this follows from the fact that (the 
slope of $W(E(\delta_3)))=0<$ (the slope of $W(\delta_1))$ and from Proposition $\ref{28}$) implies that the natural map $W(E(\delta_3))\hookrightarrow W(\delta_2)$ is an inclusion.
So we have $\tilde{\delta}_3=\delta_2\prod_{\sigma:K\hookrightarrow E}\sigma(x)^{k_{\sigma}}$ for some 
$\{k_{\sigma}\}_{\sigma}$ such that $k_{\sigma}\in\mathbb{Z}_{\ge 0}$ (here $\tilde{\delta_3}:=\delta_3\circ\mathrm{rec}_K:K^{\times}\rightarrow E^{\times}$). Then, as in the proof of the previous theorem, we can see that 
$s$ corresponds to $[(a_{\sigma}e_{\sigma})_{\sigma}]\in {S'}^{\acute{\mathrm{e}}t}(\delta_1,\delta_2)$ (or $[\overline{(a_{\sigma}e_{\sigma})_{\sigma}}]\in {S'}^{\acute{\mathrm{e}}t}(\delta_1,\delta_2)$) such that $(\sum_{\sigma,a_{\sigma}\not=0}w(\delta_1/\delta_2)_{\sigma}) +
e_K \mathrm{val}_p(\delta_2(\pi_K))\le (\sum_{\sigma:K\hookrightarrow E}k_{\sigma})
+ e_K\mathrm{val}_p(\delta_2(\pi_K))=[K:\mathbb{Q}_p]$( slope of $W(E(\delta_3)))=0$ (for some lifting $[(a_{\sigma}e_{\sigma})_{\sigma}]$ of $[\overline{(a_{\sigma}e_{\sigma})_{\sigma}}]$).
On the other hand, we have $0\le (\sum_{\sigma,a_{\sigma}\not=0}w(\delta_1/\delta_2)_{\sigma}) +
e_K \mathrm{val}_p(\delta_2(\pi_K))$ because $s\in {S'}^{\acute{\mathrm{e}}t}(\delta_1,\delta_2)$. So $s\in {S'}^{\mathrm{ord}}(\delta_1,\delta_2)$.

Next if $(\delta_1,\delta_2)\in S^+_{0}$, then $W(\delta_1)$ and $W(\delta_2)$ are pure of slope zero, hence $V(s)$ is reducible as  $E$-representation.

Next let us  assume that $(\delta_1,\delta_2)\notin S^+_{0}$ and $s\in {S'}^{\mathrm{ord}}(
\delta_1,\delta_2)$. If we put $s=[(a_{\sigma}e_{\sigma})_{\sigma}]$ (or $[\overline{(a_{\sigma}e_{\sigma})_{\sigma}}]$) such that $(\sum_{\sigma,a_{\sigma}\not= 0}w(\delta_1/\delta_2)_{\sigma})+ e_K \mathrm{val}_p(\delta_2(\pi_K))=0$ ( for some lift $[(a_{\sigma}e_{\sigma})_{\sigma}]$ of $[\overline{(a_{\sigma}e_{\sigma})_{\sigma}}]$) then, by the proof of Theorem $\ref{34}$, we have 
\begin{itemize}
\item[]$s\in\mathrm{Ker}(\mathrm{H}^1(G_K,W(\delta_1/\delta_2))\rightarrow\mathrm{H}^1(G_K, W(\prod_{\sigma,a_{\sigma}\not= 0}\sigma(x)^{-w(\delta_1/\delta_2)_{\sigma}}\delta_1/\delta_2)))$.
\end{itemize}
Then the natural inclusion $W(\prod_{\sigma,a_{\sigma}\not= 0}\sigma(x)^{w(\delta_1/\delta_2)_{\sigma}}\delta_2)\hookrightarrow W(\delta_2)$ factors through $W(\prod_{\sigma,a_{\sigma}\not= 0}\allowbreak\sigma(x)^{w(\delta_1/\delta_2)_{\sigma}}\delta_2)\hookrightarrow W(s)$. We take the saturation $W(\delta_3)$ of this map and write the cokernel by $W(\delta_4)$. Then we claim that the natural inclusion $W(\prod_{\sigma,a_{\sigma}\not= 0}\sigma(x)^{w(\delta_1/\delta_2)_{\sigma}}\delta_2)\hookrightarrow W(\delta_3)$ is an isomorphism.
If this is not an isomorphism, then by Proposition $\ref{28}$, $W(\delta_3)$ must be pure of  negative slope because 
$W(\prod_{\sigma,a_{\sigma}\not= 0}\sigma(x)^{w(\delta_1/\delta_2)_{\sigma}}\delta_2)$ is pure of slope zero by assumption. Then $W(s)$ is not  pure of slope zero, which contradicts to the assumption 
that $s\in {S'}^{\acute{\mathrm{e}}t}(\delta_1,\delta_2)$. So $W(\prod_{\sigma,a_{\sigma}\not= 0}\sigma(x)^{w(\delta_1/\delta_2)_{\sigma}}\delta_2)\hookrightarrow W(\delta_3)$ is an isomorphism. Then both $W(\delta_3)$ and $W(\delta_4)$ are pure of slope zero.
So $W(s)$ is reducible as $E$-representation.
So we have finished the proof of this proposition.
\end{proof}

\subsection{The conditions for $V(s)=V(s')$.}

Next let us discuss when two split trianguline $E$-representations $V(s)$, $V(s')$ are 
isomorphic for different parameters $s$, $s'$. Unfortunately we cannot solve this problem in the case 
where $\delta_1/\delta_2=\prod_{\sigma:K\hookrightarrow E}\sigma(x)^{k_{\sigma}}$ for some $\{k_{\sigma}\}_{\sigma}$ such that 
$k_{\sigma}\in\mathbb{Z}$ for any $\sigma$ and $s\in {S'}^{\acute{\mathrm{e}}t}(\delta_1,\delta_2)$. We can solve this problem in all the other cases. 
The following theorem is a generalization of $\cite[\mathrm{Proposition}\,4.9]
{Co07a}$. 
\begin{thm}$\label{36}$
\begin{itemize}
\item[$\mathrm{(1)}$]
Let $V(s)$, $V(s')$ be two dimensional trianguline $E$-representations 
associated to $s\in S(\delta_1,\delta_2)\setminus {S'}^{non-\acute{\mathrm{e}}t}(\delta_1,\delta_2)$, $s'\in S(\delta_3,\delta_4)\setminus {S'}^{non-\acute{\mathrm{e}}t}(\delta_3,\delta_4)$ for some 
$(\delta_1,\delta_2)$, $(\delta_3,\delta_4)\in S^+$.
Moreover we assume that $s\notin {S'}^{\acute{\mathrm{e}}t}(\delta_1,\delta_2)$.
Then $V(s)\isom V(s')$ if and only if $(\delta_1,\delta_2)=(\delta_3,\delta_4)$ and $s=s'\in S(\delta_1,\delta_2)$.
\item[$\mathrm{(2)}$]
Let $s=[(a_{\sigma}e_{\sigma})_{\sigma}]\in {S'}^{\acute{\mathrm{e}}t}(\delta_1,\delta_2)$ where 
$\delta_1/\delta_2\not= \prod_{\sigma:K\hookrightarrow E}\sigma(x)^{k_{\sigma}}$ for any $\{k_{\sigma}\}_{\sigma}$ such that $k_{\sigma}\in \mathbb{Z}$ for any $\sigma$.
Then there exists unique $((\delta_3,\delta_4), s')\not= ((\delta_1,\delta_2), s)$ $($here $(\delta_3,\delta_4)\in S^+$ and $s'\in S(\delta_3,\delta_4)\setminus {S'}^{non-\acute{\mathrm{e}}t}(\delta_3,\delta_4))$ such that 
$V(s)\isom V(s')$. Such a $((\delta_3,\delta_4), s')$ satisfies the following:
\begin{itemize}
\item[$\mathrm{(i)}$]$\delta_3=\delta_2\prod_{\sigma, a_{\sigma}\not= 0}\sigma(x)^{w(\delta_1/\delta_2)_{\sigma}}$, $\delta_4=\delta_1\prod_{\sigma, a_{\sigma}\not= 0}\sigma(x)^{-w(\delta_1/\delta_2)_{\sigma}}$.
\item[$\mathrm{(ii)}$]$s'=[(b_{\sigma}e_{\sigma})_{\sigma}]\in {S'}^{\acute{\mathrm{e}}t}(\delta_3,\delta_4)$ satisfies that $\{\sigma:K\hookrightarrow E|
b_{\sigma}\not= 0\}=\{\sigma:K\hookrightarrow E|a_{\sigma}\not= 0\}$.
\end{itemize}
\end{itemize}
\end{thm}
\begin{proof}
First, we prove (1).
If $V(s)\isom V(s')$, then we have the following short exact sequence
\begin{equation*}
0\rightarrow W(\delta_3)\rightarrow W(s)\rightarrow W(\delta_4)\rightarrow 0.
\end{equation*}
We consider the natural map $f:W(\delta_3)\rightarrow W(\delta_2)$ as in the proof of Theorem $\ref{34}$. Then we have one of the following: $\mathrm{Ker}(f)\isom W(\delta_3)$ or $\mathrm{Im}(f)\isom
W(\delta_3)$. 

If $\mathrm{Ker}(f)\isom W(\delta_3)$, then 
we have the following commutative diagram
\begin{center}
$\begin{CD}
0 @>>> W(\delta_3) @>>> W(s) @>>> W(\delta_4) @>>> 0  \\
@.     @VVV   @| @ VVV @. \\
0 @>>> W(\delta_1) @>>> W(s) @>>> W(\delta_2) @>>> 0. 
\end{CD}$
\end{center}

Then we claim that $W(\delta_3)\rightarrow W(\delta_1)$ and $W(\delta_4)
\rightarrow W(\delta_2)$ are both isomorphisms. 
We can see that $W_e(\delta_3)\isom W_e(\delta_1)$, $W_e(\delta_4)\isom W_e(
\delta_2)$ are isomorphisms in the same way as the proof of  Proposition $\ref{28}$.
For the $W^+_{\mathrm{dR}}$-part, by the snake lemma of the above diagram, we have an isomorphism $\mathrm{Ker}(W^+_{\mathrm{dR}}(\delta_4)\rightarrow W^+_{\mathrm{dR}}(\delta_2))\isom \mathrm{Cok}(W^+_{\mathrm{dR}}(\delta_3)\rightarrow 
W^+_{\mathrm{dR}}(\delta_1))$ of $B^+_{\mathrm{dR}}$-modules. But because the former is a torsion free $B^+_{\mathrm{dR}}$-module and the latter is a torsion $B^+_{\mathrm{dR}}$-module, so this must be zero. So $W(\delta_3)\rightarrow W(\delta_1)$ and $W(\delta_4)
\rightarrow W(\delta_2)$ are both isomorphisms. So $(\delta_1,\delta_2)=(\delta_3,\delta_4)$ and $s=s'\in S(\delta_1,\delta_2)$.

If $\mathrm{Im}(f)\isom W(\delta_3)$, then $\delta_3=\delta_2\prod_{\sigma:K\hookrightarrow E}\sigma(x)^{k_{\sigma}}$ for some $\{k_{\sigma}\}_{\sigma}$ such 
that $k_{\sigma}\in\mathbb{Z}_{\ge 0}$ for any $\sigma$. Then, by the proof of Theorem 
$\ref{34}$,
we can see that $s\in S'(\delta_1,\delta_2)$. 
This is a contradiction. So we have finished the proof of (1) of this theorem.

Next, we prove (2). Let us assume that $V(s)\isom V(s')$ for some $s'(\not= s)\in S(\delta_3,\delta_4)$.
Then, by the proof of (1) above and the assumption $s\not= s'$, we can see that  $W(\delta_3)\isom \mathrm{Im}(f:W(\delta_3)\rightarrow W(\delta_2))$.
So we have $\delta_3=\prod_{\sigma:K\hookrightarrow E}\sigma(x)^{k_{\sigma}}\delta_2$ for some $\{k_{\sigma}\}_{\sigma}$ such that $k_{\sigma}\in\mathbb{Z}_{\ge 0}$ for any $\sigma$. So we get $s=[(a_{\sigma}e_{\sigma})_{\sigma}]\in {S'}^{\acute{\mathrm{e}}t}(\delta_1,\delta_2)$ and $w(\delta_1/\delta_2)_{\sigma}\in \{1,2,\cdots ,k_{\sigma}\}$ for any 
$\sigma$ such that $a_{\sigma}\not= 0$. Then, by the proof of Theorem $\ref{34}$, we can see that the natural inclusion 
$W(\prod_{\sigma,a_{\sigma}\not= 0}\sigma(x)^{w(\delta_1/\delta_2)_{\sigma}}\delta_2)\hookrightarrow W(\delta_2)$ factors through an inclusion $g:W(\prod_{\sigma,a_{\sigma}\not= 0}\sigma(x)^{w(\delta_1/\delta_2)_{\sigma}}\allowbreak\delta_2)\hookrightarrow W(s)$. Because we have $\mathrm{Hom}(W(\prod_{\sigma:K\hookrightarrow E}\sigma(x)^{k_{\sigma}}\delta_2),W(\delta_1))=0$ by the assumption on $(\delta_1,\delta_2)$ and by Proposition $\ref{28}$,
we can see that $g \circ\iota:W(\prod_{\sigma:K\hookrightarrow E}\sigma(x)^{k_{\sigma}}\delta_2)\hookrightarrow W(s)$ is equal to the given inclusion $W(\prod_{\sigma:K\hookrightarrow E}\sigma(x)^{k_{\sigma}}\delta_2)=W(\delta_3)\hookrightarrow W(s)$ (here $\iota:W(\prod_{\sigma:K\hookrightarrow E}\sigma(x)^{k_{\sigma}}\delta_2)\hookrightarrow W(\prod_{\sigma,a_{\sigma}\not= 0}\allowbreak\sigma(x)^{w(\delta_1/\delta_2)_{\sigma}}\delta_2)$ is the natural inclusion). Because $W(\delta_3)$ is saturated in $W(s)$ so $W(\delta_3)\hookrightarrow W(\prod_{\sigma,a_{\sigma}\not= 0}\sigma(x)^{w(\delta_1/\delta_2)_{\sigma}}\delta_2)$ must be isomorphic. So we have $k_{\sigma}=w(\delta_1/\delta_2)_{\sigma}$ 
for $\sigma$ such that $a_{\sigma}\not= 0$ and $k_{\sigma}=0$ for other $\sigma$. Then we have $\delta_3=\prod_{\sigma,a_{\sigma}\not= 0}\sigma(x)^{w(\delta_1/\delta_2)_{\sigma}}\delta_2$ and $\delta_4=\prod_{\sigma,a_{\sigma}\not= 0}\sigma(x)^{-w(\delta_1/\delta_2)_{\sigma}}\delta_1$. We have $w(\delta_3/\delta_4)_{\sigma}=w(\delta_3)_{\sigma}-w(\delta_4)_{\sigma}=w(\delta_1/\delta_2)_{\sigma}+w(\delta_2)_{\sigma}-(-w(\delta_1/\delta_2)_{\sigma}+w(\delta_1)_{\sigma})=w(\delta_1/\delta_2)_{\sigma}$ for any $\sigma$ such that $a_{\sigma}\not= 0$ and 
$w(\delta_3/\delta_4)_{\sigma}=w(\delta_2)_{\sigma}-w(\delta_1)_{\sigma}=
-w(\delta_1/\delta_2)_{\sigma}$ for other $\sigma$. 
If we replace $s$ by $s'$ and replace $s'$ by $s$ in the above argument, 
we can see that $s'=[(b_{\sigma}e_{\sigma})_{\sigma}]\in {S'}^{\acute{\mathrm{e}}t}(\prod_{\sigma,a_{\sigma}\not= 0}\sigma(x)^{w(\delta_1/\delta_2)_{\sigma}}\delta_2,\prod_{\sigma,a_{\sigma}\not= 0}\sigma(x)^{-w(\delta_1/\delta_2)_{\sigma}}\delta_1)$ satisfies the condition that  $\{\sigma:K\hookrightarrow E| b_{\sigma}\not= 0\}=\{\sigma:K\hookrightarrow E| a_{\sigma}\not= 0\}$. If $s''=[(b'_{\sigma}e_{\sigma})_{\sigma}]\in {S'}^{\acute{\mathrm{e}}t}(\prod_{\sigma,a_{\sigma}\not= 0}\sigma(x)^{w(\delta_1/\delta_2)_{\sigma}}\delta_2,\allowbreak\prod_{\sigma,a_{\sigma}\not= 0}\sigma(x)^{-w(\delta_1/\delta_2)_{\sigma}}\delta_1)$ is another element such that $V(s)\isom V(s'')$ and $s\not= s''$, then we can see that $s'=s''$ in the same way as in the proof of (1) above. So such an $s'$ is unique.

Next let us take any $s=[a_{\sigma}e_{\sigma}]\in {S'}^{\acute{\mathrm{e}}t}(\delta_1,\delta_2)$.
Then, by Lemma $\ref{32}$, we have 
\begin{itemize}
\item[]$s\in \mathrm{Ker}(\mathrm{H}^1(G_K,W(\delta_1/\delta_2))\rightarrow \mathrm{H}^1(G_K,W(\prod_{\sigma,a_{\sigma}\not= 0}\sigma(x)^{-w(\delta_1/\delta_2)_{\sigma}}\delta_1/\delta_2)))$. 
\end{itemize}
So the natural inclusion $W(\prod_{\sigma,a_{\sigma}\not= 0}\sigma(x)^{w(\delta_1/\delta_2)_{\sigma}}\delta_2)\hookrightarrow W(\delta_2)$ factors through $W(\prod_{\sigma,a_{\sigma}\not= 0}\allowbreak\sigma(x)^{w(\delta_1/\delta_2)_{\sigma}}\delta_2)\hookrightarrow W(s)$. Then we can see that $W(\prod_{\sigma,a_{\sigma}\not= 0}\sigma(x)^{w(\delta_1/\delta_2)_{\sigma}}\delta_2)\hookrightarrow W(s)$ is saturated (if not then we have $s=[a'_{\sigma}e_{\sigma}]\in {S'}^{\acute{\mathrm{e}}t}(\delta_1,\delta_2)$ such that $\{\sigma:K\hookrightarrow E|a'_{\sigma}\not= 0\}\varsubsetneqq \{\sigma:K\hookrightarrow E| a_{\sigma}\not= 0\}$, which can not happen under the assumption on $(\delta_1,\delta_2)$ in this theorem). If we write the cokernel of this by $W(\delta_4)$, then we have $\delta_4=\prod
_{\sigma,a_{\sigma}\not= 0}\sigma(x)^{-w(\delta_1/\delta_2)_{\sigma}}\delta_1$ and we have the following short exact sequence
\begin{equation*}
0\rightarrow W(\delta_3)\rightarrow W(s)\rightarrow W(\delta_4)\rightarrow 0.
\end{equation*}
This extension corresponds to some $s'=[b_{\sigma}e_{\sigma}]\in {S'}^{\acute{\mathrm{e}}t}(\delta_3,\delta_4)$. So we have $V(s)\isom V(s')$.
We have finished the proof of this theorem.
\end{proof} 

\begin{rem}
When $K=\mathbb{Q}_p$, the exceptional cases where $s\in {S'}^{\acute{\mathrm{e}}t}(\delta_1,\delta_2)$ and $\delta_1/\delta_2=\prod_{\sigma:K\hookrightarrow E}\sigma(x)^{k_{\sigma}}$ do not appear because ${S'}^{\acute{\mathrm{e}}t}(\delta_1,\delta_2)$ is empty set in this case.
\end{rem}

\section{Classification of two dimensional potentially semi-stable 
split trianguline $E$-representations.}
In this final section, we classify all the two dimensional potentially semi-stable split trianguline $E$-representations.
We explicitly describe the $E$-filtered $(\varphi,N,G_K)$-modules associated 
to potentially semi-stable split trianguline $E$-representations.

\subsection{de Rham split trianguline $E$-representations.}
Before classifying two dimensional de Rham split trianguline 
$E$-representations, we determine all the rank one de Rham $E$-$B$-pairs.
\begin{lemma}$\label{38}$
Let $W(\delta)$ be a rank one $E$-$B$-pair.
Then the following conditions are equivalent:
\begin{itemize}
\item[$\mathrm{(1)}$]$W(\delta)$ is Hodge-Tate.
\item[$\mathrm{(2)}$]$W(\delta)$ is de Rham.
\item[$\mathrm{(3)}$]$\delta=\tilde{\delta}\prod_{\sigma:K\hookrightarrow E}\sigma(x)^{k_{\sigma}}$ such that $\tilde{\delta}:K^{\times}\rightarrow E^{\times}$ is a locally constant character and $k_{\sigma}\in\mathbb{Z}$ for any $\sigma$.
\end{itemize}
\end{lemma}
\begin{proof}
It is easy to see the implication 
$(3)\Rightarrow (2) \Rightarrow (1)$. 
We prove that $(1)$ implies $(3)$.
So we assume that $W(\delta)$ is Hodge-Tate with generalized 
Hodge-Tate weight $\{k_{\sigma}\}_{\sigma}$ such that 
$k_{\sigma}\in \mathbb{Z}$ for any $\sigma$.
We write $W(\delta)\isom W(E(\delta_0))\otimes W_0^{\otimes i}$ for 
some continuous character $\delta_0:G_K\rightarrow \mathcal{O}_E^{\times}$ and $i\in\mathbb{Z}$ as in 
 Theorem $\ref{15}$. Beause $W_0:=(W_{e,0}, W^+_{\mathrm{dR},0})$ satisfies  $W^+_{\mathrm{dR}}\isom B^+_{\mathrm{dR}}\otimes_{\mathbb{Q}_p}E$, so we have $W(\delta)^+_{\mathrm{dR}}\isom B^+_{\mathrm{dR}}\otimes_{\mathbb{Q}_p}E(\delta_0)$. Then $E(\delta_0\prod_{\sigma:K\hookrightarrow E}\sigma(\chi_{\mathrm{LT}})^{-k_{\sigma}})$ is Hodge-Tate with all Hodge-Tate weight zero. Then, by Sen's theorem $\cite[\mathrm{Proposition}\,5.24]{Be02}$, $\delta_0\prod_{\sigma:K\hookrightarrow E}\sigma(\chi_{\mathrm{LT}})^{-k_{\sigma}}:G_K\rightarrow \mathcal{O}_E^{\times}$ is a potentially unramified character. Because $\sigma(\chi_{\mathrm{LT}}):G_K\rightarrow \mathcal{O}_E^{\times}$ corresponds to the character $K^{\times}\rightarrow E^{\times}:\pi_K\mapsto 1, u\mapsto \sigma(u)$ for $u\in\mathcal{O}_K^{\times}$, we have $\delta=\tilde{\delta}\prod_{\sigma:K\hookrightarrow E}\sigma(x)^{k_{\sigma}}$ for some locally constant character $\tilde{\delta}
:K^{\times}\rightarrow E^{\times}$.
So we have finished the proof of the lemma.
\end{proof}

Next let $V(s)$ be a potentially semi-stable split trianguline $E$-representation 
such that $s\in S(\delta_1,\delta_2)\setminus {S'}^{non-\acute{\mathrm{e}}t}(\delta_1,\delta_2)$ for some $(\delta_1,\delta_2)\in S^+$.
Then we have the following short exact sequence
\begin{equation*}
0\rightarrow W(\delta_1)\rightarrow W(s)\rightarrow W(\delta_2)\rightarrow 0.
\end{equation*}
Then, as in the case of usual $E$-representations, we can see that 
$W(\delta_1)$ and $W(\delta_2)$ are also de Rham, in particular 
$W(\delta_1/\delta_2)$ is de Rham.
Next we compute $\mathrm{H}^1_e(G_K,W(\delta))$, $\mathrm{H}^1_f(G_K, W(\delta))$ and $\mathrm{H}^1_g(G_K, W(\delta))$ when $W(\delta)$ is de Rham.
\begin{lemma}$\label{39}$
Let $W(\delta)$ be a rank one de Rham $E$-$B$-pair where $\delta:=\tilde{\delta}\prod_{\sigma:K\hookrightarrow E}\sigma(x)^{k_{\sigma}}$ such that $\tilde{\delta}:K^{\times}\rightarrow E^{\times}$ is a locally constant character. 
Then we have a canonical isomorphism
\begin{itemize}
\item[]$D^K_{\mathrm{dR}}(W(\delta))/(D^K_{\mathrm{cris}}(W(\delta))^{\varphi=1}+
\mathrm{Fil}^0D^K_{\mathrm{dR}}(W(\delta)))\isom\mathrm{H}^1_e(G_K, W(\delta))$
\end{itemize}
and $\mathrm{dim}_E\mathrm{H}^1_e(G_K, W(\delta))$ is equal to 
\begin{itemize}
\item[$\mathrm{(1)}$]$\sharp\{\sigma:K\hookrightarrow E| k_{\sigma}\in \mathbb{Z}_{\ge 1}\}$ when $\tilde{\delta}$ is not the trivial character,
\item[$\mathrm{(2)}$]$\sharp\{\sigma:K\hookrightarrow E| k_{\sigma}\in \mathbb{Z}_{\ge 1}\}-1$ when $\tilde{\delta}$ is the trivial character.
\end{itemize}
\end{lemma}
\begin{proof}
For general $\delta$, we have the following short exact sequence
\begin{equation*}
0\rightarrow D^K_{\mathrm{dR}}(W(\delta))/(D^K_{\mathrm{cris}}(W(\delta))^{\varphi=1}+
\mathrm{Fil}^0D^K_{\mathrm{dR}}(W(\delta)))\rightarrow\mathrm{H}^1(G_K, W(\delta))
\end{equation*}
\begin{equation*}
\rightarrow \mathrm{Ker}(\mathrm{H}^1(G_K, W_e(\delta))\oplus\mathrm{H}^1(G_K, W^+_{\mathrm{dR}}(\delta))\rightarrow \mathrm{H}^1(G_K, W_{\mathrm{dR}}(\delta)))\rightarrow 0.
\end{equation*}
When $W(\delta)$ is de Rham, we proved in Lemma $\ref{21}$ that 
$\mathrm{H}^1(G_K, W^+_{\mathrm{dR}}(\delta))\rightarrow \mathrm{H}^1(G_K, W_{\mathrm{dR}}(\delta))$ is injective. So we have 
\begin{align*}
\mathrm{H}^1_e(G_K, W(\delta)) &=
\mathrm{Ker}(\mathrm{H}^1(G_K, W(\delta))\rightarrow \mathrm{H}^1(G_K, W_e(\delta)))  \\
 &\isom\mathrm{Ker}(\mathrm{H}^1(G_K, W(\delta))\rightarrow(\mathrm{H}^1(
G_K, W_e(\delta))\oplus\mathrm{H}^1(G_K, W^+_{\mathrm{dR}}(\delta)))) \\ 
 &=D^K_{\mathrm{dR}}(W(\delta))/(D^K_{\mathrm{cris}}(W(\delta))^{\varphi=1}+
\mathrm{Fil}^0D^K_{\mathrm{dR}}(W(\delta))).
\end{align*}
We can compute the dimension of $\mathrm{H}^1_e(G_K, W(\delta))$ by using 
this isomorphism.
\end{proof}

\begin{lemma}$\label{40}$
Let $W(\delta)$ be a rank one de Rham $E$-$B$-pair. 
Then we have 
\begin{itemize}
\item[$\mathrm{(1)}$]$\mathrm{H}^1_g(G_K, W(\delta))=\mathrm{H}^1_f(G_K, W(\delta))$ and $\mathrm{dim}_E\mathrm{H}^1_f(G_K,W(\delta))=\mathrm{dim}_E\mathrm{H}^1_e(G_K, W(\delta))+1$ when $\delta=\prod_{\sigma:K\hookrightarrow E}\sigma(x)^{k_{\sigma}}$,
\item[$\mathrm{(2)}$]$\mathrm{H}^1_f(G_K,W(\delta))=\mathrm{H}^1_e(G_K, W(\delta))$ and $\mathrm{dim}_E\mathrm{H}^1_g(G_K,W(\delta))=\mathrm{dim}_E\mathrm{H}^1_f(G_K, W(\delta))+1$ when $\delta=|N_{K/\mathbb{Q}_p}(x)|\prod_{\sigma:K\hookrightarrow E}\sigma(x)^{k_{\sigma}}$,
\item[$\mathrm{(3)}$]$\mathrm{H}^1_g(G_K, W(\delta))=\mathrm{H}^1_f(G_K, W(\delta))=\mathrm{H}^1_e(G_K, W(\delta))$ otherwise.
\end{itemize}
\end{lemma}
\begin{proof}
This follows from calculation by using Lemma $\ref{22}$ and 
Lemma $\ref{25}$.
\end{proof}
For the explicit description of potentially semi-stable split 
trianguline $E$-representations, we fix an extension whose class is contained in $\mathrm{H}^1_f(G_K,W(\delta))\setminus \mathrm{H}^1_e(G_K,W(\delta))$ when $\delta=\prod_{\sigma}\sigma(x)^{k_{\sigma}}$ and fix an extension whose class is contained in $\mathrm{H}^1_g(G_K,W(\delta))\setminus \mathrm{H}^1_f(G_K,W(\delta))$ when $\delta=|N_{K/\mathbb{Q}_p}(x)|\prod_{\sigma}\sigma(x)^{k_{\sigma}}$ as follows.
First we treat the case where $\delta=\prod_{\sigma}\sigma(x)^{k_{\sigma}}$.
In this case, we have $W(\delta)=(B_e\otimes_{\mathbb{Q}_p}E, \oplus_{\sigma}
t^{k_{\sigma}}B^+_{\mathrm{dR}}\otimes_{K,\sigma}E)$ by Lemma $\ref{26}$.
We fix an element $a_{\mathrm{cris}}\in\mathcal{O}_{\mathbb{Q}_p^{\mathrm{un}}}^{\times}$ such that 
$\mathrm{Frob}_{\mathbb{Q}_p}(a_{\mathrm{cris}})=a_{\mathrm{cris}}+1$.
We define an extension $W(s(\{k_{\sigma}\}_{\sigma})):=(W_e(s(\{k_{\sigma}\}_{\sigma})),W^+_{\mathrm{dR}}(s(\{k_{\sigma}\}_{\sigma})),\iota)$ whose class $s(\{k_{\sigma}\}_{\sigma}):=[W(s(\{k_{\sigma}\}_{\sigma}))]$ is contained in $\mathrm{H}^1_f(G_K, \allowbreak W(\prod_{\sigma}\sigma(x)^{k_{\sigma}}))\setminus \mathrm{H}^1_e(G_K, W(\prod_{\sigma}\sigma(x)^{k_{\sigma}}))$ as follows:
\begin{itemize}
\item[(1)] $W_e(s(\{k_{\sigma}\}_{\sigma})):=
(B_e\otimes_{\mathbb{Q}_p}E)e_1\oplus (B_e\otimes_{\mathbb{Q}_p}E)e_{\mathrm{cris}}$ such that $g(e_1)=e_1$, $g(e_{\mathrm{cris}})=e_{\mathrm{cris}}+f\mathrm{deg}(g)e_1$ for any $g\in G_K$, here $f=[K_0:\mathbb{Q}_p]$.
\item[(2)]$W^+_{\mathrm{dR}}(s(\{k_{\sigma}\}_{\sigma})):=(\oplus_{\sigma}t^{k_{\sigma}}B^+_{\mathrm{dR}}\otimes_{K,\sigma}E)e_1\oplus (B^+_{\mathrm{dR}}\oplus_{\mathbb{Q}_p}E)e_{\mathrm{dR}}$ such that 
$g(e_1)=e_1$ and $g(e_{\mathrm{dR}})=e_{\mathrm{dR}}$ for any $g\in G_K$.
\item[(3)]$\iota:B_{\mathrm{dR}}\otimes_{B_e}W_e(s(\{k_{\sigma}\}_{\sigma}))\isom B_{\mathrm{dR}}\otimes_{B^+_{\mathrm{dR}}}W^+_{\mathrm{dR}}(s(\{k_{\sigma}\}_{\sigma}))$ is the isomorphism defined by $\iota(e_1)=e_1$ and $\iota(e_{\mathrm{cris}})=e_{\mathrm{dR}}+a_{\mathrm{cris}}e_1$.
\end{itemize}
We can easily see that $s(\{k_{\sigma}\}_{\sigma}):=[W(s(\{k_{\sigma}\}_{\sigma}))]$ is contained in $\mathrm{H}^1_f(G_K, W(\delta))\setminus \mathrm{H}^1_e(G_K,W(\delta))$. 

Next we treat the case where $\delta=|N_{K/\mathbb{Q}_p}(x)|\prod_{\sigma}\sigma(x)^{k_{\sigma}}$. We define a continuous one cocycle $c:G_K\rightarrow \mathbb{Q}_p(\chi)$
by $g(\mathrm{log}([\tilde{p}]))=\mathrm{log}([\tilde{p}])+c(g)t$ for any $g\in G_K$. In this case, we define an extension $W(s(\{k_{\sigma}\}_{\sigma})):=(W_e(s(\{k_{\sigma}\}_{\sigma})),W^+_{\mathrm{dR}}(s(\{k_{\sigma}\}_{\sigma})),\iota)$ whose class $s(\{k_{\sigma}\}_{\sigma}):=[W(s(\{k_{\sigma}\}_{\sigma}))]$
is in $\mathrm{H}^1_g(G_K, W(|N_{K/\mathbb{Q}_p}(x)|\prod_{\sigma}\sigma(x)^{k_{\sigma}}))\setminus \mathrm{H}^1_f(G_K, W(|N_{K/\mathbb{Q}_p}(x)|\prod_{\sigma}\sigma(x)^{k_{\sigma}}))$ as follows:
\begin{itemize}
\item[(1)]$W_e(s(\{k_{\sigma}\}_{\sigma})):=(B_e\otimes_{\mathbb{Q}_p}E(\chi))e_1
\oplus (B_e\otimes_{\mathbb{Q}_p}E)e_{\mathrm{cris}}$ such that $g(e_1)=\chi(g)e_1, g(e_{\mathrm{cris}})=e_{\mathrm{cris}}+c(g)e_1$ for any $g\in G_K$.
\item[(2)]$W^+_{\mathrm{dR}}(s(\{k_{\sigma}\}_{\sigma})):=(\oplus_{\sigma}
t^{k_{\sigma}-1}B^+_{\mathrm{dR}}\otimes_{K,\sigma}E(\chi))e_1\oplus (B^+_{\mathrm{dR}}\otimes_{\mathbb{Q}_p}E)e_{\mathrm{dR}}$ such that $g(e_1)=\chi(g)e_1, 
g(e_{\mathrm{dR}})=e_{\mathrm{dR}}$ for any $g\in G_K$.
\item[(3)]$\iota:B_{\mathrm{dR}}\otimes_{B_e}W_e(s(\{k_{\sigma}\}_{\sigma}))
\isom B_{\mathrm{dR}}\otimes_{B^+_{\mathrm{dR}}}W^+_{\mathrm{dR}}(s(\{k_{\sigma}\}_{\sigma}))$ is the isomorphism defined by $\iota(e_1)=e_1, \iota(e_{\mathrm{cris}})=e_{\mathrm{dR}}+\frac{\mathrm{log}[(\tilde{p})]}{t}e_1$.
\end{itemize}
We can easily see that $s(\{k_{\sigma}\}_{\sigma}):=[W(s(\{k_{\sigma}\}_{\sigma}))]$ is contained in $\mathrm{H}^1_g(G_K, W(\delta))\setminus \mathrm{H}^1_f(G_K,W(\delta))$.

By using these classes, we can determine all the potentially cristalline split 
trianguline $E$-representations and all the potentially semi-stable split 
trianguline $E$-representations which are not potentially cristalline.
For this, let us put
\begin{itemize}
\item[$\bullet$]$S^{''}(\delta_1,\delta_2):=\{s\in S(\delta_1,\delta_2)|
s\in \mathbb{P}_E(\mathrm{H}^1_f(G_K,W(\delta_1/\delta_2))\setminus S^{'}(\delta_1,\delta_2)\}$ when $\delta_1/\delta_2=\prod_{\sigma}\sigma(x)^{k_{\sigma}}$, 
\item[$\bullet$]$S_{\mathrm{st}}(\delta_1,\delta_2):=\{s\in S(\delta_1,\delta_2)|s\in\mathbb{P}_E(\mathrm{H}^1_g(G_K,W(\delta_1/\delta_2))\setminus S^{'}(\delta_1,\delta_2))\}$ when $\delta_1/\delta_2=|N_{K/\mathbb{Q}_p}(x)|\prod_{\sigma}\sigma(x)^{k_{\sigma}}$. 
\end{itemize}
Then, by Lemma $\ref{40}$ and the above constructions, we have 
\begin{itemize}
\item[$\bullet$]$S^{''}(\delta_1,\delta_2)\isom s(\{k_{\sigma}\}_{\sigma})+
D^K_{\mathrm{dR}}(W(\delta_1/\delta_2)/(D^K_{\mathrm{cris}}(W(\delta_1/\delta_2))^{\varphi=1}+\mathrm{Fil}^0D^K_{\mathrm{dR}}(W(\delta_1/\delta_2)))$,
\item[$\bullet$]$S_{\mathrm{st}}(\delta_1,\delta_2)\isom s(\{k_{\sigma}\}_{\sigma})+
D^K_{\mathrm{dR}}(W(\delta_1/\delta_2)/(D^K_{\mathrm{cris}}(W(\delta_1/\delta_2))^{\varphi=1}+\mathrm{Fil}^0D^K_{\mathrm{dR}}(W(\delta_1/\delta_2)))$.
\end{itemize}
\begin{prop}$\label{41}$
Let $(\delta_1,\delta_2)\in S^+$ such that $W(\delta_1)$ and $W(\delta_2)$ are 
de Rham $E$-$B$-pairs. Let $V(s)$ be a split trianguline $E$-representation associated 
to $s\in S(\delta_1,\delta_2)\setminus {S'}^{non-\acute{\mathrm{e}}t}(\delta_1,\delta_2)$.
Then the following conditions  are equivalent$:$
\begin{itemize}
\item[$\mathrm{(1)}$]$V(s)$ is potentially cristalline.
\item[$\mathrm{(2)}$]$s\in S^{\acute{\mathrm{e}}t}_{\mathrm{cris}}(\delta_1,\delta_2):={S'}^{\acute{\mathrm{e}}t}(\delta_1,\delta_2)\sqcup S^{''}(\delta_1,\delta_2)$ when
$\delta_1/\delta_2=\prod_{\sigma}\sigma(x)^{k_{\sigma}}$.
\item[]$s\in S^{\acute{\mathrm{e}}t}_{\mathrm{cris}}(\delta_1,\delta_2):={S'}^{\acute{\mathrm{e}}t}(\delta_1,\delta_2)$ otherwise.
\end{itemize}
The following conditions are equivalent$:$
\begin{itemize}
\item[$\mathrm{(3)}$]$V(s)$ is potentially semi-stable and not potentially cristalline.
\item[$\mathrm{(4)}$]$\delta_1/\delta_2=|N_{K/\mathbb{Q}_p}(x)|\prod_{\sigma}\sigma(x)^{k_{\sigma}}$ and $s\in S_{\mathrm{st}}(\delta_1,\delta_2)$.
\end{itemize}
\end{prop}
\begin{proof}
This easily follows from Lemma $\ref{40}$.
We prove in the case $\delta_1/\delta_2=|N_{K/\mathbb{Q}_p}(x)|\prod_{\sigma}\sigma(x)^{k_{\sigma}}$. (Other cases can be proved in the same way.)
First we prove that (1) is equivalent to (2).
If $s\in {S'}^{\acute{\mathrm{e}}t}(\delta_1,\delta_2)$ then $V(s)$ is potentially cristalline 
by Remark $\ref{b}$. If $V(s)$ is potentially cristallline then 
$V(s)$ is de Rham. So $[W(V(s))]$ is contained in $\mathrm{H}^1_g(G_K,W(\delta_1/\delta_2))$.
If $[W(V(s))]\notin \mathrm{H}^1_f(G_K,W(\delta_1/\delta_2))$ then by the above remark, we can assume that $[W(V(s))]=s(\{k_{\sigma}\}_{\sigma})+x$ for some $x\in D^K_{\mathrm{dR}}(W(\delta_1/\delta_2)/(D^K_{\mathrm{cris}}(W(\delta_1/\delta_2))^{\varphi=1}+\mathrm{Fil}^0D^K_{\mathrm{dR}}(W(\delta_1/\delta_2)))$. From this and from the construction 
of $s(\{k_{\sigma}\}_{\sigma})$, we can easily see that $W(V(s))|_{G_L}\notin \mathrm{H}^1_f(G_L, W(\delta_1/\delta_2)|_{G_L})$ for any finite extension $L$ of $K$.
(Here, $W(\delta_1/\delta_2)|_{G_L}$ is the $E$-$B$-pair of $G_L$ obtained by restriction.) This implies that $V(s)$ is not potentially cristalline. 
This is a contradiction. So $W(V(s))$ is contained in $\mathrm{H}^1_f(G_K,W(\delta_1/\delta_2))$, i.e. $s\in {S'}^{\acute{\mathrm{e}}t}(\delta_1,\delta_2)$.
We can prove the equivalence between (3) and (4) in the same way.
\end{proof}
\begin{rem}
When $K=\mathbb{Q}_p$, $\mathcal{S}^{\mathrm{cris}}_+$ in $\cite[4.3]{Co07a}$ is equal to $\sqcup_{(\delta_1,\delta_2)\in S^+}S^{\acute{\mathrm{e}}t}_{\mathrm{cris}}(\delta_1,\delta_2)$ and $\mathcal{S}^{\mathrm{st}}_+$ is equal to 
$\sqcup_{(\delta_1,\delta_2)\in S^+}S_{\mathrm{st}}(\delta_1,\delta_2)$. 
Here, for $(\delta_1,\delta_2)\in S^+$ such that $\delta_1/\delta_2\not =|N_{K/\mathbb{Q}_p}(x)|\prod_{\sigma}\sigma(x)^{k_{\sigma}}$, we put $S_{\mathrm{st}}(\delta_1,\delta_2):=\phi$, empty set.
\end{rem}

\subsection{Potentially cristalline split trianguline
$E$-representations.}
In this subsection, we explicitly describe the filtered $(\varphi,G_K)$-modules of potentially cristalline split trianguline $E$-representations. 
First we define parameter spaces of potentially cristalline split 
trianguline $E$-representations. We put
\begin{itemize}
\item[$\bullet$]$T_{\mathrm{cris}}:=\{(\delta_1,\delta_2,\{k_{\sigma}\}_{\sigma})|\delta_1,\delta_2:K^{\times}\rightarrow E^{\times}$ locally constant characters, $k_{\sigma}\in\mathbb{Z}$ for any $\sigma$, such that $(\sum_{\sigma}k_{\sigma})+e_K \mathrm{val}_p(\delta_1(\pi_K))=-e_K\mathrm{val}_p(\delta_2(\pi_K))\ge 0\}$.
\end{itemize}
For any locally constant character $\delta:K^{\times}\rightarrow E^{\times}$ we define $n(\delta)\in\mathbb{Z}_{\ge 0}$ as the minimal $n\in\mathbb{Z}_{\ge 0}$ such that $\delta|_{1+\pi_K^n\mathcal{O}_K}$ is trivial. We write for any 
$\sigma$
\begin{itemize}
\item[$\bullet$]$G(\delta)_{\sigma}:=\sum_{\gamma\in\mathrm{Gal}(K_{n(\delta)}/K)}\gamma(\pi_{n(\delta)})\otimes \delta(\gamma^{-1})\in K_n\otimes_{K,\sigma}E$ when 
$n(\delta)\ge 1$,
\item[$\bullet$]$G(\delta)_{\sigma}:=1\in K\otimes_{K,\sigma}E$ when $n(\delta)=0$. 
 \end{itemize}
(Here $\pi_n$ is an element in $\bar{K}$ such that $[\pi_K](\pi_{n+1})=\pi_n$ for any $n\in\mathbb{N}$ and $[\pi_K](\pi_1)=0$, i.e. a system of $\pi_K^n$ torsion points of Lubin-Tate group associated to $\pi_K$, and $K_n:=K(\pi_n)$ for $n\in\mathbb{Z}_{\ge 1}$ and we identify $\mathrm{Gal}(K_n/K)\isom \mathcal{O}^{\times}_K/1+\pi_K^n\mathcal{O}_K$ via the Lubin-Tate character $\chi_{\mathrm{LT}}$.)
For any $(\delta_1,\delta_2,\{k_{\sigma}\}_{\sigma})\in T_{\mathrm{cris}}$, we define the parameter space $T_{\mathrm{cris}}^{\acute{\mathrm{e}}t}(\delta_1,\delta_2,\{k_{\sigma}\}_{\sigma})$ of weakly admissible filtrations as follows.
First, when $\delta_1\not= \delta_2$, we define 
\begin{itemize}
\item[$\bullet$]$T_{\mathrm{cris}}(\delta_1,\delta_2,\{k_{\sigma}\}_{\sigma}):=\mathbb{P}_E(\oplus_{\sigma, k_{\sigma}\ge 1}Ee_{\sigma})$,
\item[$\bullet$]$T_{\mathrm{cris}}^{\acute{\mathrm{e}}t}(\delta_1,\delta_2,\{k_{\sigma}\}_{\sigma}):=\{
[(a_{\sigma}e_{\sigma})_{\sigma}]\in T_{\mathrm{cris}}(\delta_1,\delta_2,\{k_{\sigma}\}_{\sigma})|(\sum_{\sigma,a_{\sigma}\not= 0}k_{\sigma})+e_K \mathrm{val}_p(\delta_2(\pi_K))\ge 0\}$.
\end{itemize}
When $\delta_1=\delta_2$, we define 
\begin{itemize}
\item[$\bullet$] $T_{\mathrm{cris}}'(\delta_1,\delta_2,\{k_{\sigma}\}_{\sigma}):=\mathbb{P}_E(\oplus_{\sigma, k_{\sigma}\ge 1}Ee_{\sigma}/\Delta(E))$,
\item[$\bullet$] ${T'}_{\mathrm{cris}}^{\acute{\mathrm{e}}t}(\delta_1,\delta_2,\{k_{\sigma}\}_{\sigma}):=
\{[(\overline{a_{\sigma}e_{\sigma}})_{\sigma}]\in {T'}_{\mathrm{cris}}(\delta_1,\delta_2,\{k_{\sigma}\}_{\sigma})| (\sum_{\sigma,a_{\sigma}\not= 0}k_{\sigma})+e_K\mathrm{val}_p(\delta_2(\pi_K))\ge 0$ for any lift $(a_{\sigma}e_{\sigma})_{\sigma}\in\oplus_{\sigma,k_{\sigma}\ge 1}Ee_{\sigma}$ of $[\overline{(a_{\sigma}e_{\sigma})_{\sigma}}]\}$,
\item[$\bullet$] ${T''}^{\acute{\mathrm{e}}t}_{\mathrm{cris}}(\delta_1,\delta_2,\{k_{\sigma}\}_{\sigma}):=
\oplus_{\sigma,k_{\sigma}\ge 1}Ee_{\sigma}/\Delta(E)$,
\item[$\bullet$] $T^{\acute{\mathrm{e}}t}_{\mathrm{cris}}(\delta_1,\delta_2,\{k_{\sigma}\}_{\sigma}):={T'}_{\mathrm{cris}}^{\acute{\mathrm{e}}t}(\delta_1,\delta_2,\{k_{\sigma}\}_{\sigma})\sqcup {T''}_{\mathrm{cris}}^{\acute{\mathrm{e}}t}(\delta_1,\delta_2,\{k_{\sigma}\}_{\sigma})$.
\end{itemize}
For any $x\in T^{\acute{\mathrm{e}}t}_{\mathrm{cris}}(\delta_1,\delta_2,\{k_{\sigma}\}_{\sigma})$ we want to define a rank two $E$-filtered ($\varphi,\mathrm{Gal}(K_{n(\delta_1,\delta_2)}/K))$-module $D_{(\delta_1,\delta_2,\{k_{\sigma}\}_{\sigma}),x}$. Here we put $n(\delta_1,\delta_2):=
\mathrm{max}\{n(\delta_1),n(\delta_2)\}$. But we cannot canonically construct these. 
We must fix one more parameter for any $(\delta_1,\delta_2,\{k_{\sigma}\}_{\sigma})\in T_{\mathrm{cris}}$. For this, we consider the following short exact sequence
\begin{equation}
0\rightarrow E^{\times}\xrightarrow[]{g_1} (K_0\otimes_{\mathbb{Q}_p}E)^{\times}\xrightarrow[]{g_2}(K_0\otimes_{\mathbb{Q}_p}E)^{\times}\xrightarrow[]{g_3}
 E^{\times}\rightarrow 0,
\end{equation}
here $g_1:E^{\times}\rightarrow (K_0\otimes_{\mathbb{Q}_p}E)^{\times}$ is the canonical inclusion, $g_2:(K_0\otimes_{\mathbb{Q}_p}E)^{\times}\rightarrow 
(K_0\otimes_{\mathbb{Q}_p}E)^{\times}:x\mapsto \varphi(x)/x$
and $g_3:(K_0\otimes_{\mathbb{Q}_p}E)^{\times}\rightarrow E^{\times}:x\mapsto 
\prod_{i=0}^{f-1}\varphi^i(x)$ where $\varphi:K_0\otimes_{\mathbb{Q}_p}E\rightarrow K_0\otimes_{\mathbb{Q}_p}E:\sum_{i}x_i\otimes y_i\mapsto \sum_{i}\varphi(x_i)\otimes y_i$. We can easily prove that this sequence is exact by using the 
assumption $K_0\subset E$. Then we fix $(\alpha,\beta)\in (K_0\otimes_{\mathbb{Q}_p}E)^{\times}\times (K_0\otimes_{\mathbb{Q}_p}E)^{\times}$ such 
that $g_3(\alpha)=\delta_1(\pi_K)$, $g_3(\beta)=\delta_2(\pi_K)$ when $\delta_1\not= \delta_2$. We fix $\alpha\in (K_0\otimes_{\mathbb{Q}_p}E)^{\times}$ such that $g_3(\alpha)=\delta_1(\pi_K)$ and put $\beta=\alpha$ when $\delta_1=\delta_2$. 
Then we define $D_{(\delta_1,\delta_2,\{k_{\sigma}\}_{\sigma}),x}$ as follows:
\begin{itemize}
\item[$\bullet$]$D_{(\delta_1,\delta_2,\{k_{\sigma}\}_{\sigma}),x}:=(K_{0}\otimes_{\mathbb{Q}_p}E)e_1\oplus
(K_0\otimes_{\mathbb{Q}_p}E)e_2$.
\item[(1)] $N(e_1)=0, N(e_2)=0$.
\item[(2)]When $\delta_1\not= \delta_2$ or when $\delta_1=\delta_2$ 
and $x\in {T'}^{\acute{\mathrm{e}}t}_{\mathrm{cris}}(\delta_1,\delta_2,\{k_{\sigma}\}
_{\sigma})$, we put $\varphi(e_1)=\alpha e_1,  \varphi(e_2)=\beta e_2$.
 (Then $\varphi^f(e_1)=\delta_1(\pi_K)e_1, \varphi^f(e_2)=\delta_2(\pi_K)e_2$.)
\item[(2)'] When $\delta_1=\delta_2$ and $x\in {T''}^{\acute{\mathrm{e}}t}_{\mathrm{cris}}(\delta_1,\delta_2,\{k_{\sigma}\}_{\sigma})$, we put $\varphi(e_1)=\alpha e_1$,  $\varphi(e_2)=\alpha(e_2+e_1)$.
\item[(3)]$g(e_1)=\delta_1(\chi_{\mathrm{LT}}(g))e_1$, $g(e_2)=
\delta_2(\chi_{\mathrm{LT}}(g))e_2$ for any $g\in G_K$.
 \item[(4)] We put $K_{n(\delta_1,\delta_2)}\otimes_{K_0}D_{(\delta_1,\delta_2,\{k_{\sigma}\}_{\sigma}),x}=(K_{n(\delta_1,\delta_2)}\otimes_{\mathbb{Q}_p}E)e_1\oplus (K_{n(\delta_1,\delta_2)}\otimes_{\mathbb{Q}_p}E)e_2\isom\oplus_{\sigma:K\hookrightarrow E}(K_{n(\delta_1,\delta_2)}\otimes_{K,\sigma}E)e_{1,\sigma}\oplus (K_{n(\delta_1,\delta_2)}\otimes_{K,\sigma}E)e_{2,\sigma}=:\oplus_{\sigma}D_{\sigma}$, then 
 \begin{itemize}
 \item[(i)]For $\sigma$ such that $k_{\sigma}\le -1$, we put $\mathrm{Fil}^0D_{\sigma}=D_{\sigma}$
 , $\mathrm{Fil}^1D_{\sigma}=\cdots =\mathrm{Fil}^{-k_{\sigma}}D_{\sigma}=K_{n(\delta_1,\delta_2)}\otimes_{
 K,\sigma}Ee_{1,\sigma}$, $\mathrm{Fil}^{-k_{\sigma}+1}D_{\sigma}=0$.
 \item[(ii)]For $\sigma$ such that $k_{\sigma}=0$, we put $\mathrm{Fil}^0D_{\sigma}=D_{\sigma}$,
 $\mathrm{Fil}^1D_{\sigma}=0$.
 \item[(iii)]For $\sigma$ such that $k_{\sigma}\ge 1$, we put $\mathrm{Fil}^{-k_{\sigma}}D_{\sigma}
 =D_{\sigma}$, $\mathrm{Fil}^{-k_{\sigma}+1}D_{\sigma}=\cdots =\mathrm{Fil}^0D_{\sigma}
 =K_{n(\delta_1,\delta_2)}\otimes_{K,\sigma}E(a_{\sigma}G(\delta_2/\delta_1)_{\sigma}e_{1,\sigma}+e_{2,\sigma})
 $, $\mathrm{Fil}^1D_{\sigma}=0$.
 \end{itemize}
 \end{itemize}
Here when $\delta_1\not= \delta_2$, then $(a_{\sigma}e_{\sigma})_{\sigma}
\in \oplus_{\sigma,k_{\sigma}\ge 1}Ee_{\sigma}$ is a lift of 
$x=[(a_{\sigma}e_{\sigma})_{\sigma}]\in {T}^{\acute{\mathrm{e}}t}_{\mathrm{cris}}(\delta_1,\delta_2,\{k_{\sigma}\}_{\sigma})=\mathbb{P}_E(\oplus_{\sigma,k_{\sigma}\ge 1}Ee_{\sigma})$ and when $\delta_1=\delta_2$, then $(a_{\sigma}e_{\sigma})_{\sigma}\in\oplus_{\sigma,k_{\sigma}\ge 1}Ee_{\sigma}$ is a lift of $[\overline{(a_{\sigma}e_{\sigma})_{\sigma}}]\in {T'}^{\acute{\mathrm{e}}t}_{\mathrm{cris}}(\delta_1,\delta_2,
\{k_{\sigma}\}_{\sigma})=\mathbb{P}_E(\oplus_{\sigma,k_{\sigma}\ge 1}Ee_{\sigma}/\Delta(E))$ 
or is a lift of $\overline{(a_{\sigma}e_{\sigma})_{\sigma}}\in {T''}^{\acute{\mathrm{e}}t}_{\mathrm{cris}}(\delta_1,\delta_2,
\{k_{\sigma}\}_{\sigma})=\oplus_{\sigma,k_{\sigma}\ge 1}Ee_{\sigma}/\Delta(E)$. 
Then we claim that the isomorphism class of $D_{(\delta_1,\delta_2,\{k_{\sigma}\}_{\sigma}),x}$ does not depend on the choice of a lift $(a_{\sigma}e_{\sigma})_{\sigma}$. We prove this claim in the case where $\delta_1=\delta_2$ and $x=\overline{(a_{\sigma}e_{\sigma})_{\sigma}}\in {T''}^{\acute{\mathrm{e}}t}_{\mathrm{cris}}(\delta_1,\delta_2,\{k_{\sigma}\}_{\sigma})$. 
Let us take two lifts $(a_{\sigma}e_{\sigma})_{\sigma}$, $(a'_{\sigma}e_{\sigma})_{\sigma}\in \oplus_{\sigma, k_{\sigma}\in \mathbb{Z}_{\ge 1}}Ee_{\sigma}$ of $x=\overline{(a_{\sigma}e_{\sigma})_{\sigma}}\in {T''}^{\acute{\mathrm{e}}t}_{\mathrm{cris}}(\delta_1,\delta_2,\{k_{\sigma}\}_{\sigma})$. Then there exist $a\in E$ and $b\in E^{\times}$ such that $a'_{\sigma}=ba_{\sigma}+a$ for any $\sigma$ such that $k_{\sigma}\in 
\mathbb{Z}_{\ge 1}$. We denote $D_0:=(K_0\otimes_{\mathbb{Q}_p}E)e_1\oplus
 (K_0\otimes_{\mathbb{Q}_p}E)e_2$ and $D'_0:=(K_0\otimes_{\mathbb{Q}_p}E)e'_1\oplus (K_0\otimes_{\mathbb{Q}_p}E)e'_2$ the filtered $(\varphi, G_K)$-modules 
$D_{(\delta_1,\delta_2,\{k_{\sigma}\}_{\sigma}),x}$ defined by using the lifts $(a_{\sigma}e_{\sigma})_{\sigma}$ and 
$(a'_{\sigma}e_{\sigma})_{\sigma}$ respectively. 
Then it is easy to see that the map 
$D_0\rightarrow D'_0:e_1\mapsto be'_1, e_2\mapsto e_2'+ae_1'$ is 
an isomorphism of filtered $(\varphi,G_K)$-modules.
We can easily prove this claim in other cases.

Next we note the dependence of the choice of 
$(\alpha, \beta)$ as above. We only treat the case $\delta_1\not= \delta_2$.
Let us take another $(\alpha',\beta')\in (K_0\otimes_{\mathbb{Q}_p}E)^{\times}
\times (K_0\otimes_{\mathbb{Q}_p}E)^{\times}$ such that 
$g_3(\alpha')=\delta_1(\pi_K)$, $g_3(\beta')=\delta_2(\pi_K)$. 
Then, by the above exact sequence (2), there exist $\alpha_0,\beta_0\in (K_0\otimes_{\mathbb{Q}_p}E)^{\times}$ such that $\alpha=\alpha' g_2(\alpha_0)$, $\beta=\beta' g_2(\beta_0)$. 
For stressing the dependence of the choice of $(\alpha,\beta)$, we write 
$T^{\acute{\mathrm{e}}t}_{\mathrm{cris}}(\delta_1,\delta_2,\{k_{\sigma}\}_{\sigma},(\alpha,\beta))$, $D_{(\delta_1,\delta_2,\{k_{\sigma}\}_{\sigma},(\alpha,\beta)),x}$, instead of $T^{\acute{\mathrm{e}}t}_{\mathrm{cris}}(\delta_1,\delta_2,\{k_{\sigma}\}_{\sigma})$, $D_{(\delta_1,\delta_2,\{k_{\sigma}\}_{\sigma}),x}$.Then we can easily see that the map 
$D_{(\delta_1,\delta_2,\{k_{\sigma}\}_{\sigma},(\alpha,\beta)),[(a_{\sigma}e_{\sigma})_{\sigma}]}\rightarrow D_{(\delta_1,\delta_2,\{k_{\sigma}\}_{\sigma},}\allowbreak{}_{(\alpha',\beta')), [(a_{\sigma}\frac{\alpha_{0,\sigma}}{\beta_{0,\sigma}}e_{\sigma})_{\sigma}]}:e_1\mapsto \alpha_0 e'_1, 
e_2\mapsto \beta_0 e'_2$ is an isomorphism of filtered $(\varphi, G_K)$-modules
 for any $[(a_{\sigma}e_{\sigma})_{\sigma}]\in T^{\acute{\mathrm{e}}t}_{\mathrm{cris}}(\delta_1,\delta_2,\{k_{\sigma}\}_{\sigma},(\alpha,\beta))$. (Here $e'_i$ is the basis of $D_{(\delta_1,\delta_2,\{k_{\sigma}\}_{\sigma},(\alpha',\beta')),[(a_{\sigma}\frac{\alpha_{0,\sigma}}{\beta_{0,\sigma}}e_{\sigma})_{\sigma}] }$ defined in the same way as in the case $(\alpha,\beta)$ and 
$\alpha_{0,\sigma},\beta_{0,\sigma}\in E^{\times}$ is the $\sigma$-components of $\alpha_0,\beta_0\in (K_0\otimes_{\mathbb{Q}_p}E)^{\times}
\hookrightarrow (K\otimes_{\mathbb{Q}_p}E)^{\times}\isom \prod_{\sigma:K\hookrightarrow E}Ee_{\sigma}^{\times}$.) Other cases are same.
So we have an isomorphism $\iota:T^{\acute{\mathrm{e}}t}_{\mathrm{cris}}(\delta_1,\delta_2,\{k_{\sigma}\}_{\sigma}, (\alpha,\beta))\isom T^{\acute{\mathrm{e}}t}_{\mathrm{cris}}(\delta_1,\delta_2,\{k_{\sigma}\}_{\sigma}, (\alpha',\beta'))$ between parameter spaces such that $D_{(\delta_1,\delta_2,\{k_{\sigma}\}_{\sigma},(\alpha,\beta)), x}\isom D_{(\delta_1,\delta_2,\{k_{\sigma}\}_{\sigma},(\alpha',\beta')), \iota(x)}$ for any $x\in T^{\acute{\mathrm{e}}t}_{\mathrm{cris}}(\delta_1,\delta_2,\{k_{\sigma}\}_{\sigma}, (\alpha,\beta))$.

Henceforth, we fix a $(\alpha,\beta)$ for any $(\delta_1,\delta_2,\{k_{\sigma}\}_{\sigma})\in T_{\mathrm{cris}}$. 
In the proof of the next theorem, we will prove that these are weakly admissible filtered $(\varphi, G_K)$-modules. So, by the ``weakly admissible imply admissible" theorem $\cite[\mathrm{Theorem} B]{Be04}$, for any $x\in T^{\acute{\mathrm{e}}t}_{\mathrm{cris}}(\delta_1,\delta_2,\{k_{\sigma}\}_{\sigma})$, there exists unique (up to isomorphism) two dimensional  potentially cristalline $E$-representation $ V_{(\delta_1,\delta_2,\{k_{\sigma}\}_{\sigma}),x } $ such that $D^{K_{n(\delta_1,\delta_2)}}_{\mathrm{cris}}(V_{(\delta_1,\delta_2,\{k_{\sigma}\}_{\sigma}),x})\isom D_{(\delta_1,\delta_2,\{k_{\sigma}\}_{\sigma}),x}$.
Then our main result of classification concerning  potentially cristalline split trianguline $E$-representations 
is as follows.
\begin{thm}$\label{42}$
Let $V$ be a two dimensional $E$-representation.
Then the following conditions are equivalent:
\begin{itemize}
\item[$\mathrm{(1)}$]$V$ is split trianguline and potentially cristalline.
\item[$\mathrm{(2)}$]There exist $(\delta_1,\delta_2,\{k_{\sigma}\}_{\sigma})\in T_{\mathrm{cris}}$ and $x\in T^{\acute{\mathrm{e}}t}_{\mathrm{cris}}(\delta_1,\delta_2,\{k_{\sigma}\}_{\sigma})$ and $\{w_{\sigma}\}_{\sigma}$ such that $w_{\sigma}\in\mathbb{Z}$ for any $\sigma$, such that $V\isom V_{(\delta_1,\delta_2,\{k_{\sigma}\}_{\sigma}),
x}(\prod_{\sigma}\sigma(\chi_{\mathrm{LT}}^{w_{\sigma}}))$.
\end{itemize}

\end{thm}
\begin{proof}
First we prove that $(1)$ implies $(2)$.
Let us assume that $V$ is a split trianguline and potentially cristalline $E$-representation.
Then, by Proposition $\ref{41}$, there exists a $(\delta_1,\delta_2)\in S^+$ and an $s\in S^{\acute{\mathrm{e}}t}_{\mathrm{cris}}(\delta_1,\delta_2)$ such that 
$V\isom V(s)$.  Twisting $V$ by a suitable $\prod_{\sigma}\sigma(\chi_{\mathrm{LT}}(x)^{l_{\sigma}})$, we may assume that $W(\delta_2)$ has generalized Hodge-Tate weight $\{0\}_{\sigma}$. If we write the generalized Hodge-Tate weight of $W(\delta_1)$ as $\{k_{\sigma}\}_{\sigma}$, then 
we have $\delta_1=\tilde{\delta}_1\prod_{\sigma}\sigma(x)^{k_{\sigma}}$ and 
$\delta_2=\tilde{\delta}_2$ for some locally constant characters $\tilde{\delta_1}, \tilde{\delta_2}:
K^{\times}\rightarrow E^{\times}$. Then the condition 
that $(\delta_1,\delta_2)\in S^+$ is equivalent to the condition that $(\tilde{\delta}_1,\tilde{\delta}_2,\{k_{\sigma}\}_{\sigma})\in T_{\mathrm{cris}}$.
Now we explicitly calculate $D^{K_{n(\tilde{\delta}_1, \tilde{\delta}_2)}}_{\mathrm{cris}}(V)$.
First we  compute this in the case where $\tilde{\delta}_1\not= \tilde{\delta}_2$ or $\tilde{\delta}_1=\tilde{\delta}_2$ and $s\in {S'}^{\acute{\mathrm{e}}t}_{\mathrm{cris}}(\delta_1, \delta_2)$. 
Then, by definition of ${S'}^{\acute{\mathrm{e}}t}_{\mathrm{cris}}(\delta_1, \delta_2)\isom
\mathbb{P}_E(D^K_{\mathrm{dR}}(W(\delta_1/\delta_2))/(D^K_{\mathrm{cris}}(W(\delta_1/\delta_2))^{\varphi=1}+\mathrm{Fil}^0D^K_{\mathrm{dR}}(W(\delta_1/\delta_2)))$ and by definition of the boundary map $D^K_{\mathrm{dR}}(W(\delta_1/\delta_2))/(D^K_{\mathrm{cris}}(W(\delta_1/\delta_2))^{\varphi=1}+\mathrm{Fil}^0D^K_{\mathrm{dR}}(W(\delta_1/\delta_2))\hookrightarrow \mathrm{H}^1(G_K, W(\delta_1/
\delta_2))$, $W(s)$ is given as follows.
\begin{itemize}
\item[(1)]$W_e(s):=W_e(\delta_1)\oplus W_e(\delta_2)$ such that $g(x,y):=(gx,gy)$ for any 
$g\in G_K, x\in W_e(\delta_1), y\in W_e(\delta_2)$.
\item[(2)]$W^+_{\mathrm{dR}}(s):=W^+_{\mathrm{dR}}(\delta_1)\oplus W^+_{\mathrm{dR}}(\delta_2)$ such that $g(x,y):=(gx,gy)$ for any 
$g\in G_K, x\in W^+_{\mathrm{dR}}(\delta_1), y\in W^+_{\mathrm{dR}}(\delta_2)$.
\item[(3)]$\iota:B_{\mathrm{dR}}\otimes_{B_e}W_e(s)\isom B_{\mathrm{dR}}\otimes_{B^+_{\mathrm{dR}}}W^+_{\mathrm{dR}}(s)$ is given by $\iota(x,y):=(x+a\otimes y, y)$ for any 
$ x\in B_{\mathrm{dR}}\otimes_{B_e}W_e(\delta_1), y\in B_{\mathrm{dR}}\otimes_{B_e}W_e(\delta_2)$. Here, $a\in D^K_{\mathrm{dR}}(\delta_1/\delta_2)$ is a lifing of $\bar{a}\in D^K_{\mathrm{dR}}(W(\delta_1/\delta_2))/(D^K_{\mathrm{cris}}(W(\delta_1/\delta_2))^{\varphi=1}+\mathrm{Fil}^0D^K_{\mathrm{dR}}(W(\delta_1/\delta_2))$ corresponding to $s$.
\end{itemize}
By definition of $W(\delta)$ for any $\delta:K^{\times}\rightarrow E^{\times}$,  we have $D^{K_{n(\tilde{\delta}_1,\tilde{\delta}_2)}}_{\mathrm{cris}}(W(\delta_1))\isom (K_0\otimes_{\mathbb{Q}_p}E)e_1$
, $D^{K_{n(\tilde{\delta}_1,\tilde{\delta}_2)}}_{\mathrm{cris}}(W(\delta_2))\isom (K_0\otimes_{\mathbb{Q}_p}E)e_2$ such that 
\begin{itemize}
\item[$\bullet$]$\varphi^f(e_1)=\tilde{\delta_1}(\pi_K)e_1, \varphi^f(e_2)=\tilde{\delta_2}(\pi_K)e_2$, 
\item[$\bullet$]$g(e_1)=\tilde{\delta_1}(\chi_{\mathrm{LT}}(g))e_1, g(e_2)=\tilde{\delta_2}(\chi_{\mathrm{LT}}(g))e_2$ for any $g\in G_K$.
\end{itemize}
Then we can take a basis $e_1$, $e_2$ such that $\varphi(e_1)=\alpha e_1$, $\varphi(e_2)=\beta e_2$, here $(\alpha,\beta)$ is the fixed pair as in before Theorem $\ref{42}$ for $(\tilde{\delta}_1,\tilde{\delta}_2,\{k_{\sigma}\}_{\sigma})\in T_{\mathrm{cris}}$
By $(1)$ in the above description of $W(s)$, we have $D^{K_{n(\tilde{\delta_1}, \tilde{\delta_2})}}_{\mathrm{cris}}(V)=D^{K_{n(\tilde{\delta}_1,\tilde{\delta}_2)}}_{\mathrm{cris}}(W(\delta_1))\oplus D^{K_{n(\tilde{\delta}_1,\tilde{\delta}_2)}}_{\mathrm{cris}}(W(\delta_2))$ as $(\varphi,G_K)$-module.
 We compute the filtration on $D^{K_{n(\tilde{\delta}_1,\tilde{\delta}_2)}}_{\mathrm{dR}}(W(s))\isom 
 K_{n(\tilde{\delta}_1,\tilde{\delta}_2)}\otimes_{K_0}D^{K_{n(\tilde{\delta}_1,\tilde{\delta}_2)}}_{\mathrm{cris}}(W(s))
 \isom (K_{n(\tilde{\delta}_1,\tilde{\delta}_2)}\otimes_{\mathbb{Q}_p}E)e_1\oplus  (K_{n(\tilde{\delta}_1,\tilde{\delta}_2)}\otimes_{\mathbb{Q}_p}E)e_2\isom \oplus_{\sigma:K\hookrightarrow E} (K_{n(\tilde{\delta}_1,\tilde{\delta}_2)}\otimes_{K,\sigma}E)e_{1,\sigma}\oplus  (K_{n(\tilde{\delta}_1,\tilde{\delta}_2)}\otimes_{K,\sigma}E)e_{2,\sigma}=:\oplus_{\sigma}
 D_{\sigma}$ (Here we put $D_{\sigma}:=(K_{n(\tilde{\delta}_1,\tilde{\delta}_2)}\otimes_{K,\sigma}E)e_{1,\sigma}\oplus  (K_{n(\tilde{\delta}_1,\tilde{\delta}_2)}\otimes_{K,\sigma}E)e_{2,\sigma}$). 
 Because $D^K_{\mathrm{dR}}(W(\delta_1/\delta_2))=(W_{\mathrm{dR}}(\delta_1/\delta_2))^{G_K}
 \isom \oplus_{\sigma} (W_{\mathrm{dR}}(\delta_1/\delta_2)\otimes_{B_{\mathrm{dR}}\otimes_{\mathbb{Q}_p}E}B_{\mathrm{dR}}\otimes_{K,\sigma}E)^{G_K}$ and $\mathrm{Fil}^0D^K_{\mathrm{dR}}(W(\delta_1/\delta_2))=\oplus_{\sigma}(W^+_{\mathrm{dR}}(\delta_1/\delta_2)\otimes_{B^+_{\mathrm{dR}}\otimes_{\mathbb{Q}_p}E}B^+_{\mathrm{dR}}\otimes_{K,\sigma}E)^{G_K}
 \isom \oplus_{\sigma, k_{\sigma}\le 0}(W_{\mathrm{dR}}(\delta_1/\delta_2)\otimes_{B_{\mathrm{dR}}\otimes_{\mathbb{Q}_p}E}B_{\mathrm{dR}}\otimes_{K,\sigma}E)^{G_K}$ , so 
 we can take $a\in D^K_{\mathrm{dR}}(W(\delta_1/\delta_2))$ in (3) as $a=(a_{\sigma})\in\oplus_{\sigma}(W_{\mathrm{dR}}(\delta_1/\delta_2)\otimes_{B_{\mathrm{dR}}\otimes_{\mathbb{Q}_p}E}B_{\mathrm{dR}}\otimes_{K,\sigma}E)^{G_K}$ such that $a_{\sigma}=0$ for any $\sigma$ such that 
 $k_{\sigma}\le 0$. Then, for any $\sigma$ such that $k_{\sigma}\le 0$ or $k_{\sigma}\ge 1$ and $a_{\sigma}=0$, we have 
\begin{itemize}
\item[$\bullet$]$\mathrm{Fil}^iD_{\sigma}= \mathrm{Fil}^iD^{K_{n(\tilde{\delta}_1,\tilde{\delta}_2)}}_{\mathrm{dR}}(W(\delta_1))_{\sigma}\oplus 
 \mathrm{Fil}^iD^{K_{n(\tilde{\delta}_1,\tilde{\delta}_2)}}_{\mathrm{dR}}(W(\delta_2))_{\sigma}$.
\end{itemize}
 ($D^{K_{n(\tilde{\delta}_1,\tilde{\delta}_2)}}_{\mathrm{dR}}(W(\delta_i))\isom \oplus _{\sigma}D^{K_{n(\tilde{\delta}_1,\tilde{\delta}_2)}}_{\mathrm{dR}}(W(\delta_i))\otimes_{K_{n(\tilde{\delta}_1,\tilde{\delta}_2)}\otimes_{\mathbb{Q}_p}E}K_{n(\tilde{\delta}_1,\tilde{\delta}_2)}\otimes_{K,\sigma}E=:\oplus_{\sigma}D^{K_{n(\tilde{\delta}_1,\tilde{\delta}_2)}}_{\mathrm{dR}}(W(\delta_1))_{\sigma}$ for $i\in\{1,2\}$.) Because we have
 \begin{itemize}
\item[$\bullet$]$\mathrm{Fil}^{-k_{\sigma}}D^{K_{n(\tilde{\delta}_1,\tilde{\delta}_2)}}_{\mathrm{dR}}(W(\delta_1))_{\sigma}=D^{K_{n(\tilde{\delta}_1,\tilde{\delta}_2)}}_{\mathrm{dR}}(W(\delta_1))_{\sigma}$, $\mathrm{Fil}^{-k_{\sigma}+1}D^{K_{n(\tilde{\delta}_1,\tilde{\delta}_2)}}_{\mathrm{dR}}(W(\delta_1))_{\sigma}=0$, 
\item[$\bullet$]$\mathrm{Fil}^0D^{K_{n(\tilde{\delta}_1,\tilde{\delta}_2)}}_{\mathrm{dR}}(W(\delta_2))_{\sigma}=D^{K_{n(\tilde{\delta}_1,\tilde{\delta}_2)}}_{\mathrm{dR}}(W(\delta_2))_{\sigma}$, $\mathrm{Fil}^1D^{K_{n(\tilde{\delta}_1,\tilde{\delta}_2)}}_{\mathrm{dR}}(W(\delta_2))_{\sigma}=0$, 
\end{itemize}
we get the filtration on $D_{\sigma}$ as in (i), (ii), (iii) of (4) given before Theorem $\ref{42}$ for such $\sigma$. 
 Finally we calculate the filtration on $D_{\sigma}$ for $\sigma$ such that 
 $k_{\sigma}\ge 1$ and $a_{\sigma}\not= 0$. 
 Then we have $\mathrm{Fil}^{-k_{\sigma}}D_{\sigma}=D_{\sigma}$ and $\mathrm{Fil}^1D_{\sigma}
 =0$. For any $xe_{1,\sigma}+ye_{2,\sigma}\in D_{\sigma}$ ($x,y\in  K_{n(\tilde{\delta}_1,\tilde{\delta}_2)}\otimes_{K,\sigma}E$), we have $xe_{1,\sigma}+ye_{2,\sigma}\in \mathrm{Fil}^0D_{\sigma}$ if and only if 
 $xe_{1,\sigma}+y a_{\sigma}\otimes e_{2,\sigma}\in \mathrm{Fil}^0D^{K_{n(\tilde{\delta}_1,\tilde{\delta}_2)}}_{\mathrm{dR}}(W(\delta_1))_{\sigma}= 0$ and $ye_{2,\sigma}\in \mathrm{Fil}^0D_{\mathrm{dR}}(W(\delta_2))_{\sigma}=
 D^{K_{n(\tilde{\delta}_1,\tilde{\delta}_2)}}_{\mathrm{dR}}(W(\delta_2))_{\sigma}$. 
 Because $a_{\sigma}\in D^K_{\mathrm{dR}}(W(\delta_1/\delta_2))_{\sigma}\isom E$ is non-zero, it is a base of $D^K_{\mathrm{dR}}(W(\delta_1/\delta_2))_{\sigma}$ as an $E$-vector space. So it is a base of $D^{K_{n(\tilde{\delta}_1,\tilde{\delta}_2)}}_{\mathrm{dR}}(W(\delta_1/\delta_2))_{\sigma}$ as a $K_{n(\tilde{\delta}_1,\tilde{\delta}_2)}\otimes_{K,\sigma}E$-module. So $a_{\sigma}\otimes e_{2,\sigma}$ is a base of $D^{K_{n(\tilde{\delta}_1,\tilde{\delta}_2)}}_{\mathrm{dR}}(W(\delta_1))$. So
 there exists unique $z\in (K_{n(\tilde{\delta}_1,\tilde{\delta}_2)}\otimes_{K,\sigma}E)^{\times}$ such that 
 $a_{\sigma}\otimes e_{2,\sigma}=ze_{1,\sigma}$. Then $xe_{1,\sigma}+ya_{\sigma}\otimes e_{2,\sigma}=(x+zy)e_{1,\sigma}$, so $\mathrm{Fil}^{-k_{\sigma}+1}D_{\sigma}=\mathrm{Fil}^0D_{\sigma}=K_{n(\tilde{\delta}_1,\tilde{\delta}_2)}\otimes_{K,\sigma}E(-ze_{1,\sigma}+e_{2,\sigma})$. 
 By comparing $G_K$-actions on $a_{\sigma}\otimes e_{2,\sigma}=ze_{1,\sigma}$, we can see that 
 $g(z)=\frac{\tilde{\delta}_2(\chi_{\mathrm{LT}}(g))}{\tilde{\delta}_1(\chi_{\mathrm{LT}}(g))}
 z$ for any $g\in G_K$. So there exists unique $b_{\sigma}\in E^{\times}$ such that 
 $z=-b_{\sigma}G(\tilde{\delta_2}/\tilde{\delta_1})_{\sigma}$. 
 By combining all these calculations, we get $D^{K_{n(\tilde{\delta}_1,\tilde{\delta}_2)}}_{\mathrm{cris}}(V(s))\isom D_{(\tilde{\delta}_1,\tilde{\delta}_2,\{k_{\sigma}\}_{\sigma}), x}$ (here $x=[(b_{\sigma}e_{\sigma})_{\sigma}]\in T^{\acute{\mathrm{e}}t}_{\mathrm{cris}}(\tilde{\delta}_1,\tilde{\delta}_2,\{k_{\sigma}\}_{\sigma})$ or $x=[(\overline{b_{\sigma}e_{\sigma}})_{\sigma}]\in {T'}^{\acute{\mathrm{e}}t}_{\mathrm{cris}}(\tilde{\delta}_1,\tilde{\delta}_2,\{k_{\sigma}\}_{\sigma}))$.
So $V(s)\isom V_{(\tilde{\delta_1},\tilde{\delta_2}, \{k_{\sigma}\}_{\sigma}),x}$. 
 (Here we put $b_{\sigma}=0$ if $a_{\sigma}=0$.)
From this, we know that $D_{(\tilde{\delta}_1,\tilde{\delta}_2,\{k_{\sigma}\}_{\sigma}),x}$ is weakly admissible for such $x$.
 So we can prove that $(1)$ implies $(2)$ in the case $\tilde{\delta}_1\not= \tilde{\delta}_2$ or 
 $\tilde{\delta}_1=\tilde{\delta}_2$ and $s\in {S'}^{\acute{\mathrm{e}}t}(\delta_1,\delta_2)$. 

 Next we compute in the case where $\tilde{\delta_1}=\tilde{\delta_2}$ and $s\in {S''}(\delta_1,\delta_2)$. Then, by the construction of $s(\{k_{\sigma}\}_{\sigma})$ and by definition of $S''(\delta_1,\delta_2)$ before Proposition $\ref{41}$, $W(s):=(W_e(s), W^+_{\mathrm{dR}}(s),\iota)$ is given as follows:
 \begin{itemize}
 \item[(1)]$W_e(s):=W_e(\delta_1)\oplus W_e(\delta_2)$ such that $g(x,y):=
 (gx+f\mathrm{deg}(g)gy,gy)$ for any $g\in G_K$, $x\in W_e(\delta_1)$, $y\in W_e(\delta_2)$. (Here we use the fact that $\delta_1/\delta_2=\prod_{\sigma}\sigma(x)^{k_{\sigma}}$ implies that $W_e(\delta_1)=W_e(\delta_2)$ (Proposition $\ref{28}$)).
 \item[(2)]$W^+_{\mathrm{dR}}(s):=W^+_{\mathrm{dR}}(\delta_1)\oplus W^+_{\mathrm{dR}}(\delta_2)$ such that $g(x,y):=(gx,gy)$ for any $g\in G_K$, $x\in W^+_{\mathrm{dR}}(\delta_1)$, $y\in W^+_{\mathrm{dR}}(\delta_2)$.
 \item[(3)]$\iota:B_{\mathrm{dR}}\otimes_{B_e}W_e(s)\isom B_{\mathrm{dR}}\otimes_{B^+_{\mathrm{dR}}}W^+_{\mathrm{dR}}(s)$ is given by $\iota(x,y):=
 (x+(a_{\mathrm{cris}}+a)y,y)$ for any $x\in B_{\mathrm{dR}}\otimes_{B_e}W_e(\delta_1)$, $y\in B_{\mathrm{dR}}\otimes_{B_e}W_e(\delta_2)$. Here, $a\in D^K_{\mathrm{dR}}(W(\delta_1/\delta_2))\isom (B_{\mathrm{dR}}\otimes_{\mathbb{Q}_p}E)^{G_K}=K\otimes_{\mathbb{Q}_p}E$ is a lifting of $\bar{a}\in D^K_{\mathrm{dR}}(W(\delta_1/\delta_2))/(D^K_{\mathrm{cris}}(W(\delta_1/\delta_2))+\mathrm{Fil}^0D^K_{\mathrm{dR}}(W(\delta_1/\delta_2)))$ coresponding to $s=s(\{k_{\sigma}\}_{\sigma})+\bar{a}$.
 \end{itemize}
Then we have $D^{K_{n(\tilde{\delta}_1)}}_{\mathrm{cris}}(W(\delta_1))=(B_{\mathrm{cris}}\otimes_{B_e}W_e(\delta_1))^{G_{K_{n(\tilde{\delta}_1)}}}\isom (B_{\mathrm{cris}}\otimes_{B_e}W_e(\delta_2))^{G_{K_{n(\tilde{\delta}_1)}}}=D^{K_{n(\tilde{\delta}_1)}}_{\mathrm{cris}}(W(\delta_2))$ because $W_e(\delta_1)\isom W_e(\delta_2)$.
So we have following isomorphisms of $(\varphi,G_K)$-modules $
D^{K_{n(\tilde{\delta}_1)}}_{\mathrm{cris}}(W(\delta_1))\isom (K_0\otimes_{\mathbb{Q}_p}E)e_1$, $D^{K_{n(\tilde{\delta}_2)}}_{\mathrm{cris}}(W(\delta_2))\isom (K_0\otimes_{\mathbb{Q}_p}E)e'_2$  such that 
\begin{itemize}
\item[(1)]$\varphi(e_1)=\alpha e_1$, $\varphi(e'_2)=\alpha e'_2$,
\item[(2)]$g(e_1)=\tilde{\delta_1}(\chi_{\mathrm{LT}}(g))e_1$, 
$g(e'_2)=\tilde{\delta_1}(\chi_{\mathrm{LT}}(g))e'_2$ for any $g\in G_K$.
\end{itemize}
(Here $\alpha\in (K_0\otimes E)^{\times}$ is the fixed one such that 
$g_3(\alpha)=\tilde{\delta}_1(\pi_K)$ defined before Theorem 
$\ref{42}$.)
Then, by (1) in the definition of $W(s)$, we have $D^{K_{n(\tilde{\delta}_1)}}_{\mathrm{cris}}(W(s))=(K_0\otimes_{\mathbb{Q}_p}E)e_1\oplus (K_0\otimes_{\mathbb{Q}_p}E)(-e'_2+a_{\mathrm{cris}}e_1)$. So if we put $e_2:=-e'_2+a_{\mathrm{cris}}e_1$, we have $D^{K_{n(\tilde{\delta}_1)}}_{\mathrm{cris}}(W(s))=(K_0\otimes_{\mathbb{Q}_p}E)e_1\oplus (K_0\otimes_{\mathbb{Q}_p}E)e_2$ such that
\begin{itemize}
\item[(1)]$\varphi(e_1)=\alpha e_1$, $\varphi(e_2)=\alpha(e_2+e_1)$,
\item[(2)]$g(e_1)=\tilde{\delta_1}(\chi_{\mathrm{LT}}(g))e_1$, 
$g(e_2)=\tilde{\delta_1}(\chi_{\mathrm{LT}}(g))e_2$ for any $g\in G_K$.
\end{itemize}
We compute the filtration on $K_{n(\tilde{\delta}_1)}\otimes_{K_0}D^{K_{n(\tilde{\delta}_1)}}_{\mathrm{cris}}(W(s))=(K_{n(\tilde{\delta}_1)}\otimes_{\mathbb{Q}_p}E)e_1\oplus (K_{n(\tilde{\delta}_1)}\otimes_{\mathbb{Q}_p}E)e_2$ as follows.
We can take $a:=(a_{\sigma})\in D^K_{\mathrm{dR}}(\delta_1/\delta_2)=\oplus_{\sigma}(D^K_{\mathrm{dR}}\otimes_{K\otimes_{\mathbb{Q}_p}E}B_{\mathrm{dR}}\otimes_{K,\sigma}E)$ such that $a_{\sigma}=0$ for any $\sigma$ such that $k_{\sigma}\in
\mathbb{Z}_{\le 0}$. Then we can compute the filtration as in the previous case. Then we have $D^{K_{n(\tilde{\delta}_1)}}_{\mathrm{cris}}(W(s))\isom D_{(\tilde{\delta}_1,\tilde{\delta}_2,\{k_{\sigma}\}_{\sigma}),x}$ for some $x=\overline{(b_{\sigma}e_{\sigma})_{\sigma}}\in {T''}^{\acute{\mathrm{e}}t}_{\mathrm{cris}}(\tilde{\delta}_2,\{k_{\sigma}\}_{\sigma})\isom \oplus_{\sigma, k_{\sigma}\ge 1}Ee_{\sigma}/\Delta(E)$ such that a lift $(b_{\sigma}e_{\sigma})_{\sigma}$ of $\overline{(b_{\sigma}e_{\sigma})_{\sigma}}$ satisfies 
$\{\sigma|b_{\sigma}\not= 0\}=\{\sigma|a_{\sigma}\not= 0\}$. 
So, for this $x$, $D_{(\tilde{\delta}_1,\tilde{\delta}_2,\{k_{\sigma}\}_{\sigma}),x}$ is weakly admissible and we have $V(s)\isom V_{(\tilde{\delta}_1,\tilde{\delta}_2,\{k_{\sigma}\}_{\sigma}),x}$.
So we have proved that (1) implies (2) in all cases.

Next we prove that (2) implies (1).
Let us take $(\delta_1,\delta_2,\{k_{\sigma}\}_{\sigma})\in T_{\mathrm{cris}}$.
We prove this in the case $\delta_1\not=\delta_2$. (We can prove the other cases in the same way.) Let us take $[(b_{\sigma}e_{\sigma})_{\sigma}]\in T^{\acute{\mathrm{e}}t}_{\mathrm{cris}}(\delta_1,\delta_2,\{k_{\sigma}\}_{\sigma})\subset \mathbb{P}_E(\oplus_{\sigma,k_{\sigma}\ge 1}Ee_{\sigma})$. 
If we put $\tilde{\delta}_1=\delta_1\prod_{\sigma}\sigma(x)^{k_{\sigma}}$, $\tilde{\delta}_2=\delta_2$, we can see that $(\tilde{\delta}_1, \tilde{\delta}_2)\in S^+$. From the argument of the proof of the claim that (1) implies (2), we can see that there exists an element $s=[(a_{\sigma}e_{\sigma})_{\sigma}]\in S'(\tilde{\delta}_1,\tilde{\delta}_2)$ such that $D^{K_{n(\delta_1,\delta_2)}}_{\mathrm{cris}}(W(s))\isom D_{(\delta_1,\delta_2,\{k_{\sigma}\}_{\sigma}),[(b_{\sigma}e_{\sigma})_{\sigma}]}$. Moreover, from the previous argument, we can also see that $\{\sigma|a_{\sigma}\not= 0\}=\{\sigma|b_{\sigma}\not= 0\}$. This implies that $s$ is contained in ${S'}^{\acute{\mathrm{e}}t}(\tilde{\delta}_1,\tilde{\delta}_2)$, i.e. $W(s)=W(V(s))$ for some split trianguline $E$-representation. So we have $D_{(\delta_1,\delta_2,\{k_{\sigma}\}_{\sigma}),[(b_{\sigma}e_{\sigma})_{\sigma}]}$ is weakly admissible and $V(s)\isom 
V_{(\delta_1,\delta_2,\{k_{\sigma}\}_{\sigma}),[(b_{\sigma}e_{\sigma})_{\sigma}]}$.
We have finished the proof of the theorem.
 \end{proof}

The following is the corollary of Theorem $\ref{36}$.
\begin{corollary}$\label{43}$
\begin{itemize}
\item[$\mathrm{(1)}$]
Let $(\delta_1,\delta_2,\{k_{\sigma}\}_{\sigma})\in T_{\mathrm{cris}}$ such that $\delta_1\not=\delta_2$. Let $[(a_{\sigma}e_{\sigma})_{\sigma}]\in T^{\acute{\mathrm{e}}t}_{\mathrm{cris}}(\delta_1,\delta_2,\{k_{\sigma}\}_{\sigma})$. Then there exists unique $((\delta'_1,\delta'_2,\{k'_{\sigma}\}_{\sigma}), x, \{w_{\sigma}\}_{\sigma})$ where $(\delta'_1,\delta'_2,\{k'_{\sigma}\}_{\sigma})\in T_{\mathrm{cris}}$, $x\in T^{\acute{\mathrm{e}}t}_{\mathrm{cris}}(\delta'_1,\delta'_2,\allowbreak\{k'_{\sigma}\}_{\sigma})$, $w_{\sigma}\in \mathbb{Z}$, satisfying 
$V_{(\delta_1,\delta_2,\{k_{\sigma}\}_{\sigma}),[(a_{\sigma}e_{\sigma})_{\sigma}]}\isom 
V_{(\delta'_1,\delta_2',\{k'_{\sigma}\}_{\sigma}),x}(\prod_{\sigma}\sigma(\chi_{\mathrm{LT}})^{w_{\sigma}})$ and $((\delta'_1,\delta'_2,\allowbreak\{k'_{\sigma}\}_{\sigma}), x, \{w_{\sigma}\}_{\sigma})\not= ((\delta_1,\delta_2,\{k_{\sigma}\}_{\sigma}), [(a_{\sigma}e_{\sigma})_{\sigma}], \{0\}_{\sigma})$.
Such unique $((\delta'_1,\delta'_2,\{k'_{\sigma}\}_{\sigma}), \allowbreak 
x, \{w_{\sigma}\}_{\sigma})$ satisfies
 \begin{itemize}
 \item[$\mathrm{(i)}$]$\delta_1'|_{\mathcal{O}_K^{\times}}=\delta_2|_{\mathcal{O}_K^{\times}}, \delta_2'|_{\mathcal{O}_K^{\times}}=\delta_1|_{\mathcal{O}_K^{\times}}$,
 \item[$\mathrm{(ii)}$]$\delta_1'(\pi_K)=\delta_2(\pi_K)\prod_{\sigma, a_{\sigma}=0 \,\mathrm{or}\, k_{\sigma}\le 0}\sigma(\pi_K)^{k_{\sigma}},\delta_2'(\pi_K)=\delta_1(\pi_K)
 \prod_{\sigma,a_{\sigma}=0\,\mathrm{or}\, k_{\sigma}\le 0}\sigma(\pi_K)^{k_{\sigma}}$,
 \item[$\mathrm{(iii)}$]$k'_{\sigma}=k_{\sigma}$ if $a_{\sigma}\not= 0$, $k'_{\sigma}=-k_{\sigma}$ for other $\sigma$,
 \item[$\mathrm{(iv)}$]$w_{\sigma}=0$ if $a_{\sigma}\not= 0$, $w_{\sigma}=k_{\sigma}$ for other $\sigma$,
 \item[$\mathrm{(v)}$]unique suitable $x=[(a'_{\sigma}e_{\sigma})_{\sigma}]\in T^{\acute{\mathrm{e}}t}_{\mathrm{cris}}(\delta'_1,\delta'_2,\{k'_{\sigma}\}_{\sigma})$ such that $\{\sigma|a'_{\sigma}\not= 0\}=\{\sigma|a_{\sigma}\not= 0\}$.
 \end{itemize}
\item[$\mathrm{(2)}$]
Let $(\delta_1,\delta_2,\{k_{\sigma}\}_{\sigma})$, $(\delta'_1,\delta'_2,\{k'_{\sigma}\}_{\sigma})\in T_{\mathrm{cris}}$ such that 
$\delta_1=\delta_2$. Let $x\in {T''}^{\acute{\mathrm{e}}t}_{\mathrm{cris}}(\delta_1,\delta_2,\{k_{\sigma}\}_{\sigma})$, $y\in T^{\acute{\mathrm{e}}t}_{\mathrm{cris}}(\delta'_1,\delta'_2,\{k'_{\sigma}\}_{\sigma})$, $\{w_{\sigma}\}_{\sigma}$ such that $w_{\sigma}\in\mathbb{Z}$ for any $\sigma$.  Then we have an isomorphism 
$V_{(\delta_1,\delta_2,\{k_{\sigma}\}_{\sigma}),x}\isom V_{(\delta'_1,\delta'_2,\{k'_{\sigma}\}_{\sigma}),y}(\prod_{\sigma}\sigma(\chi_{\mathrm{LT}})^{w_{\sigma}})$ if and only if $((\delta'_1,\delta'_2,\{k'_{\sigma}\}_{\sigma}),y)\allowbreak=((\delta_1,\delta_2,\{k_{\sigma}\}_{\sigma}),x)$ and $w_{\sigma}=0$ for any $\sigma$.
\end{itemize}
 
 \end{corollary}
 \begin{proof}
(1) is the corollary of (2) of Theorem $\ref{36}$.
 Because $\delta_1\not=\delta_2$, then by the proof of the above theorem, 
 we have $V_{(\delta_1,\delta_2,\{k_{\sigma}\}_{\sigma}),[(a_{\sigma}e_{\sigma})_{\sigma}]}\isom V(s)$ where $s=[(b_{\sigma}e_{\sigma})_{\sigma}]\in {S'}^{\acute{\mathrm{e}}t}_{\mathrm{cris}}(\delta_1\prod_{\sigma}\sigma(x)^{k_{\sigma}},\delta_2)$ for some $[(b_{\sigma}e_{\sigma})_{\sigma}]$ such that $\{\sigma|b_{\sigma}\not= 0\}=\{\sigma|a_{\sigma}\not= 0\}$. Then, by Theorem $\ref{36}$, we have $V(s)\isom V(s')$ for unique $s'\not=s$ such that 
 $s'=[(b'_{\sigma}e_{\sigma})_{\sigma}]\in {S'}^{\acute{\mathrm{e}}t}_{\mathrm{cris}}(\prod_{\sigma,b_{\sigma}
 \not= 0}\sigma(x)^{k_{\sigma}}\delta_2, \prod_{\sigma,b_{\sigma}\not= 0}
 \sigma(x)^{-k_{\sigma}}\allowbreak\prod_{\sigma}\sigma(x)^{k_{\sigma}}\delta_1)={S'}^{\acute{\mathrm{e}}t}_{\mathrm{cris}}(\prod_{\sigma,b_{\sigma}
 \not= 0}\sigma(x)^{k_{\sigma}}\delta_2,  \allowbreak\prod_{\sigma,k_{\sigma}\le 0\,\mathrm{ or} \,b_{\sigma}= 0}\sigma(x)^{k_{\sigma}}\delta_1)$ such that $\{\sigma|b'_{\sigma}\not= 0\}=\{\sigma|b_{\sigma}\not= 0\}$. By the proof of the above theorem, we can see that $s'$ corresponds to $V_{(\delta_1',\delta_2',\{k'_{\sigma}\}_{\sigma}),[(a'_{\sigma}e_{\sigma})_{\sigma}]}\allowbreak(\prod_{\sigma,a_{\sigma}=0\,\mathrm{or}\, k_{\sigma}\le 0}\sigma(\chi_{\mathrm{LT}})^{k_{\sigma}})$ such that 
 $\delta_1'|_{\mathcal{O}_K^{\times}}=\delta_2|_{\mathcal{O}_K^{\times}}, 
 \delta_2'|_{\mathcal{O}_K^{\times}}=\delta_1|_{\mathcal{O}_K^{\times}}$, 
 $\delta_1'(\pi_K)= \prod_{\sigma,a_{\sigma}=0 \,\mathrm{or} \,k_{\sigma}\le 0}\allowbreak\sigma(\pi_K)^{k_{\sigma}}\delta_2(\pi_K), \delta_2'(\pi_K)=\prod_{\sigma,a_{\sigma}=0\,\mathrm{or}\, k_{\sigma}\le 0}\sigma(\pi_K)^{k_{\sigma}}\delta_1(\pi_K)$, 
 $k'_{\sigma}=k_{\sigma}$ if $a_{\sigma}\not= 0$, $k'_{\sigma}=-k_{\sigma}$ for other $\sigma$ and $\{\sigma|a'_{\sigma}\not= 0\}=\{\sigma|b'_{\sigma}\not= 0\}$. So we have finished the proof of (1).

(2) follows from the above theorem and (1) of Theorem $\ref{36}$.
 
\end{proof}

\subsection{Potentially semi-stable and non-cristalline split trianguline $E$-representations.}
In this subsection, we explicitly describe the $E$-filtered $(\varphi,N,G_K)$-modules 
of potentially semi-stable split trianguline $E$-representations which are not 
potentially cristalline. First we define the parameter spaces of potentially semi-stable split 
triangulilne $E$-representations. We put
\begin{itemize}
\item[$\bullet$]$T_{\mathrm{st}}:=\{(\delta,\{k_{\sigma}\}_{\sigma})| \delta:K^{\times}\rightarrow E^{\times}$ a locally constant character, $k_{\sigma}\in\mathbb{Z}$ for any $\sigma$ such that $ (\sum_{\sigma}k_{\sigma})+2e_K\mathrm{val}_p(\delta(\pi_K))+[K:\mathbb{Q}_p] =0\}$.
\end{itemize}
For any $(\delta,\{k_{\sigma}\}_{\sigma})\in T_{\mathrm{st}}$, we define the 
parameter space $T_{\mathrm{st}}(\delta,\{k_{\sigma}\}_{\sigma})$ of weakly admissible filtrations as follows. We define
\begin{itemize}
\item[$\bullet$]$T_{\mathrm{st}}(\delta,\{k_{\sigma}\}_{\sigma}):=\oplus_{\sigma,k_{\sigma}\ge 1}Ee_{\sigma}$.
\end{itemize}
For any $(a_{\sigma}e_{\sigma})_{\sigma}\in T_{\mathrm{st}}(\delta,\{k_{\sigma}\}_{\sigma})$ we define a rank two $E$-filtered $(\varphi,N,\mathrm{Gal}(K_{n(\delta)}/K))$-module 
$D_{(\delta,\{k_{\sigma}\}_{\sigma}),(a_{\sigma}e_{\sigma})_{\sigma}}$. However, as in the potentially cristallilne case, the construction depend on the choice of one more parameter. For any $(\delta,\{k_{\sigma}\}_{\sigma})\in T_{\mathrm{st}}$, we fix 
$\alpha\in (K_0\otimes_{\mathbb{Q}_p}E)^{\times}$ such that $g_3(\alpha)=\delta(\pi_K)$, where 
$g_3$ is as before. Then we define $D_{(\delta,\{k_{\sigma}\}_{\sigma}),(a_{\sigma}e_{\sigma})_{\sigma}}:=(K_0\otimes_{\mathbb{Q}_p}E)e_1
\oplus (K_0\otimes_{\mathbb{Q}_p}E)e_2$ as follows:
\begin{itemize}
\item[(1)]$N(e_2)=e_1, N(e_1)=0$.
\item[(2)]$\varphi(e_1)=\frac{\alpha}{p}e_1, \varphi(e_2)=\alpha e_2$ for the fixed $\alpha\in (K_0\otimes_{\mathbb{Q}_p}E)^{\times}$,
\item[]so we have $\varphi^f(e_1)=\frac{\delta(\pi_K)}{q}e_1, \varphi^f(e_2)=\delta(\pi_K)e_2$.
\item[(3)]$g(e_1)=\delta(\chi_{\mathrm{LT}}(g))e_1, g(e_2)=\delta(\chi_{\mathrm{LT}}(g))e_2$ for any $g\in G_K$.
\item[(4)]We put $K_{n(\delta)}\otimes_{K_0}D_{(\delta,\{k_{\sigma}\}_{\sigma}),(a_{\sigma}e_{\sigma})_{\sigma}}=(K_{n(\delta)}\otimes_{\mathbb{Q}_p}E)e_1\oplus (K_{n(\delta)}\otimes_{\mathbb{Q}_p}E)e_2\isom\oplus_{\sigma}(K_{n(\delta)}\otimes_{K,\sigma}E)e_{1,\sigma}\oplus (K_{n(\delta)}\otimes_{K,\sigma}E)e_{2,\sigma}=:\oplus_{\sigma}D_{\sigma}$, then
\begin{itemize}
\item[(i)] For $\sigma$ such that $k_{\sigma}\le -1$, we put $\mathrm{Fil}^0D_{\sigma}=D_{\sigma}, \mathrm{Fil}^1D_{\sigma}=\cdots =\mathrm{Fil}^{-k_{\sigma}}D_{\sigma}=(K_{n(\delta)}\otimes_{K,\sigma}E)e_{2,\sigma}, \mathrm{Fil}^{-k_{\sigma}+1}D_{\sigma}=0$,
\item[(ii)]For $\sigma$ such that $k_{\sigma}=0$, we put $\mathrm{Fil}^0D_{\sigma}=D_{\sigma}, \mathrm{Fil}^1D_{\sigma}=0$,
\item[(iii)]For $\sigma$ such that $k_{\sigma}\ge 1$, we put $\mathrm{Fil}^{-k_{\sigma}}D_{\sigma}=D_{\sigma}, 
\mathrm{Fil}^{-k_{\sigma}+1}D_{\sigma}=\cdots =\mathrm{Fil}^0D_{\sigma}=(K_{n(\delta)}\otimes_{K,\sigma}E)(a_{\sigma}e_{1,\sigma}+e_{2,\sigma}), \mathrm{Fil}^1D_{\sigma}=0$.
\end{itemize}
\end{itemize}
We will prove in the proof of the next theorem that these are weakly admissible. So, by ``weakly admissible implies admissible" theorem ($\cite[\mathrm{Theorem}B]{Be04}$), there exists 
 unique (up to isomorphism) two dimensional potentially semi-stable, not potentially cristalline $E$-representation $V_{(\delta,\{k_{\sigma}\}),(a_{\sigma}e_{\sigma})_{\sigma}}$ such that $D^{K_{n(\delta)}}_{\mathrm{st}}(V_{(\delta,\{k_{\sigma}\}_{\sigma}),(a_{\sigma}e_{\sigma})_{\sigma}})\isom D_{(\delta,\{k_{\sigma}\}_{\sigma}),(a_{\sigma}e_{\sigma})_{\sigma}}$.
 Then our main result on the classification of potentially semi-stable split 
trianguline $E$-representations is as follows.
 \begin{thm}$\label{45}$
 Let $V$ be a two dimensional $E$-representation.
 Then the following conditions are equivalent.
\begin{itemize}
\item[$\mathrm{(1)}$]$V$ is split trianguline, potentially semi-stable and not potentially cristalline.
\item[$\mathrm{(2)}$]There exist $(\delta,\{k_{\sigma}\}_{\sigma})\in T_{\mathrm{st}}$, $(a_{\sigma}e_{\sigma})_{\sigma}\in T_{\mathrm{st}}(\delta,\{k_{\sigma}\}_{\sigma})$ and $\{w_{\sigma}\}_{\sigma}$ with $w_{\sigma}\in\mathbb{Z}$ for any $\sigma$, such that $V\isom V_{(\delta,\{k_{\sigma}\}_{\sigma}),(a_{\sigma}e_{\sigma})_{\sigma}}(\prod_{\sigma}\sigma(\chi_{\mathrm{LT}}^{w_{\sigma}}))$.
\end{itemize}
Moreover we have an isomorphism $V_{(\delta,\{k_{\sigma}\}_{\sigma}),(a_{\sigma}e_{\sigma})_{\sigma}}\isom V_{(\delta',\{k'_{\sigma}\}_{\sigma}),(a'_{\sigma}e_{\sigma})_{\sigma}}(\prod_{\sigma}\sigma(\chi_{\mathrm{LT}}^{w'_{\sigma}}))$ if and only if $((\delta',\{k'_{\sigma}\}_{\sigma}),(a'_{\sigma}e_{\sigma})_{\sigma})=
((\delta,\{k_{\sigma}\}_{\sigma}),(a_{\sigma}e_{\sigma})_{\sigma})$ and $w'_{\sigma}=0$ for any $\sigma$.
\end{thm} 
\begin{proof}
First we prove that $(1)$ implies $(2)$.
Let $V$ be as in $(1)$. Then, by Proposition $\ref{41}$, we have $V\isom V(s)$ for some $s\in S_{\mathrm{st}}(
\delta_1,\delta_2)$ such that $\delta_1/\delta_2=|N_{K/\mathbb{Q}_p}(x)|\prod_{\sigma}\sigma(x)^{k_{\sigma}}$ for some $k_{\sigma}\in\mathbb{Z}$ for any $\sigma$. Twisting $V$ by some suitable $\prod_{\sigma}\sigma(\chi_{\mathrm{LT}})^{l_{\sigma}}$, we may assume that 
$(\delta_1,\delta_2)=(\prod_{\sigma}\sigma(x)^{k_{\sigma}}|N_{K/\mathbb{Q}_p}(x)|\delta, \delta)$ 
for some locally constant character $\delta:K^{\times}\rightarrow E^{\times}$ and for some 
$k_{\sigma}\in\mathbb{Z}$ for any $\sigma$.  By definition of $S_{\mathrm{st}}(\delta_1,\delta_2)\isom s(\{k_{\sigma}\}_{\sigma})+ D^K_{\mathrm{dR}}(W(\delta_1/\delta_2))/\mathrm{Fil}^0D^K_{\mathrm{dR}}(W(\delta_1/\delta_2))$, $s$ corresponds to $s(\{k_{\sigma}\}_{\sigma})+\overline{b}$ 
for some $\overline{b}\in D^K_{\mathrm{dR}}(W(\delta_1/\delta_2))/\mathrm{Fil}^0D^K_{\mathrm{dR}}(W(\delta_1/\delta_2))$. So $W(s)$ is given as follows:
\begin{itemize}
\item[(1)]$W_e(s):=W_e(\delta_1)\oplus W_e(\delta_2)$ such  that $g(x,y)=(gx+c(g)gy,gy)$ for any 
$g\in G_K, x\in W^+_{\mathrm{cris}}(\delta_1), y\in 
W^+_{\mathrm{cris}}(\delta_2)$. 
(Here $c(g)$ is defined in the definition of $s(\{k_{\sigma}\}_{\sigma})$ and we use the fact that $W_e(\delta_1)\isom W_e(\delta_2)(\chi)$ which we can see from Lemma $\ref{27}$.)
\item[(2)]$W^+_{\mathrm{dR}}(s):=W^+_{\mathrm{dR}}(\delta_1)\oplus W^+_{\mathrm{dR}}(\delta_2)$ such that $g(x,y)=(gx,gy)$ for any $g\in G_K, x\in W^+_{\mathrm{dR}}(\delta_1), y\in 
W^+_{\mathrm{dR}}(\delta_2)$.
\item[(3)]$\iota:B_{\mathrm{dR}}\otimes_{B_e}W_e(s)\isom B_{\mathrm{dR}}\otimes_{B^+_{\mathrm{dR}}}W^+_{\mathrm{dR}}(s)$ is an isomorphism given by $\iota(x,y)=(x+\frac{\mathrm{log}([\tilde{p}])}{t}y+b\otimes y, y)$, here $b=(b_{\sigma})\in D^K_{\mathrm{dR}}(W(\delta_1/\delta_2))$ is the lift of $\overline{b}$ such that $b_{\sigma}=0$ for any $\sigma$ such that $k_{\sigma}\le 0$.
\end{itemize}
We calculate $D^{K_{n(\delta)}}_{\mathrm{st}}(V(s))$ as follows. First, if we take a base $e'\in D^{K_{n(\delta)}}_{\mathrm{cris}}(W(\delta_2))\subseteq B_{\mathrm{cris}}\otimes_{B_e}W_e(\delta_2)$, then we have $e'':=\frac{1}{t}e'\in D^{K_{n(\delta)}}_{\mathrm{cris}}(W(\delta_1))$ because $W_e(\delta_1)\isom W_e(\delta_2)(\chi)$. Then $\varphi^f(e')=\delta(\pi_K)e'$, $\varphi^f(e'')=\frac{\delta(\pi_K)}{q}e''$ and $g(e')=\delta(\chi_{\mathrm{LT}}(g))e', g(e'')=\delta(\chi_{\mathrm{LT}}(g))e''$ for any $g\in G_K$. Moreover 
we can take a base $e'$ such that $\varphi(e')=\alpha e'$ for the fixed 
$\alpha$. Then $D^{K_{n(\delta)}}_{\mathrm{st}}(W(s))\isom 
(K_0\otimes_{\mathbb{Q}_p}E)e''\oplus (K_0\otimes_{\mathbb{Q}_p}E)(e'-\mathrm{log}([\tilde{p}])e'')$. If we put $e_1:=e''$ , $e_2:=e'-\mathrm{log}([\tilde{p}])e''$, then $N(e_2)=e_1$. Next we calculate the filtration on $K_{n(\delta)}\otimes_{K_0}D^{K_{n(\delta)}}_{\mathrm{st}}(W(s))=(K_{n(\delta)}\otimes_{\mathbb{Q}_p}E)e_1\oplus(K_{n(\delta)}\otimes_{\mathbb{Q}_p}E)e_2$. By definition of $\iota$ in $(3)$ of the definition of $W(s)$, we have $xe_1+ye_2\in \mathrm{Fil}^i$ for $x,y\in K_{n(\delta)}\otimes_{\mathbb{Q}_p}E$ if and only if $xe''+y(\iota(e')-\mathrm{log}([\tilde{p}])e'')=xe''+y(e'+b\otimes e''+\mathrm{log}([\tilde{p}])e''-\mathrm{log}([\tilde{p}])e'')= (xe''+yb\otimes e'')+ye'\in \mathrm{Fil}^iD^{K_{n(\delta)}}_{\mathrm{dR}}(W(\delta_1))\oplus \mathrm{Fil}^iD^{K_{n(\delta)}}_{\mathrm{dR}}(W(\delta_2))$. So we can calculate the filtration similarly as in the proof of Theorem 
$\ref{42}$. Then we can show that $D^{K_{n(\delta)}}_{\mathrm{st}}(V(s))\isom D_{(\delta,\{k_{\sigma}\}_{\sigma}),(a_{\sigma}e_{\sigma})_{\sigma}}$ for some $(a_{\sigma}e_{\sigma})_{\sigma}\in T_{\mathrm{st}}(\delta,\{k_{\sigma}\}_{\sigma})\isom\oplus_{\sigma, k_{\sigma}\ge 1}Ee_{\sigma}$ such that $\{\sigma|a_{\sigma}\not= 0\}=
\{\sigma|b_{\sigma}\not= 0\}$. So $D_{(\delta,\{k_{\sigma}\}_{\sigma}),(a_{\sigma}e_{\sigma})_{\sigma}}$ is 
weakly admissible for such a $(a_{\sigma}e_{\sigma})_{\sigma}$ and we have $V(s)\isom V_{(\delta,\{k_{\sigma}\}_{\sigma}),(a_{\sigma}e_{\sigma})_{\sigma}}$.
So we have proved that $(1)$ implies $(2)$. 

Next we prove that (2) implies (1).
Let us take $(\delta,\{k_{\sigma}\}_{\sigma})\in T_{\mathrm{st}}$ and 
$(a_{\sigma}e_{\sigma})_{\sigma}\in T_{\mathrm{st}}(\delta,\{k_{\sigma}\}_{\sigma})$. Then we can see that $(\delta_1,\delta_2):=(\prod_{\sigma}\sigma(x)^{k_{\sigma}}|N_{K/\mathbb{Q}_p}(x)|\delta, \delta)\in S^+$.
Moreover it is easy to see from the above argument that there exist some $s:=(b_{\sigma}e_{\sigma})_{\sigma}\in S_{\mathrm{st}}(\delta_1,\delta_2)$ such that 
$D_{(\delta,\{k_{\sigma}\}_{\sigma}),(a_{\sigma}e_{\sigma})_{\sigma}}\isom D^{K_{n(\delta)}}_{\mathrm{st}}(V(s))$. So $D_{(\delta,\{k_{\sigma}\}_{\sigma}),(a_{\sigma}e_{\sigma})_{\sigma}}$ is weakly admissible and $V(s)\isom V_{(\delta,\{k_{\sigma}\}_{\sigma}),(a_{\sigma}e_{\sigma})_{\sigma}}$.

Finally we prove the uniqueness of $((\delta,\{k_{\sigma}\}_{\sigma}),(a_{\sigma}e_{\sigma})_{\sigma},\{w_{\sigma}\}_{\sigma})$. This easily follows from the above 
argument and from Theorem $\ref{36}$.
We have finished the proof of the theorem.

\end{proof}

\appendix

\section{A relation between two dimensional potentially semi-stable 
trianguline $E$-representations and classical local Langlands 
correspondence for $\mathrm{GL}_2(K)$.}

In this appendix, we show a simple relation between two dimensional potentially semi-stable 
trianguline $E$-representations and classical local Langlands 
correspondence for $\mathrm{GL}_2(K)$.

\subsection{Classical local Langlands correspondence for $\mathrm{GL}_2(K)$.}
First we briefly recall classical local Langlands correspondence for $\mathrm{GL}_2(K)$.
Let $W_K\subset G_K$ be the Weil group of $K$, i.e. the inverse image of $<\mathrm{Frob}_k>_{\mathbb{Z}}\subseteq \mathrm{Gal}(\bar{k}/k)$ by the natural surjection $G_K\rightarrow \mathrm{Gal}(\bar{k}/k)$. (Here $<\mathrm{Frob}_k>_{\mathbb{Z}}$ is the subgroup of $\mathrm{Gal}(\bar{k}/k)$ generated by the $q=p^f$-th power 
Frobenius $\mathrm{Frob}_k$ of $k$.) $W_K$ is a topological group such that the inertia $I_K:=\mathrm{Ker}(G_K\rightarrow \mathrm{Gal}(\bar{k}/k))$ equipped with usual profinite  topology is an open subgroup of 
$W_K$. 
\begin{defn}
Let $L$ be a field of characteristic zero.
We say that a finite dimensional $L$-vector space $D$ with discrete topology is an $L$-Weil-Deligne representation of $K$ if $D$ is equipped with
\begin{itemize}
\item[(1)] a continuous $L$-linear action of $W_K$, i.e., there is a continuous 
morphism $\rho:W_K\rightarrow \mathrm{Aut}_L(D)$,
\item[(2)] a nilpotent $L$-linear operator $N:D\rightarrow D$ such 
that $N\rho(g)=q^{-{\mathrm{deg}(g)}}\rho(g)N$ for any $g\in W_K$.
Here we define $\mathrm{deg}(g)\in\mathbb{Z}$ such that $g=\mathrm{Frob_k}^{\mathrm{deg}(g)}\in\mathrm{Gal}(\bar{k}/k)$ for any $g\in W_K$.
\end{itemize}
\end{defn}
We say that an $L$-Weil-Deligne representation $D:=(D,\rho,N)$ is semi-simple if 
$(D,\rho)$ is a semi-simple representation of $W_K$.

Next we recall Fontaine's recipe to construct a $\bar{K}$-Weil-Deligne representation
 from a potentially semi-stable $E$-representation of $G_K$.
Let $V$ be a potentially semi-stable $d$ dimensional $E$-representation of $G_K$.Then $D_{\mathrm{pst}}(V):=\cup_{K\subset L,\mathrm{finite}}(B_{\mathrm{log}}\otimes
_{\mathbb{Q}_p}V)^{G_L}$ is a free $K_0^{\mathrm{un}}\otimes_{\mathbb{Q}_p}E$-
module of rank $d$ equipped with the natural $(\varphi, N, G_K)$-action on which $\varphi$ and $G_K$ act semi-linerly and $N$ acts linearly.
We define a continuous action of $W_K$ on $D_{\mathrm{pst}}(V)$ by $\rho(g)x:=
\varphi^{-f\mathrm{deg}(g)}gx$ for any $g\in W_K$ and $x\in D_{\mathrm{pst}}(V)$.
We can see that this action is $K_0^{\mathrm{un}}\otimes_{\mathbb{Q}_p}E$-linear 
and satsfies $N\rho(g)=q^{-\mathrm{deg}(g)}\rho(g)N$ for any $g\in W_K$.
By using the fixed embeddings $E\hookrightarrow \bar{K}$ and $K_0^{\mathrm{un}}\hookrightarrow \bar{K}$, we define a map $K^{\mathrm{un}}\otimes_{\mathbb{Q}_p}E\rightarrow \bar{K}$. Extending scalar by this map, we get a $\bar{K}$-Weil-Deligne representation $\bar{D}_{\mathrm{pst}}(V):=D_{\mathrm{pst}}(V)\otimes
_{K^{\mathrm{un}}\otimes_{\mathbb{Q}_p}E}\bar{K}$ of $K$ of dimension $d$.
So we get a functor $V\mapsto \bar{D}_{\mathrm{pst}}(V)$ from the category of potentially semi-stable $E$-representations
of $G_K$ to the category of $\bar{K}$-Weil-Deligne representations of $K$.

Next, for a $\bar{K}$-Weil-Deligne representation $\bar{D}$ of $K$, we recall the definition of $L$-factor of  $\bar{D}$. Let us fix an isomorphism 
$\bar{K}\isom \mathbb{C}$. By this isomorphism, we can see $\bar{D}$ as a $\mathbb{C}$-Weil-Deligne representation of $K$. Put $\bar{J}:=\bar{D}^{N=0, I_K=1}$, the subspace of $\bar{D}$ on which $N=0$ and $I_K$ acts trivially.
Let $\sigma\in W_K$ be an element such that $\mathrm{deg}(\sigma)=-1$.
Then $\sigma$ acts on $\bar{J}$ and this action does not depend on the 
choice of a lifting $\sigma$.
Then we define the $L$-factor of $\bar{D}$ by $L(\bar{D}, s):=\mathrm{det}(1-\sigma q^{-s}|_{\bar{J}})^{-1}$. 

For a two dimensional semi-simple $\bar{K}$-Weil-Deligne representation 
$\bar{D}$, we can define an irreducible smooth admissible representation $\pi(\bar{D})$ of $\mathrm{GL}_2(K)$ as in $\cite[33.1]{Bu-He}$ (if we fix a nontrivial additive smooth character $\psi:K\rightarrow \mathbb{C}^{\times}$). The irreducible smooth admissible representations of $\mathrm{GL}_2(K)$ are classified into non supercuspidal ones (i.e. one dimensional 
representations or principal series representations or special series representations) and supercuspidal ones. We do not recall these definitions here.
We only need the following proposition concerning this correspondence.
\begin{prop}$\label{17}$
Let $\bar{D}$ be a $\bar{K}$-Weil-Deligne representation of $K$. 
Let $\bar{D}^{ss}$ be the semi-simplification of $\bar{D}$.
Then the following conditions are equivalent.
\begin{itemize}
\item[$\mathrm{(1)}$]$\pi(\bar{D}^{ss})$ is non supercuspidal.
\item[$\mathrm{(2)}$]The representation $(\bar{D},\rho)$ of $W_K$ is reducible.
\item[$\mathrm{(3)}$]There is a continuous character $\delta:W_K\rightarrow 
\bar{K}^{\times}$ such that $L(\bar{D}\otimes_{\bar{K}}\bar{K}(\delta), s)
\not= 0$, here $\bar{K}$ is equipped with discrete topology.

\end{itemize}
\end{prop}
\begin{proof}
This follows from $\cite[\mathrm{Proposition}\, 33.2]{Bu-He}$.
\end{proof}

\subsection{Two dimensional trianguline representations and non supercuspidal representations.}
Let $V$ be a two dimensional potentially semi-stable $E$-representation of $G_K$. Then $\pi(\bar{D}_{\mathrm{pst}}(V)^{ss})$ is an irreducible smooth 
admissible representation of $\mathrm{GL}_2(K)$.
In this subsection, we prove a relation between two dimensional potentially 
semi-stable trianguline $E$-representations and non supercuspidal representations of $\mathrm{GL}_2(K)$.

Before stating the main result of this appendix, we gives another useful 
characterization of two dimensional split trianguline $E$-representations in terms of 
$D^K_{\mathrm{cris}}$. This characterization is also a generalization 
of $\cite[\mathrm{Proposition}\, 5.3]{Co07a}$.
\begin{prop}$\label{18}$
Let $V$ be a two dimensional $E$-representation of $G_K$.
Then the following conditions are equivalent:
\begin{itemize}
\item[$\mathrm{(1)}$]$V$ is a split trianguline $E$-representation.
\item[$\mathrm{(2)}$]There exists a continuous character $\delta:G_K\rightarrow E^{\times}$ and an element $\alpha\in (K_0\otimes_{\mathbb{Q}_p}E)^{\times}$ such that $D^K_{\mathrm{cris}}(V(\delta))^{\varphi=\alpha}\not= 1$.
\end{itemize}

\end{prop}
\begin{proof}
First we prove that (1) implies (2).
Assume that $V$ is a split trianguline $E$-representation.
By definition, $W(V)$ sits in a following short exact sequence of 
$E$-$B$-pairs
\begin{equation*}
0\rightarrow W_1\rightarrow W(V)\rightarrow W_2\rightarrow 0.
\end{equation*}
Here $W_i$ are rank one $E$-$B$-pairs for $i=1,2$. By Theorem $\ref{15}$ and 
by the construction $W(\delta)$ for any $\delta:K^{\times}\rightarrow E^{\times}$, if we fix a uniformizer $\pi_E$ of E then there exist $(\delta_i, k_i)$ where $\delta_i:G_K\rightarrow E^{\times}$ are continuous characters and $k_i\in\mathbb{Z}$ for $i=1,2$ such that 
$W_i\isom W(E(\delta_i))\otimes W_0^{\otimes k_i}$, here $W_0$ is defined before Theorem $\ref{15}$.
If we twist the above exact sequence by the character $\delta_1^{-1}$ and apply the  left exact functor $D^K_{\mathrm{cris}}$, we get an inclusion $D_0^{\otimes k_1}\subseteq D^K_{\mathrm{cris}}(V(\delta_1^{-1}))$ of $E$-$\varphi$-modules of $K$. If we put $K_0\otimes_{\mathbb{Q}_p}E\isom \oplus_{0\le i\le f-1}K_0\otimes
_{E_0,\varphi^i}Ee_{i}$,
$D_0^{\otimes k_1}\isom \oplus_{0\le i\le f-1}K_0\otimes_{E_0,\varphi^i}Ee_{i,k_1}$, then we have $\varphi(e_{0,k_1}+\cdots e_{f-1,k_1})=(\pi_E^{k_1}e_0+e_1+\cdots
e_{f-1})(e_{0,k_1}+\cdots + e_{f-1,k_1})$. In particular, $(e_{0,k_1}+\cdots +
e_{f-1, k_1})\in D^K_{\mathrm{cris}}(V(\delta^{-1}_1))^{\varphi=(\pi_E^{k_1}e_0+e_1+\cdots +e_{f-1})}$ and $(\pi_E^{k_1}e_0+e_1+\cdots +e_{f-1})\in (K_0\otimes_{\mathbb{Q}_p}E)^{\times}$. So we have proved that (1) implies (2).
Next we prove that (2) implies (1). 
Assume that $D^K_{\mathrm{cris}}(V(\delta))^{\varphi=\alpha}\not= 0$ for 
a character $\delta:G_K\rightarrow E^{\times}$ and $\alpha\in (K_0\otimes_{\mathbb{Q}_p}E)^{\times}$. Because $V$ is split trianguline if and only if $V(\delta)$ is spllit trianguline, we may assume that $V=V(\delta)$.
Let us take a nonzero $x\in D^K_{\mathrm{cris}}(V)^{\varphi=\alpha}$. We 
consider the sub $E$-filtered $\varphi$-module $D_1\subseteq D^K_{\mathrm{cris}}(V)$ of rank one which is generated by $x$. Then by Theorem $\ref{4}$ (3), we have a natural inclusion of $E$-$B$-pairs $W(D_1)\hookrightarrow W(D^K_{\mathrm{cris}}(V))\hookrightarrow W(V)$. $W(D_1)$ is an $E$-$B$-pair of rank one by Theorem $\ref{4}$ (3). Taking the saturation $W(D_1)^{\mathrm{sat}}$ of $W(D_1)$ in $W(V)$ (Lemma $\ref{3}$), we get a short exact sequence of $E$-$B$-pairs
\begin{equation*}
0\rightarrow W(D_1)^{\mathrm{sat}}\rightarrow W(V)\rightarrow W_2\rightarrow 0.
\end{equation*}
Here $W_2$ is the cokernel of $W(D_1)^{\mathrm{sat}}\hookrightarrow 
W(D)$, which is an $E$-$B$-pair of rank one. So $V$ is a split trianguline $E$-representation. We have finished the proof of this proposition.
\end{proof}

The main theorem of this appendix is as follows.
\begin{thm}$\label{19}$
Let $V$ be a two dimensional potentially semi-stable $E$-representation of $G_K$.
Then the following conditions are equivalent:
\begin{itemize}
\item[$\mathrm{(1)}$] $V$ is trianguline, i.e. $V\otimes_{E}E'$ is a split 
trianguline $E'$-representation
 for some finite extension of $E'$ of $E$.
\item[$\mathrm{(2)}$] The representation $(\bar{D}_{\mathrm{pst}}(V),\rho)$ of 
$W_K$ is reducible.
\item[$\mathrm{(2)}$]$\pi(\bar{D}_{\mathrm{pst}}(V)^{ss})$ is a non supercuspidal representation of $\mathrm{GL}_2(K)$.
\end{itemize}
\end{thm}
\begin{proof}
By Proposition $\ref{17}$, it suffices to show equivalence between 
(1) and (3). The proof of this is a modified version of the proof 
of $\cite[\mathrm{Lemma}\,1.3]{Ki}$.
First we prove that (1) implies (3). Assume that $V\otimes_{E}E'$ is a split 
trianguline 
$E'$-representation
 for a finite extension $E'$ of $E$. By construction, we have $\pi(\bar{D}_{\mathrm{pst}}(V\otimes_{E}E')^{ss})=\pi(\bar{D}_{\mathrm{pst}}(V)^{ss})$. So we may assume that $E'=E$ and $V$ is a spllit trianguline $E$-representation.
Then, by Proposition $\ref{18}$, there exist $\alpha\in (K_0\otimes_{\mathbb{Q}_p}E)^{\times}$ and $\delta:G_K
\rightarrow E^{\times}$ such that $D^K_{\mathrm{cris}}(V(\delta))^{\varphi=\alpha}\not= 0$.
In this case, we claim that we can take $\delta$ which is a potentially cristalline character. We show this claim 
as follows. In the proof of Proposition $\ref{18}$, $W(V)$ sits in a following short exact
 sequence of $E$-$B$-pairs
 \begin{equation*}
 0\rightarrow W_1\rightarrow W(V)\rightarrow W_2\rightarrow 0,
 \end{equation*} 
 where $W_1:=W_0^{\otimes k_1}\otimes W(E(\delta'))$  and $W_2$ are rank one $E$-$B$-pairs.
 In this case, because $W(V)$ is a potentially semi-stable $E$-$B$-pair, $W_1$ is also a potentially semi-stable $E$-$B$-pair.
Because $W_0$ is a cristalline $E$-$B$-pair by definition, so $W(E(\delta'))\isom W_1\otimes W_0^{\otimes -k_1}$ is also potentially semi-stable. 
So $\delta'$ is a potentially cristalline character because rank one semi-stable $E$-representations are cristalline. 
 Then, as in the proof of Proposition $\ref{18}$, if we put $\delta:=\delta^{' -1}$, then 
 we have $D^K_{\mathrm{cris}}(V(\delta))^{\varphi=\alpha}\not= 0$ for some $\alpha\in (K_0\otimes
 _{\mathbb{Q}_p}E)^{\times}$. So we have proved the claim. For this $\delta$, we have 
$ \bar{D}_{\mathrm{pst}}(V(\delta))\isom \bar{D}_{\mathrm{pst}}(V)\otimes_{\bar{K}}
 \bar{D}_{\mathrm{pst}}(E(\delta))$ because both $V$ and $E(\delta)$ are potensially 
 semi-stable. So, by Proposition $\ref{17}$, it suffices to show that $L(\bar{D}_{\mathrm{pst}}
 (V(\delta)), s)\not= 1$. Put $\bar{J}:=\bar{D}_{\mathrm{pst}}(V(\delta))^{N=0, I_K=1}$. 
 If we put $J:=D_{\mathrm{pst}}(V(\delta))^{N=0, I_K=1}$, then we have $J\isom
 \cup_{K\subset L\subset K^{\mathrm{un}}}D^L_{\mathrm{cris}}(V(\delta))$, here $K^{\mathrm{un}}$ is the maximal unramified extension of $K$.
By using the fact that $N$ is nilpotent on $D_{\mathrm{pst}}(V(\delta))$ and 
that $I_K$ acts discretely on it, we can easily see that $\bar{J}\isom J\otimes_{K_0^{\mathrm{un}}
 \otimes_{\mathbb{Q}_p}E}\bar{K}$. 
 Moreover, because $\mathrm{Gal}(K^{\mathrm{un}}/K)\isom\mathrm{Gal}(K_0^{\mathrm{un}}/K_0)$ 
 acts on $J$ discretely and semi-linearly, 
 so we have $J\isom J^{\mathrm{Gal}(K^{\mathrm{un}}/K)}\otimes_{K_0\otimes_{\mathbb{Q}_p}E}(K_0^{\mathrm{un}}\otimes_{\mathbb{Q}_p}E)$ by Hilbert's theorem 90. Because we have $J^{\mathrm{Gal}
 (K^{\mathrm{un}}/K)}=D^K_{\mathrm{cris}}(V(\delta))$,  we have $\bar{J}\isom D_{\mathrm{cris}}
 (V(\delta))\otimes_{K_0\otimes_{\mathbb{Q}_p}E}\bar{K}$. Finally, if we take a non zero element 
 $x\in D_{\mathrm{cris}}(V(\delta))^{\varphi=\alpha}$ and consider the action of $\sigma$ on $x$
  as Weil-Deligne representation (here $\sigma$ is an element in $W_K$ such that $\mathrm{deg}(\sigma)=-1$), then we have $\rho(\sigma)(x)=\varphi^{f}\sigma(x)=\varphi^{f}(x)=\varphi^{f-1}(\alpha)\cdots \varphi(\alpha)
  \alpha x =\beta x$, here we put $\beta:=\varphi^{f-1}(\alpha)\cdots
  \varphi(\alpha)\alpha\in ((K_0\otimes_{\mathbb{Q}_p}E)^{\times})^{\varphi=1}=E^{\times}$.
  So $\bar{J}$ has a non zero eigenvector of $\sigma$, hence we get $L(\bar{D}_{\mathrm{pst}}(V(\delta)), s)\not= 1$. So we have proved that (1) implies  (3) 
  by Proposition $\ref{17}$.
  
  Next we prove that (3) implies (1).
  Let $V$ be a two dimensional potentially semi-stable $E$-representation such that $\pi(\bar{D}_{\mathrm{pst}}(V)^{ss})$ is non supercuspidal. 
  Then, by Proposition $\ref{17}$, $L(\bar{D}_{\mathrm{pst}}(V)\otimes_{\bar{K}}\bar{K}(\delta), s)\not= 1$ for a character $\delta:W_K\rightarrow \bar{K}^{\times}$. 
  By this we can take a non zero eigenvector $x\in(\bar{D}_{\mathrm{pst}}(V)\otimes_{\bar{K}}\bar{K}(\delta))^{N=0, I_K=1}$ of $\sigma\in W_K$ with a non zero eigenvalue $\beta\in \bar{K}^{\times}$. If we take $E$ large enough, 
  we may assume that $\beta\in E^{\times}$ and that there exsists a potentially cristalline character $\delta':G_K\rightarrow E^{\times}$ such that $\bar{D}_{\mathrm{pst}}(E(\delta'))\isom \bar{K}(\delta)$.
   Then we can prove that there is an isomorphism $(\bar{D}_{\mathrm{pst}}(V)\otimes_{\bar{K}}\bar{K}(\delta))^{N=0, I_K=1}\isom (\bar{D}_{\mathrm{pst}}(V(\delta')))^{N=0, I_K=1}\isom 
   D^K_{\mathrm{cris}}(V(\delta'))\otimes_{K_0\otimes_{\mathbb{Q}_p}E} \bar{K}$ 
   in the same way as the argument in the proof of (1) $\Rightarrow$ (3).  If we decompose $D^K_{\mathrm{cris}}(V(\delta'))
   \isom \oplus_{\tau:K_0\hookrightarrow E}D_{\tau}e_{\tau}$, then we have $D^K_{\mathrm{cris}}(V(\delta'))\otimes_{K_0\otimes_{\mathbb{Q}_p}E}\bar{K}=D_{id}\otimes_E \bar{K}$. So, if we take $E$ large enough, we may assume that $x$ is contained in $D_{id}e_{id}$. If we put $e:=x+\varphi(x)+\cdots +\varphi^{f-1}(x)\in \oplus_{\tau:K_0\hookrightarrow E}D_{\tau}e_{\tau}=D_{\mathrm{cris}}^K(V(\delta'))$, then we have $\varphi^i(x)\in
   D_{\varphi^{-i}}e_{\varphi^{-i}}$ and $\varphi(e)=(\beta e'_{id}+e'_{\varphi^{-1}}+\cdots+e'_{\varphi^{-(f-1)}})e$, where $(\beta e'_{id}+\cdots+e'_{\varphi^{-(f-1)}})\in (K_0\otimes_{\mathbb{Q}_p}E)^{\times}
   =(\oplus_{0\le i\le f-1}Ee'_{\varphi^{-i}})^{\times}$. So, by Proposition 
   $\ref{18}$, $V(\delta')$ is a split trianguline $E$-representation. So 
   $V$ is also a split trianguline $E$-representation.
   We have finished the proof of the theorem.
   \end{proof}

\section{List of notations}
Here is a list of the main notations of the article, in the order in which 
they appear.

\begin{itemize}
\item[0.1:]$p$, $K$, $G_K$, $\varphi$, $K_{\infty}$, $\Gamma_K$.
\item[0.2:]$E$, $k$, $K_0$, $K_0'$, $\chi$, $\mathbb{C}_p$, 
$\chi_{\mathrm{LT}}$, $\mathrm{rec}_K$.
\item[1.1:] $\widetilde{\mathbb{E}}^+$, $\widetilde{\mathbb{A}}^+$, 
$\widetilde{\mathbb{A}}$, $\widetilde{\mathbb{B}}^+$, $\widetilde{\mathbb{B}}$, $B^+_{\mathrm{dR}}$, $A_{\mathrm{max}}$, $[\tilde{p}]$, $B^+_{\mathrm{max}}$, $t$, $B_{\mathrm{max}}$, $B_{\mathrm{dR}}$, $B_e$, $\widetilde{B}^+_{\mathrm{rig}}$, $\mathrm{Fil}^iB_{\mathrm{dR}}$, $\mathrm{log}[\tilde{p}]$, $B_{\mathrm{log}}$, $N$, $V$, $W$, $W_e$, $W^+_{\mathrm{dR}}$, $W_{\mathrm{dR}}$, 
$W(V)$, $\widetilde{B}^{\dagger}_{\mathrm{rig}}$, 
$\mathrm{rank}(W)$, $W_1\otimes W_2$, $W^{\vee}$, 
$W_1^{\mathrm{sat}}$, 
$D^L_{\mathrm{cris}}(W)$, $D^L_{\mathrm{st}}(W)$, $D^L_{\mathrm{dR}}(W)$, $W(D)$.

\item[1.2:]$B^{\dagger}_{\mathrm{rig},K}$, 
$\mathcal{R}_L$, $\phi_L$, $\mathcal{R}_L^{\mathrm{bd}}$, 
$\mathcal{R}_L^{\mathrm{int}}$, $v_L$, $\omega$, $\mathrm{deg}(M)$, 
$\mu(M)$, $[a]_* M$, $\widetilde{A}_{[r,s]}$, $\widetilde{B}_{[r,s]}$, 
$\widetilde{B}^{\dagger,r}$, $\widetilde{B}^{\dagger}$, $\widetilde{B}
^{\dagger,r}_{\mathrm{rig}}$, $\widetilde{B}
^{\dagger}_{\mathrm{rig}}$, $i_n$, $A_{K_0}$, $E_{K_0}$, 
$\mathbb{A}$, $\mathbb{B}$, $A_K$, $B_K$, $B_K^{\dagger}$, 
$B_K^{\dagger,r}$, $E_K$,  
$B^{\dagger,r}_{\mathrm{rig},K}$, $B^{\dagger}_{\mathrm{rig},K}$, 
\item[1.3:]$W_e(D)$, $W^+_{\mathrm{dR}}(D)$, $W(D)$, $D(W)$.
\item[1.4:]$D_0$, $W_0$, $\delta$, 
$W(\delta)$, $U_{\mathrm{fini}}^{H_K}$, $\nabla_{U}$, $D_{\mathrm{Sen}}(W)$, 
$\Theta_{\mathrm{Sen},W}$, $\{w_{\sigma}\}_{\sigma}$.
\item[2:]$B_E$.
\item[2.1:]$C^{\bullet}(W)$, $C^0(W)$, $C^1(W)$, $\mathrm{H}^*(G_K, W)$, 
$\mathrm{Ext}^1(W_2, W_1)$, $\mathrm{H}^1_e(G_K, W)$, $\mathrm{H}^1_f(G_K, W)$, $\mathrm{H}^1_g(G_K, W)$.
\item[2.3:]$\sigma(x)$, $N_{K/\mathbb{Q}_p}(x)$, $|-|$, 
$|N_{K/\mathbb{Q}_p}(x)|$, $\partial$.

\item[3.1:]$S^+$, $S(\delta_1,\delta_2)$, $\mathbb{P}_E(M)$, 
$S'(\delta_1,\delta_2)$, $S^{'\acute{\mathrm{e}}t}(\delta_1,\delta_2)$, 
$S^{' non-\acute{\mathrm{e}}t}(\delta_1,\delta_2)$, $W(s)$, $V(s)$.
\item[3.2:]$S^+_0$, $S^+_*$, $S^{'\mathrm{ord}}(\delta_1,\delta_2)$.
\item[4.1:]$S''(\delta_1,\delta_2)$, $S_{\mathrm{st}}(\delta_1,\delta_2)$, 
$S_{\mathrm{cris}}^{\acute{\mathrm{e}}t}(\delta_1,\delta_2)$.
\item[4.2:]$T_{\mathrm{cris}}$, $G(\delta)_{\sigma}$, 
$T_{\mathrm{cris}}(\delta_1,\delta_2,\{k_{\sigma}\}_{\sigma})$, 
$T^{\acute{\mathrm{e}}t}_{\mathrm{cris}}(\delta_1,\delta_2,\{k_{\sigma}\}_{\sigma})$, $D_{(\delta_1,\delta_2,\{k_{\sigma}\}_{\sigma}), x}$, $V_{(\delta_1,\delta_2,\{k_{\sigma}\}_{\sigma}), x}$.
\item[4.3:]$T_{\mathrm{st}}$, $T_{\mathrm{st}}(\delta, \{k_{\sigma}\}_{\sigma})$, $D_{(\delta,\{k_{\sigma}\}_{\sigma}), (a_{\sigma}e_{\sigma})_{\sigma}}$, 
$V_{(\delta,\{k_{\sigma}\}_{\sigma}), (a_{\sigma}e_{\sigma})_{\sigma}}$.
\item[A.1:]$W_K$, $I_K$, $D_{\mathrm{pst}}(V)$, 
$\bar{D}_{\mathrm{pst}}(V)$, $L(\bar{D},s)$, $\pi(\bar{D})$.
\end{itemize}

\end{document}